\documentclass[11pt]{amsart}

\usepackage{amsfonts, amstext, amsmath, amsthm, amscd, amssymb}
\usepackage{epsfig, graphics, psfrag}
\usepackage{color}

\renewcommand{\setminus}{{\smallsetminus}}

\newcommand{\HH}{{\mathbb{H}}}
\newcommand{\RR}{{\mathbb{R}}}
\newcommand{\ZZ}{{\mathbb{Z}}}

\newcommand{\CC}{{\mathbb{C}}}

\newcommand{\vol}{{\rm vol}}
\newcommand{\area}{{\rm area}}

\newcommand{\bdy}{{\partial}}
\newcommand{\abs}[1]{{\left\vert #1 \right\vert}}
\newcommand{\GA}{{\mathbb{G}_A}}
\newcommand{\GB}{{\mathbb{G}_B}}
\newcommand{\GRA}{{\mathbb{G}'_A}}
\newcommand{\GRB}{{\mathbb{G}'_B}}
\newcommand{\farey}{{\mathcal{F}}}
\newcommand{\tgen}{{t_{\mathrm{gen}}(D)}}

\newcommand{\p}{{\bf p}}
\newcommand{\q}{{\bf q}}
\newcommand{\ep}{\varepsilon_{\p}}
\newcommand{\eq}{\varepsilon_{\q}}

\def\co{\colon\thinspace}

\theoremstyle{plain}
\newtheorem{theorem}{Theorem}[section]
\newtheorem{corollary}[theorem]{Corollary}
\newtheorem{lemma}[theorem]{Lemma}
\newtheorem{prop}[theorem]{Proposition}

\newtheorem{question}[theorem]{Question}
\newtheorem*{no-num-theorem}{Theorem}
\newtheorem{conjecture}[theorem]{Conjecture}

\newtheorem*{namedtheorem}{\theoremname}
\newcommand{\theoremname}{testing}
\newenvironment{named}[1]{\renewcommand{\theoremname}{#1}\begin{namedtheorem}}{\end{namedtheorem}}

\theoremstyle{definition}
\newtheorem{define}[theorem]{Definition}
\newtheorem*{remark}{Remark}

\begin{document}
\title[Cusp areas of Farey manifolds]{Cusp areas of Farey manifolds \\
	and applications to knot theory}

\author[D. Futer]{David Futer}
\author[E. Kalfagianni]{Efstratia Kalfagianni}
\author[J. Purcell]{Jessica S. Purcell}

\address[]{Department of Mathematics, Temple University, Philadelphia, PA 19122, USA}

\email[]{dfuter@temple.edu}

\address[]{Department of Mathematics, Michigan State University, East
Lansing, MI 48824, USA}

\email[]{kalfagia@math.msu.edu}

\address[]{ Department of Mathematics, Brigham Young University,
Provo, UT 84602, USA}

\email[]{jpurcell@math.byu.edu }

\thanks{{Kalfagianni is supported in part by NSF--FRG grant DMS-0456155 and by NSF
grant DMS--0805942. Purcell is supported in part by NSF grant
DMS-0704359.}}

\thanks{ \today}

\begin{abstract}
This paper gives the first explicit, two--sided estimates on the cusp area of 
once--punctured torus bundles, 4--punctured sphere bundles, and
2--bridge link complements. The input for these estimates is 
purely combinatorial data coming from the Farey tesselation of the hyperbolic plane.
  The bounds on cusp
area lead to explicit bounds on the volume of Dehn fillings of these
manifolds, for example 
sharp bounds on volumes of hyperbolic closed 3--braids in terms of the
Schreier normal form of the associated braid word.  Finally, these
results are applied to derive relations between the Jones polynomial
and the volume of hyperbolic knots, and to disprove a related
conjecture.  
\end{abstract}

\maketitle

\begin{flushright}
\emph{Dedicated to the memory of Xiao-Song Lin}
\end{flushright}

\section{Introduction}\label{sec:intro}

Around 1980, Thurston proved that 3--manifolds with torus boundary
decompose into pieces that admit locally homogeneous geometric structures
\cite{thurston:survey}, and that in an appropriate sense the most
common such structure is hyperbolic \cite{thurston:notes}.  By
Mostow--Prasad rigidity, a hyperbolic structure is unique for such a
manifold, and thus the geometry of a hyperbolic manifold ought to give
a wealth of information to aid in its classification.  However, in
practice it has been very difficult to determine geometric properties
of a hyperbolic manifold from a combinatorial or topological
description.

In this paper, we address this problem for a class of 3--manifolds
that we call \emph{Farey manifolds}: punctured torus bundles,
4--punctured sphere bundles, and 2--bridge link complements.  The
combinatorial and geometric structure of these manifolds can be neatly
described in terms of the Farey tesselation of the hyperbolic
plane. For each type of Farey manifold, we use purely combinatorial
data coming from this tesselation to give the first explicit,
two-sided estimates on the area of a maximal cusp.

The bounds on cusp areas lead to explicit bounds on the volume of Dehn
fillings of Farey manifolds.  An example of such a Dehn filling is the
complement of a closed 3--braid.  We bound the volumes of such
manifolds, and in particular give sharp bounds on volumes of
hyperbolic closed 3--braids in terms of the Schreier normal form of
the associated braid word.  These results are applied to derive
relations between the Jones polynomial and the volume of hyperbolic
knots and to disprove a related conjecture.

\subsection{Cusp shapes and areas}

In a finite--volume hyperbolic 3--manifold $M$, a horoball
neighborhood of a torus boundary component becomes a \emph{cusp},
homeomorphic to $T^2\times[0,\infty)$.  Mostow--Prasad rigidity
implies that each cross-sectional torus $T^2$ is endowed with a flat
metric, or \emph{cusp shape}, that is determined up to similarity by
the topology of $M$.  When $M$ has a single torus boundary component,
we may expand a horoball neighborhood until it meets itself.  This
\emph{maximal horoball neighborhood} completely determines a flat
metric on the torus, and one can measure lengths of curves and area on
the torus using this metric.  We will refer to such a metric as a
\emph{maximal cusp metric}.  Similarly, when a 3--manifold has
multiple cusps, a maximal horoball neighborhood is given by expanding
a collection of horoball neighborhoods until none can be exapanded
further while keeping their interiors disjoint.  In the case of
multiple cusps, the choice of horoball neighborhoods is no longer
unique.  However, if we are required to expand the cusp neighborhoods 
in a fixed order, this expansion recipe once again determines a collection 
of maximal cusp metrics.

It is known, due to Nimershiem, that the set of similarity classes of
tori that can be realized as cusps of hyperbolic 3--manifolds is dense
in the moduli space of 2--tori \cite{nimershiem:flat-tori}.  However,
in general it is not known how to determine the cusp shape of a
manifold.  For simple manifolds, for example those built of a small
number of ideal tetrahedra, or links with a small number of crossings,
Weeks' computer program SnapPea will determine shapes of cusps and
maximal cusp metrics \cite{weeks:snappea}.  For other, larger classes
of 3--manifolds, some bounds on cusp shape have been obtained.
Aitchison, Lumsden, and Rubinstein found cusp shapes of certain
alternating links, but the metrics they used were singular
\cite{alr:alternating}.  For non-singular hyperbolic metrics, Adams
\emph{et al.} found upper bounds on the cusp area of knots, in terms
of the crossing number of a diagram \cite{adams-students:cusp-area}.
Purcell found that for ``highly twisted'' knots, the lengths of
shortest arcs on a maximal cusp metric are bounded above and below in
terms of the twist number of a diagram \cite{purcell:cusps}.  These
results were obtained using cusp estimates on a class of links called
fully augmented links, whose cusp shapes and lengths of slopes on
maximal cusp metrics were also worked out by Purcell
\cite{purcell:cusps} and Futer and Purcell \cite{futer-purcell}.

In this paper, we prove explicit, readily applicable bounds on cusp
shapes and maximal cusp metrics of punctured torus bundles and
4--punctured sphere bundles, as well as of 2--bridge knot complements.
These manifolds have a natural ideal triangulation, first discovered
for punctured torus bundles by Floyd and Hatcher \cite{floyd-hatcher},
and later studied by many others \cite{akiyoshi, aswy,
gueritaud:thesis, lackenby:punct}.  One feature that makes these
3--manifolds particularly attractive is that their geometry can be
described in terms of the combinatorics of the Farey tessellation of
$\HH^2$.  Hence, we refer to these manifolds as \emph{Farey
manifolds}.

To state an example of our results in this direction, let $M$ be a
hyperbolic once--punctured torus bundle.  The monodromy of $M$ can be
thought of as a conjugacy class in $SL_2(\ZZ)$.  As such, it has a
(unique up to cyclic permutation of factors) presentation of the form
$$ \pm \begin{bmatrix} 1 & a_1 \\ 0 & 1 \end{bmatrix}
\begin{bmatrix} 1 & 0 \\ b_1& 1 \end{bmatrix} \cdots \cdots \cdots
\begin{bmatrix} 1 & a_s \\ 0 & 1 \end{bmatrix}
\begin{bmatrix} 1 & 0 \\ b_s & 1\end{bmatrix},$$
where $a_i$, $b_i$ are positive integers.  The integer $s$ is called
the \emph{length} of the monodromy.

\begin{named}{Theorem \ref{thm:ptorus-area}}
Let $M$ be a punctured--torus bundle with monodromy of length $s$.
Let $C$ be a maximal horoball neighborhood about the cusp of $M$. Then
$$\frac{16 \sqrt{3}}{147} \, s \: \leq \: \area(\bdy C) \:<\: 2
\sqrt{3} \, \frac{v_8}{v_3} \, s,$$
where $v_3 = 1.0149...$ is the volume of a regular ideal tetrahedron
and $v_8 = 3.6638...$ is the volume of a regular ideal octahedron.

Furthermore, if $\gamma$ is a non-trivial simple closed curve on $\bdy
C$ (that is, any simple closed curve that is transverse to the
fibers), then its length $\ell(\gamma)$ satisfies
$$ \ell(\gamma) \: \geq \: \frac{4\sqrt{6}}{147} \, s.$$
\end{named}

The proof of Theorem \ref{thm:ptorus-area} contains two main
steps. First, we derive an estimate for the size of horoballs in the
universal cover of punctured torus bundles (see Proposition
\ref{prop:horosize}).  Then, we pack the cusp torus with the shadows
of these horoballs.

Our estimate on horoball size should be compared with J{\o}rgensen's
work on quasifuchsian punctured torus groups, which appears in a
well--known but unfinished manuscript \cite{jorgensen}.  A careful
exposition of J{\o}rgensen's work was given by Akiyoshi, Sakuma, Wada,
and Yamashita \cite{aswy-book}.  J{\o}rgensen's results can be applied
in our setting to show that the universal cover of a punctured torus
bundle contains a number of maximal horoballs whose size is bounded
from below (see \cite[Lemma 4.3]{jorgensen} and \cite[Lemma
8.1.1]{aswy-book}).  J{\o}rgensen conjectured the existence of a much
better lower bound for this horoball size; and indeed, our Proposition
3.6 improves J{\o}rgensen's lower bound by a factor of more than
$10$. This improvement is very important in our setting, since a
10--fold improvement in horoball size yields a 100--fold improvement
in the cusp area estimate.  See the end of Section
\ref{subsec:large-horosphere} for a more detailed discussion.

\subsection{Cusp area and link diagrams}\label{sec:diagram-intro}

A closely related class of manifolds are complements of 2--bridge
links.  Using similar techniques, in this paper we are also able to
bound the lengths of slopes on maximal cusps in hyperbolic 2--bridge
links.  Since all 2--bridge links can be represented by an alternating
diagram, our results give further evidence for a conjectural picture
of the cusp shapes and maximal cusp metrics of alternating knots.

For general alternating knots and links, there is increasing evidence
that the cusp shape and maximal cusp metric ought to be bounded in
terms of the twist number of a reduced diagram. We say that a link
diagram is \emph{reduced} if it does not contain any crossings that
separate the diagram: that is, any crossings in the projection plane
such that there is a simple closed curve meeting the diagram
transversely in only that crossing.  Similarly, two crossings are said
to be equivalent if there exists a simple closed curve meeting the
knot diagram transversely in those two crossings, disjoint from the
knot diagram elsewhere.  The \emph{twist number} is the number of
equivalence classes of crossings (called \emph{twist regions}).

\begin{conjecture}
The area of a maximal cusp metric on an alternating knot is bounded
above and below by linear functions of the twist number of a reduced,
alternating diagram.  Similarly, the length of the shortest
non-meridional slope of an alternating knot is bounded above and below
by a linear function of the twist number of the diagram.
\label{conj:alt-knot-cusp}
\end{conjecture}

We first became aware of this conjecture several years ago by viewing
slides of a talk by Thistlethwaite, in which he showed using SnapPea
that the conjecture holds for many simple alternating knots.  Lackenby
proved a close variant the conjecture, relating the twist number of a
diagram to the \emph{combinatorial length} of slopes
\cite{lackenby:word}.  However, Lackenby's methods are purely
combinatorial and cannot be applied to give the geometric information
of the conjecture.

In this paper we prove the conjecture for 2--bridge link complements.
In particular, we show the following.

\begin{named}{Theorem \ref{thm:2bridge}}
Let $K$ be a 2--bridge link in $S^3$, whose reduced alternating
diagram has twist number $t$.
Let $C$ be a maximal neighborhood about the cusps of $S^3 \setminus
K$, in which the two cusps have equal volume if $K$ has two
components.  Then
$$\frac{8\sqrt{3}}{147}\,(t-1) \: \leq \: \area(\bdy C) \:
<\: 2\sqrt{3}\,\, \frac{v_8}{v_3}\,(t-1).$$
Furthermore, if $K$ is a knot, let $\gamma$ be any non-trivial arc
that starts on a meridian and comes back to the same meridian (for
example, a non-meridional simple closed curve). Then its length
satisfies
$$\ell(\gamma) \: \geq \: \frac{4\sqrt{6\sqrt{2}}}{147}\, (t-1).$$
\end{named}

This result should be compared to that of Adams \emph{et al.}
\cite{adams-students:cusp-area}, where they prove upper bounds on cusp
area in terms of the crossing number $c$ of a knot.  For alternating
knots, including 2--bridge knots, they show that the cusp area
satisfies $\area(\bdy C) \leq 9c - 36 + 36/c$.  For those 2--bridge
knots whose diagrams have very few crossings per twist region (in
particular, when $c/t < 1.39$), the bound of Adams \emph{et al.} is
sharper than the upper bound of Theorem \ref{thm:2bridge}.  For more
general 2--bridge knots that have more crossings per twist region, the
upper bound of Theorem \ref{thm:2bridge} is a significant improvement.
To the best of our knowledge, the lower bound of Theorem
\ref{thm:2bridge} does not have any predecessors in the literature.

\subsection{Applications to hyperbolic Dehn filling}

The shapes of the cusps and their actual metrics give information not
just on the 3--manifold itself, but also on the Dehn fillings of that
manifold.

For example, modulo the geometrization conjecture, several theorems
imply that Dehn fillings on slopes of sufficient length yield
hyperbolic manifolds (these are the $2\pi$--Theorem, due to Gromov and
Thurston \cite{bleiler-hodgson}; the 6--Theorem, due to Agol
\cite{agol:bounds} and Lackenby \cite{lackenby:word}; and the
$7.515$--Theorem, due to Hodgson and Kerckhoff \cite{hk:univ}).  When
we combine these theorems with the results on maximal cusp areas and
slope lengths above, we find that Farey manifolds with long monodromy
admit no non-trivial Dehn fillings, where ``long'' is explicit.

In particular, Bleiler and Hodgson \cite{bleiler-hodgson} note that
the work of J{\o}rgensen \cite{jorgensen} combined with the
$2\pi$--Theorem implies that there is a constant $N$ such that every
non-trivial Dehn filling of a punctured torus bundle with monodromy
length $s>N$ gives a hyperbolic 3--manifold.  However, they remark on
the lack of an explicit value for the constant $N$.  Now Theorem
\ref{thm:ptorus-area}, coupled with the 6--Theorem, allows the
estimate $N\leq 90$.

More recently, the authors proved a result that bounds the volume of
manifolds obtained by Dehn filling along a slope of length at least
$2\pi$, in terms of the length of that slope \cite{fkp-07}.  Thus we
may combine Theorem \ref{thm:ptorus-area} with this recent result to
estimate the volumes of the manifolds obtained by Dehn filling.  For
example, if $M$ is a punctured--torus bundle with monodromy of length
$s>94$, then the length of any non-trivial slope $\gamma$ on the cusp
of $M$ (i.e. any slope transverse to the fibers) will be at least
$2\pi$.  Then by \cite[Theorem 1.1]{fkp-07}, the volume of the
manifold $M(\gamma)$ obtained by Dehn filling $M$ along $\gamma$ will
be bounded explicitly below.  See Corollary \ref{cor:ptorus-volume}.

One large class of examples obtained by Dehn filling 4--punctured
sphere bundles is the class of closed 3--braids, which has been
extensively studied by others (see e.g.\ Murasugi
\cite{murasugi:3braids}, Birman and Menasco
\cite{birman-menasco:3braids}).  In this paper, we classify the
hyperbolic links that are closed 3--braids (see Theorem
\ref{thm:hyperbolic-3braid}), and obtain the first estimates on
volumes of these links.

\begin{figure}[b]
\psfrag{s}{$\sigma_1$}
\psfrag{t}{$\sigma_2$}
\begin{center}
\includegraphics{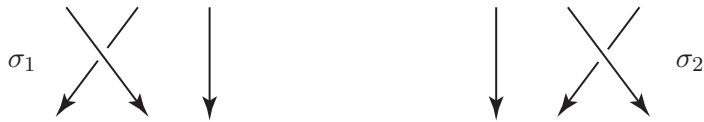}
\end{center}
\caption{Braid group generators $\sigma_1$ and $\sigma_2$.}
\label{fig:3braid-generators}
\end{figure}
%%%%%%%%%%%%%%%%%%%%%%%%%%%%%%%%%%%%%%%%%%%%%%%%%%%%%%%%%%%%%%%%%

To state these results, let $\sigma_1$, $\sigma_2$ denote the
generators for $B_3$, the braid group on three strands, as in Figure
\ref{fig:3braid-generators}.  Let $C= (\sigma_1 \sigma_2)^3$.  It is
known, by work of Schreier \cite{schreier}, that most $3$--braids are
conjugate to words of the form $w = C^k \sigma_1^{p_1} \sigma_2^{-q_1}
\cdots \sigma_1^{p_s} \sigma_2^{-q_s}$, where $p_i$, $q_i$ are all
positive. In particular, all $3$--braids with hyperbolic closures are
of this form, as we shall show in Theorem \ref{thm:hyperbolic-3braid}.
Following Birman and Menasco \cite{birman-menasco:3braids}, we call
such braids \emph{generic}.

\begin{named}{Theorem \ref{thm:3braid-volume}}
Let $K= \hat{w}$ be the closure of a generic 3--braid $w = C^k
\sigma_1^{p_1} \sigma_2^{-q_1} \cdots \sigma_1^{p_s} \sigma_2^{-q_s}$,
where $C= (\sigma_1 \sigma_2)^3$, and $p_i$, $q_i$ are all positive.
Suppose, furthermore, that $w$ is not conjugate to $\sigma_1^{p}
\sigma_2^q$ for arbitrary $p$, $q$.  Then $K$ is hyperbolic, and
$$
4v_3 \, s - 277 \: < \: \vol(S^3 \setminus K) \: < \: 4v_8 \, s,
$$
where $v_3 = 1.0149...$ is the volume of a regular ideal tetrahedron
and $v_8 = 3.6638...$ is the volume of a regular ideal octahedron.
Furthermore, the multiplicative constants in both the upper and lower
bounds are sharp.
\end{named}

\subsection{Volume and Jones polynomial invariants}

The volume estimate of Theorem \ref{thm:3braid-volume} has a very
interesting application to conjectures on the relationship of the
volume to the Jones polynomial invariants of hyperbolic knots.

For a knot $K$, let
$$J_K(t)= \alpha_K t^m+ \beta_K t^{m-1}+ \ldots + \beta'_K t^{r+1}+
\alpha'_K t^r$$
denote the Jones polynomial of $K$.  We will always denote the second
and next--to--last coefficients of $J_K(t)$ by $\beta_K$ and
$\beta'_K$, respectively.

The Jones polynomial fits into an infinite family of knot invariants:
the \emph{colored Jones} polynomials.  These are Laurent polynomial
knot invariants $J^n_K(t)$, $n>1$, where $J^2_K(t)=J_K(t)$.  The
volume conjecture \cite{kashaev:volume-conj, murakami-squared} states
that for a hyperbolic knot $K$,
$$2\pi \lim_{n\to \infty} { \frac{\log \abs {J^n_K(e^{2\pi i / n} )
}} {n}}= \vol(S^3 \setminus K),$$
where $e^{2\pi i /n}$ is a primitive $n$-th root of unity.  If the
volume conjecture is true, then one expects correlations between
$\vol(S^3 \setminus K)$ and the coefficients of $J_K^n(t)$, at least
for large values of $n$.  For example, for $n\gg 0$ one would have
$$\vol(S^3 \setminus K)\: < \: C || J_K^n ||,$$
where $|| J_K^n ||$ denotes the sum of absolute values of the
coefficients of $J_K^n(t)$ and $C$ is a constant independent of $K$.
At the same time, several recent results and much experimental
evidence \cite{ckp:simplest-knots, dasbach-lin:volumeish, fkp-07,
fkp:conway} actually indicate that there may be a correlation between
$\vol(S^3\setminus K)$ and the coefficients of the Jones polynomial
itself.  These results prompt the following question.

\begin{question}\label{question}
Do there exist constants $C_i>0$, $i=1, \ldots, 4$, and a function
$B_K$ of the coefficients of $J_K(t)$, such that all hyperbolic knots
satisfy
\begin{equation}
	C_1\ B_K-C_2 \: < \: \vol(S^3 \setminus K) \: < \: C_3\ B_K+C_4?
\label{eqn:question}
\end{equation}
\end{question}

Dasbach and Lin \cite{dasbach-lin:volumeish} showed that for
alternating knots, equation (\ref{eqn:question}) holds for $B_K:=\abs
{\beta_K}+\abs {\beta'_K}$.  They also presented experimental evidence
suggesting linear correlations between $\abs {\beta_K}+\abs
{\beta'_K}$ and the volume of non-alternating knots; their data is
based on knots with a low numbers of crossings.  The authors of the
current paper have shown that the same function works for several
large families in the class of \emph{adequate knots}, which are a vast
generalization of alternating knots \cite{fkp-07,fkp:conway}.  In
fact, Dasbach and Lin \cite{dasbach-lin:head-tail} and Stoimenow 
\cite{stoimenow:amphichiral}
showed that for adequate knots, the second and
next--to--last coefficients of the colored Jones polynomial $J_K^n(t)$
are independent of $n$, equal to those of the Jones polynomial
$J_K(t)$. So these results establish strong versions of
relations between volume and coefficients of the colored Jones
polynomials for these knots, as predicted by the volume conjecture.
This led to some hope that not only would Question \ref{question} be
answered in the affirmative, but also that $B_K = \abs{\beta_K} +
\abs{\beta'_K}$ could always work in equation (\ref{eqn:question}).

In this paper we show that  a slightly modified function,  involving
the first two and the last two coefficients of the Jones polynomial, 
satisfies equation (\ref{eqn:question}) for hyperbolic closed 3-braids.
Building on Theorem \ref{thm:3braid-volume}, we prove the following.

\begin{named}{Theorem \ref{thm:relations}}
Let $K$ be a hyperbolic closed $3$--braid. 
From the Jones polynomial $J_K(t)$, we define ${\zeta_K}, {\zeta'_K}$ as follows. Let
$$\zeta_K = \left\{ \begin{array}{r l}
\beta_K, & \mbox{if } \abs{\alpha_K} = 1 \\
0, & \mbox{otherwise}
\end{array} \right.
\qquad \mbox{and} \qquad
\zeta'_K = \left\{ \begin{array}{r l}
\beta'_K, & \mbox{if } \abs{\alpha'_K} = 1 \\
0, & \mbox{otherwise.}
\end{array} \right.$$
Define $\zeta = \max \left\{ \abs{\zeta_K}, \abs{\zeta'_K} \right\}$. Then
$$4v_3 \cdot \zeta - 281 \: < \: \vol(S^3 \setminus K) \: < \:
4v_8 \, ( \zeta +1 ).$$
Furthermore, the multiplicative constants in both the upper and lower bounds are sharp.
\end{named}

We remark that the quantity $\zeta = \max \left\{ \abs{\zeta_K}, \abs{\zeta'_K}  \right\} $, as defined in 
Theorem \ref{thm:relations}, will also serve to estimate the volumes of alternating links,
sums of alternating tangles, and highly twisted adequate links. In other words, this quantity
estimates the volume of every family of knots and links where the volume is known to be bounded
above and below in terms of the Jones polynomial. The quantity $\zeta$ depends only 
on the first two and last two coefficients of $J_K(t)$, and can be taken as positive evidence for
Question \ref{question}. On the other hand, we also  show that no function of $\beta_K$ and $\beta'_K$ alone can
satisfy equation (\ref{eqn:question}) for all hyperbolic knots.

\begin{named}{Theorem \ref{thm:counterexample}}
There does not exist a function $f(\cdot, \cdot)$ of two variables,
together with constants $C_i>0$, $i=1, \ldots, 4$, such that
$$C_1 f(\beta_K, \beta'_K) - C_2 \: < \: \vol(S^3\setminus K) \:<\:
C_3 f(\beta_K, \beta'_K) + C_4$$
for every hyperbolic knot $K$. In other words, the second and
next--to--last coefficients of the Jones polynomial do not coarsely
predict the volume of a knot.
\end{named}

Theorem \ref{thm:counterexample} relies on two families of examples:
adequate Montesinos knots and closed 3--braids.  For both of these families,
Stoimenow found upper bounds on volume in terms of outer coefficients of the 
Jones polynomial \cite{stoimenow:amphichiral}. While
Theorem \ref{thm:counterexample} implies that equation (\ref{eqn:question}) 
cannot hold for any function of $\beta_K$ and $\beta'_K$ alone, there might still be 
an affirmative answer to Question \ref{question} that uses other coefficients.

\subsection{Organization}

We begin by discussing Farey manifolds.  In Section
\ref{sec:farey-ford}, we describe the canonical triangulations of the
three families of Farey manifolds.  In Section \ref{sec:horosphere},
we show that the universal cover of one of these manifolds must
contain a number of maximal horospheres whose size is bounded below.
This leads to the cusp area estimates of Section \ref{sec:cusp-area}.

The later sections give applications of these cusp area estimates.  In
Section \ref{sec:volume3braids}, we apply the results on cusp area to
estimate the volumes of closed 3--braids.  Finally, in Section
\ref{sec:jones3braids}, we combine this with a discussion of Jones
polynomials.

\subsection{Acknowledgements}

We thank Ian Agol for a number of helpful comments, and in particular
for getting us started in the right direction towards Proposition
\ref{prop:horosize}.  We thank Fran\c{c}ois Gu\'eritaud for
enlightening discussions about punctured torus bundles.  In order to generate 
Table \ref{table:jones} on page \pageref{table:jones},
we used 
% Jones polynomial 
software written by Dror Bar-Natan and Nathan Broaddus, 
as well as a handy script by Ilya Kofman. 
 We are grateful to all of them. 
Finally, we thank 
the referees for suggesting a number of revisions that improved 
this paper.

\section{The canonical triangulation of a Farey manifold}\label{sec:farey-ford}

In this section, we review the canonical triangulations of Farey
manifolds. We begin by recalling the definition of the Ford domain and
the canonical polyhedral decomposition that is its dual. We then
describe the combinatorics of the canonical polyhedral decomposition for
each of the three families of Farey manifolds; for the manifolds in
question, it is always a triangulation. Along the way, we introduce a
number of terms and notions that will be needed in the ensuing
arguments.

\subsection{The Ford domain and its dual}\label{sec:ford-rev}

For a hyperbolic manifold $M$ with a single cusp, expand a horoball
neighborhood about the cusp.  In the universal cover $\HH^3$, this
neighborhood lifts to a disjoint collection of horoballs.  In the
upper half space model for $\HH^3$, we may ensure that one of these
horoballs is centered on the point at infinity.  Select vertical
planes in $\HH^3$ that cut out a fundamental region for the action of
the $\ZZ \times \ZZ$ subgroup of $\pi_1(M)$ that fixes the point at
infinity.  The Ford domain is defined to be the collection of points
in such a fundamental region that are at least as close to the
horoball about infinity as to any other lift of the horoball
neighborhood of the cusp.

The Ford domain is canonical, except for the choice of fundamental
region of the action of the subgroup fixing infinity.  It is a
finite--sided polyhedron, with one ideal vertex.  The faces glue
together to form the manifold $M$.

If the manifold $M$ has several cusps, the above construction still
works, but is less canonical.  Once one chooses a horoball
neighborhood of each cusp, as well as a fundamental domain for each
cusp torus, the nearest--horoball construction as above produces a
fundamental domain for $M$.  This fundamental domain is a disjoint union 
of finite--sided polyhedra, with one polyhedron for each cusp of $M$
 and one ideal vertex per polyhedron.
We refer to this fundamental domain as a Ford domain determined by the choice
of horoball neighborhood.

Dual to the Ford domain is a decomposition of $M$ into ideal
polyhedra.  This decomposition, first studied by Epstein and Penner
\cite{epstein-penner}, is canonically determined by the relative
volumes of the cusp neighborhoods.  In particular, if $M$ has only one
cusp, the decomposition dual to the Ford domain is completely
canonical.  We refer to it as the \emph{canonical polyhedral
decomposition}.

One of the few infinite families for which the canonical polyhedral
decomposition is completely understood is the family of Farey
manifolds.  For once--punctured torus bundles and 4--punctured sphere
bundles, the combinatorial structure of this ideal triangulation was
first described by Floyd and Hatcher \cite{floyd-hatcher}.  Akiyoshi
\cite{akiyoshi} and Lackenby \cite{lackenby:punct} gave distinct and
independent proofs that the combinatorial triangulation is
geometrically canonical, i.e. dual to the Ford domain.  Gu\'eritaud
used the combinatorics of the triangulation to determine by direct
methods those punctured torus bundles that admit a hyperbolic
structure \cite{gf:two-bridge}; he also re-proved that the
Floyd--Hatcher triangulation is canonical \cite{gueritaud:thesis}.

For two-bridge link complements, the analogue of the Floyd--Hatcher
triangulation was described by Sakuma and Weeks \cite{sakuma-weeks}.
Following Gu\'eritaud's ideas, Futer used this triangulation to find a
hyperbolic metric for all the 2--bridge link complements that admit
one \cite[Appendix]{gf:two-bridge}.  Akiyoshi, Sakuma, Wada, and
Yamashita \cite{aswy} and (independently) Gu\'eritaud
\cite{gueritaud:thesis} showed that the Sakuma--Weeks triangulation is
geometrically canonical.  For all of the Farey manifolds, our
exposition below follows that of Gu\'eritaud and Futer, and we refer
the reader to reference \cite{gf:two-bridge} for more details.

\subsection{Once--punctured torus bundles}\label{sec:ptorus}

Let $V_\varphi$ be a hyperbolic punctured torus bundle with monodromy
$\varphi$.  The mapping class group of the punctured torus is
isomorphic to $SL_2(\ZZ)$.  By a well--known argument that we recall
below, either $\varphi$ or $-\varphi$ is conjugate to an element of
the form
$$\Omega = R^{a_1}L^{b_1}\cdots R^{a_s}L^{b_s},$$
where $a_i, b_i$ are positive integers, and $R$ and $L$ are the
matrices
$$R := \begin{bmatrix} 1& 1 \\ 0&1\end{bmatrix}, \quad
	L := \begin{bmatrix} 1& 0 \\ 1& 1 \end{bmatrix}.$$
Moreover, $\Omega$ is unique up to cyclic permutation of its letters.

By projecting $\varphi$ down to $PSL_2(\ZZ) \subset Isom(\HH^2) $, we
may view the matrix $\pm \varphi$ as an isometry of $\HH^2$ in the
upper half--plane model, where the boundary at infinity of $\HH^2$ is
$\RR \cup \{ \infty \}$.  Then the slopes of the eigenvectors of
$\varphi$ are the fixed points of its action on $\overline{\HH^2}$.

Now, subdivide $\HH^2$ into ideal triangles, following the \emph{Farey
tesselation} $\farey$.  In this tesselation, every vertex is a
rational number (or $\infty$) in $\bdy \HH^2$.  Each such vertex
corresponds to a \emph{slope} on the punctured torus $T$, that is, an
isotopy class of arcs running from the puncture to itself.  Two
vertices are connected by an edge in $\farey$ if and only if the
corresponding arcs can be realized disjointly.  Thus an ideal triangle
of $\farey$ corresponds to a triple of disjoint arcs, which gives an
ideal triangulation of $T$.  The monodromy $\varphi$ naturally acts on
$\farey$.

There is an oriented geodesic $\gamma_\varphi$ running from the
repulsive fixed point of $\varphi$ to its attractive fixed point.
This path crosses an infinite sequence of triangles of the Farey graph
$(\dots, t_{-1}, t_0, t_1, t_2, \dots)$.  We can write down a
bi-infinite word corresponding to $\varphi$, where the $k$-th letter
is $R$ (resp. $L$) if $\gamma_\varphi$ exits the $k$-th triangle $t_k$
to the right (resp. left) of where it entered.  This bi-infinite word
will be periodic of period $m$, where $m$ is some integer such that
$t_0$ is taken by $\varphi$ to $t_m$.  Then letting $\Omega$ be any
subword of length $m$, and substituting the matrices above for $R$ and
$L$, we find that $\pm \Omega$ is conjugate to $\varphi$.

\begin{figure}
\centerline{\includegraphics{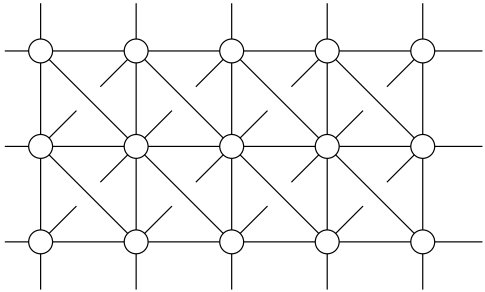}}
\caption{Copies of ideal tetrahedra in the cover $(\RR^2\setminus
	\ZZ^2)\times \RR$ of $T\times \RR$.}
\label{fig:tetrahedra}
\end{figure}

Next, we review the relation between the word $\Omega$ and the
triangulation of $V_\varphi$.  The path $\gamma_\varphi$ through the
Farey graph determines a sequence of triangulations of the punctured
torus $T$.  Every time $\gamma_\varphi$ crosses an edge $e \subset
\farey$, moving from one triangle of $\farey$ to an adjacent triangle,
we change one ideal triangulation of $T$ (call it $\tau_-(e)$) into a
different ideal triangulation $\tau_+(e)$, replacing a single edge
with another.  In other words, we are performing a diagonal exchange
in a quadrilateral of $T$.  This diagonal exchange determines an ideal
tetrahedron $\Delta(e)$ as follows.  The boundary of the tetrahedron
is made up of two pleated surfaces homotopic to $T$, with
triangulations corresponding to $\tau_-(e)$ and $\tau_+(e)$.  These
two pleated surfaces are glued together along the two edges in $T$
where $\tau_-(e)$ and $\tau_+(e)$ agree.  The result is an ideal
tetrahedron.  See Figure \ref{fig:tetrahedra}.

%%%%%%%%%%%%%%%%%%%%%%%%%%%%%%%%%%%%%%%%%%%%%%%%%%%%%%%%%%%%%%%%%

If $\gamma_\varphi$ crosses the edges $e_i$, $e_{i+1}$, then we may
glue $\Delta(e_i)$ to $\Delta(e_{i+1})$ top to bottom, since
$\tau_+(e_i) = \tau_-(e_{i+1})$.  Thus $\gamma_\varphi$ determines a
bi-infinite stack $U$ of tetrahedra.  $U$ is homeomorphic to $T\times
\RR$, and there is an orientation--preserving homeomorphism $\Phi$ of
$U$, taking the $i$-th tetrahedron to the $(i+m)$-th tetrahedron,
acting as $\varphi$ on $T$.  The quotient $U/\Phi$ is homeomorphic to
$V_\varphi$, and gives a triangulation of $V_\varphi$ into $m$ ideal
tetrahedra.  This is the Floyd--Hatcher triangulation of $V_\varphi$,
also called the \emph{monodromy triangulation}.

We summarize the discussion above as follows.
\begin{enumerate}
\item The monodromy $\varphi$ of the bundle is conjugate to a word
	$$\Omega = \pm R^{a_1}L^{b_1}\cdots R^{a_s}L^{b_s}.$$
\item Each letter $R$ or $L$ corresponds to a triangle in the Farey
	tesselation of $\HH^2$.
\item Each letter $R$ or $L$ corresponds to a pleated surface
homotopic to $T$, pleated along arcs whose slopes are the vertices of
the corresponding triangle of the Farey graph.  This pleated surface
forms the boundary between two tetrahedra of the canonical
triangulation of $V_\varphi$.
  \end{enumerate}

\begin{define}
Let $\Omega = \pm R^{a_1}L^{b_1}\dots R^{a_s}L^{b_s}$.  A
\emph{syllable} of $\Omega$ is defined to be a subword $R^{a_i}$ or
$L^{b_i}$.  That is, a syllable is a maximal string of $R$'s or $L's$
in the word $\Omega$.
\end{define}

\begin{figure}%[h]
\begin{center}
\input{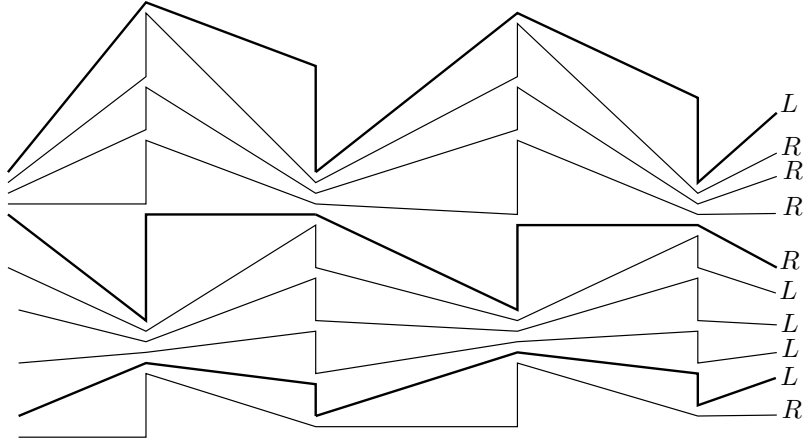}
\end{center}
\caption{Zigzags on the boundary torus of a punctured torus
bundle. The pleated surfaces in this figure correspond to a sub-word
$R L^4 R^4 L \subset \Omega$.}
\label{fig:stack}
\end{figure}

A punctured torus bundle is a manifold with a single torus boundary
component.  It is often convenient to work with the universal cover
$\HH^3$ of the bundle, seen as the upper half space model, with the
boundary lifting to the point at infinity in this model.  Each of the
pleated surfaces corresponding to the letters $R$ and $L$ will lift to
$\HH^3$.  Their intersection with the boundary of a maximal cusp gives
a triangulation of the boundary which is well understood.  In
particular, these intersections give a collection of \emph{zigzags}
that determine a triangulation of the boundary with combinatorics
specified by the word $\Omega$.  See Figure \ref{fig:stack}.

\begin{define}
A \emph{zigzag} is defined to be the lift of one of the pleated
surfaces corresponding to $R$ or $L$ to the universal cover $\HH^3$,
with the cusp lifting to infinity.
\end{define}

Note that in pictures of zigzags, as in Figure \ref{fig:stack}, the
vertices of the zigzag correspond to edges in $\HH^3$ along which
these zigzags meet.  To distinguish separate zigzags, it is
conventional to split them apart at the vertices.  (In the manifold
$M$, a sequence of pleated surfaces corresponding to a syllable of
$\Omega$ will meet along a single edge.  Thus, in a more topologically
accurate but less enlightening picture, one would collapse together
the split--apart vertices in Figure \ref{fig:stack}.  See also
\cite[Figure 4]{gf:two-bridge}.)

Akiyoshi \cite{akiyoshi}, Lackenby \cite{lackenby:punct}, and
Gu\'eritaud \cite{gueritaud:thesis} have independently proved that
this triangulation is geometrically canonical, i.e.\ dual to the Ford
domain.  As a result, each edge of the triangulation runs through the
geometric center of a face of the Ford domain.
(More precisely, each face of the Ford domain lifts to a hemisphere in
$\HH^3$, and each edge of the triangulation runs through the geometric
center of the hemisphere.)  Thus, when viewed from infinity, the
``corners'' of the zigzag lie over centers of hemispheres projecting
to faces of the Ford domain.  We will use this extensively below.
See, for example, Figure \ref{fig:zigzag} below.

\subsection{4--punctured sphere bundles}\label{sec:4ps}

Consider the universal abelian cover $X := \RR^2 \setminus \ZZ^2$ of
the punctured torus, and define the following transformations of $X$:
$$\alpha(x,y) = (x+1, y), \quad \beta(x,y) = (x,y+1), \quad
\sigma(x,y) = (-x,-y).$$
Then one obtains the punctured torus as $T = X/\langle\alpha, \beta
\rangle$ and the $4$--punctured sphere as $S = X/\langle\alpha^2,
\beta^2, \sigma \rangle$.  Both $S$ and $T$ are covered by the
4--punctured torus $R = X/\langle\alpha^2, \beta^2 \rangle$.  Then the
action of $SL_2(\ZZ)$ on $T$ lifts to an action on $X$, and descends
to an action on both $R$ and $S$.  As a result, every hyperbolic
punctured torus bundle $M$ is commensurable with a hyperbolic
4--punctured sphere bundle $N$, whose monodromy can be described by
the same word $\Omega$.  The common cover is a 4--punctured torus
bundle $P$.  In Figure \ref{fig:tetrahedra}, we see lifts of two
pleated surfaces to the common cover.

The $4$--punctured sphere bundle $N$ can have anywhere from one to
four cusps, depending on the action of its monodromy on the punctures
of $S$.  Thus, for the purpose of discussing Ford domains, it is
important to choose the right horoball neighborhood of the cusps.
Unless stated otherwise (e.g.\ in Theorem \ref{thm:4ps-area}), we
shall always choose the cusp neighborhood in $N$ that comes from
lifting a maximal cusp of the corresponding punctured torus bundle $M$
to the 4--punctured torus bundle $P$, and then projecting down to $N$.
We call this the \emph{equivariant} cusp neighborhood of a
$4$--punctured sphere bundle.

By lifting the canonical monodromy triangulation of $M$ to $P$, and
projecting down to $N$, we obtain the layered monodromy triangulation
of a $4$--punctured sphere bundle.  Every tetrahedron $\Delta(e)$ of
this triangulation lifts to a layer of four tetrahedra in $P$, and
projects down to a layer of two tetrahedra in $N$.  (See \cite[Figure
16]{gf:two-bridge}.)  This triangulation is still geometrically
canonical: it is dual to the Ford domain determined by the the
equivariant cusp neighborhood.  We refer to this Ford domain as an
\emph{equivariant} Ford domain.  In particular, it still makes sense
to talk about ``syllables'', ``zigzags'', etc.\ in relation to
4--punctured sphere bundles.  Note that because of the rotational
action of $\sigma$, a loop around a puncture of the fiber will only
cross three edges in zigzag of $N$, instead of six edges as in a
zigzag of the punctured torus bundle $M$.

\begin{define}\label{def:fibered-meridian}
In a punctured torus bundle $M$ or a $4$--punctured sphere bundle $N$,
call the loop about a puncture of the fiber the \emph{meridian} of the
corresponding manifold.  We shall denote the length of a meridian in a
maximal cusp of a punctured torus bundle $M$ by $2\mu$.  With this
convention, the meridian of the corresponding $4$--punctured sphere
bundle $N$ will have length $\mu$ in a maximal equivariant cusp.
\end{define}

In the discussion of the geometry below, we will switch between
descriptions of 4--punctured sphere bundles and punctured tori,
depending on which leads to the simplest discussion.  Because of the
covering property, results on the geometry of the universal cover will
apply immediately to both types of manifolds.

\subsection{Two-bridge links}\label{sec:2bridge}

If a $4$--punctured sphere bundle $N$ is cut along a pleated fiber
$S$, the result is a manifold homeomorphic to $S \times I$, equipped
with an ideal triangulation. To recover $N$, we reglue the top of this
\emph{product region} $S \times I$ to the bottom along faces of this
triangulation. Meanwhile, the complement of a two-bridge link $K$ also
contains a product region $S \times I$: namely, the complement of the
$4$--string braid that runs between the minima and maxima in a diagram
of $K$. It turns out that the combinatorics of this braid once again
defines a layered triangulation of the product region, and that a
particular folding of the top and bottom faces of $S \times I$ yields
the canonical triangulation of $S^3 \setminus K$.

\begin{figure}%[h]
\begin{center}
\psfrag{a}{(a)}
\psfrag{b}{(b)}
\includegraphics{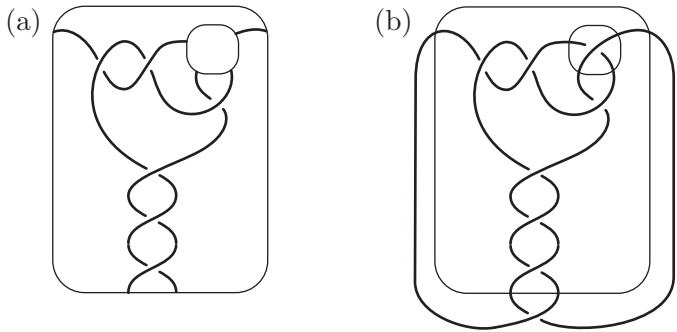}
\end{center}
\caption{(a) An alternating braid between two pillowcases, described
by the word $\Omega = R^3 L^2 R$. (b) The corresponding two-bridge
link $K(\Omega)$.}
\label{fig:pillow-braid}
\end{figure}

A 4--punctured sphere $S$ can be viewed as a square pillowcase with
its corners removed.  Consider two such nested pillowcases, with an
alternating $4$--string braid running between them, as in Figure
\ref{fig:pillow-braid}(a). The combinatorics of this braid, as well as
of the complementary product region $S \times I$, may be described by
a (finite) monodromy word of the form $\Omega = R^{p_1}L^{q_1}\dots
R^{p_s}L^{q_s}$, as above, where the $p_i$, $q_i$ are all positive,
except $p_1$ and $q_s$ are non-negative.  Each syllable $R^{p_i}$ or
$L^{q_i}$ determines a string of crossings in a vertical or horizontal
band, corresponding to a \emph{twist region} in which two strands of
the braid wrap around each other $p_i$ times.  To complete this
picture to a link diagram, we connect two pairs of punctures of the
outside pillowcase together with a crossing, and connect two pairs of
punctures of the inside pillowcase together with a crossing, as in
Figure \ref{fig:pillow-braid}(b).  This creates an alternating diagram
of a 2--bridge link $K(\Omega)$. It is well--known that any 2--bridge
link can be created in this manner (see, for example, Murasugi
\cite[Theorems 9.3.1 and 9.3.2]{murasugi}).

Just as in Sections \ref{sec:ptorus} and \ref{sec:4ps} above, the
monodromy word $\Omega$ describes a layered ideal triangulation of the
product region $S \times I$.  To form a 4--punctured sphere bundle,
one would glue the outer pillowcase $S_1$ to the inner pillowcase
$S_c$.  To obtain a 2--bridge link complement, we fold the surface
$S_1$ onto itself, identifying its four ideal triangles in pairs.
(See \cite[Figure 17]{gf:two-bridge}.)  We perform the same folding
for the interior pillowcase $S_{c}$.  This gives the desired canonical
triangulation of $S^3 \setminus K(\Omega)$.

Now consider the combinatorics of the cusp triangulation.  The pleated
surfaces between $S_1$ and $S_c$ are 4--punctured spheres with
combinatorics identical to that of the 4--punctured sphere bundle with
the same monodromy.  The universal cover of the product region looks
like a stack of zigzags, as in Figure \ref{fig:stack}.  (Just as with
$4$--punctured sphere bundles, a meridian of $K$ crosses three edges
of a zigzag -- so Figure \ref{fig:stack} shows two meridians.)  The
folding along $S_1$ and $S_c$ creates ``hairpin turns'', as in
\cite[Figure 19]{gf:two-bridge}.

Note that when $K$ is a two-component link, we shall always choose the
two cusp neighborhoods of $K$ to have equal volume, following the same
principle as in Section \ref{sec:4ps}.  This \emph{equivariant cusp}
neighborhood is the one whose Ford domain is dual to the layered
triangulation described above.  Also, because the symmetry group of $K$
interchanges the two cusps, it does not matter which cusp we look at
in the calculations of Section \ref{sec:horosphere}.

Finally, it is worth remarking that every surface $S_i$, lying between
two layers of tetrahedra, is a bridge sphere for the link $K$, and is
thus compressible in $S^3 \setminus K$.  Despite being compressible,
$S_i$ can nevertheless be realized as a pleated surface in the
geometry of $S^3 \setminus K$.  With the exception of the folded
surfaces $S_1$ and $S_c$, every other pleated $S_i$ is embedded, and
carries the same geometric information as the incompressible fiber in
a $4$--punctured sphere bundle.

\section{Geometric estimates for Ford domains}\label{sec:horosphere}

This section contains a number of geometric estimates on the Ford
domains of Farey manifolds.  We begin with a few estimates (Lemmas
\ref{lemma:triangle-ineq}--\ref{lemma:length-non-canonical}) that
apply to all triangulated cusped hyperbolic manifolds, and are
generally known to hyperbolic geometers.  We then restrict our
attention to Farey manifolds, and establish several estimates about
their Ford domains.  The main result of this section is Proposition
\ref{prop:horosize}: every zigzag contains an edge whose length
outside a maximal cusp is universally bounded.

\subsection{Estimates for triangulated hyperbolic 3--manifolds}

Recall from Section \ref{sec:ford-rev} that the Ford domain of a
cusped hyperbolic manifold $M$ is a union of finite--sided polyhedra, with one
ideal vertex for each cusp of $M$.  Consider those faces of the Ford
domain which do not meet an ideal vertex.  These consist of points
that are equidistant from two or more lifts of a cusp into $\HH^3$.
Each such face is the portion of a geodesic plane in $\HH^3$ which can
be ``seen'' from infinity.  That is, the geodesic planes are Euclidean
hemispheres centered on points of $\CC$ (here we are considering the
boundary at infinity of $\HH^3$ to be $\CC \cup \{\infty\}$), of some
Euclidean radius.  These overlap to cover all of $\CC$.  Looking down
from infinity, one sees portions of these Euclidean spheres.  These
are the faces.  The intersections of two adjacent faces give edges.  The
intersections of edges are vertices.

These faces of the Ford domain glue together in pairs.  Each pair of
faces consists of two hemispheres with identical Euclidean radii,
which glue together by some isometry of $\HH^3$.  In fact this
isometry can be taken to be a reflection in the face of the Ford
domain, followed by a Euclidean reflection (i.e. reflection in the
vertical plane that is the perpendicular bisector of the geodesic
connecting centers of the two hemispheres), followed by a rotation.
See, for example, Maskit's book \cite[Chapter IV, Section G]{maskit}.

We will be interested in the sizes of the radii, as well as distances
between centers of the Euclidean hemispheres that give the faces of
the Ford domain.  For our applications, the cell decomposition dual to
the Ford domain is always an actual triangulation, hence we shall talk
about triangles and tetrahedra.

Now, suppose $S_1$ and $S_2$ are two adjacent faces, which are
Euclidean hemispheres of radius $R_1$ and $R_2$, respectively, and
whose centers are Euclidean distance $D$ apart.

\begin{lemma}\label{lemma:triangle-ineq}
$R_1$, $R_2$, and $D$ as above satisfy the triangle inequality:
	$$R_1 + R_2 > D, \quad R_2 + D > R_1, \quad D+R_1 > R_2.$$
\end{lemma}

\begin{proof}
If $D \geq R_1+R_2$, the two faces $S_1$ and $S_2$ do not meet,
contradicting the fact that they are adjacent.

If $R_2+ D < R_1$, then the hemisphere $S_2$ lies completely inside
the region bounded by the complex plane and the hemisphere $S_1$.
Thus $S_2$ cannot be a face of the Ford domain.  This is a
contradiction.  By a symmetric argument, $D+R_1 > R_2$.
\end{proof}

Now consider the geometric dual of the Ford domain, which we will
assume is an actual triangulation.  In the universal cover, this dual
is given by taking an ideal vertex at the center of each Euclidean
hemisphere face of the Ford domain, and one at infinity.  There is one
edge for each hemisphere face of the Ford domain: a geodesic running
from infinity down to the center of the hemisphere.  For each
intersection of two faces of the Ford domain, there is a 2--cell.  By
assumption, when we project to the manifold these 2--cells become
ideal triangles.  Similarly, the intersection of three adjacent faces
is dual to a 3--cell which projects to an ideal tetrahedron.  Finally,
note the geometric dual may not be realized as a combinatorial dual
since, for example, the top of a face $S$ of the Ford domain may be
covered by another face of the Ford domain, and thus the geodesic dual
to $S$ will run through this other face in the universal cover before
meeting $S$.  However, this will not affect our arguments below.

As above, let $S_1$ and $S_2$ denote adjacent faces of the Ford
domain, which are Euclidean hemispheres of radius $R_1$ and $R_2$,
respectively, and whose centers are Euclidean distance $D$ apart.  Let
$S_1'$ and $S_2'$ denote the faces that glue to $S_1$ and $S_2$,
respectively.  So $S_1'$ and $S_2'$ are Euclidean hemispheres of
radius $R_1$ and $R_2$ in the universal cover.

Because $S_1$ is adjacent to $S_2$, we may consider the 2--cell which
is the geometric dual of their intersection.  By assumption, this is
an ideal triangle in the manifold $M$.  One edge of this triangle is
dual to $S_1$ and its paired face $S_1'$.  We take a lift to $\HH^3$
such that this edge runs from infinity straight down the vertical
geodesic with endpoints infinity and the center of $S_1$.  When it
meets $S_1$, it is identified with the corresponding point (at the
center) of $S_1'$, and then runs up the vertical geodesic from the
center of $S_1'$ to infinity.  Another edge is dual to $S_2$ and
$S_2'$, and can be seen in $\HH^3$ similarly.

This triangle will have a third edge, by assumption, dual to a pair of
faces $S_3$ and $S_3'$.  Here $S_3$ will be a sphere adjacent to
$S_1'$, and $S_3'$ will be a sphere adjacent to $S_2'$.

\begin{lemma}
\begin{enumerate}
	\item[(a)] The radius of the spheres $S_3$ and $S_3'$ is $R_1R_2/D$.
	\item[(b)] The distance between the center of $S_3$ and the center
		of $S_1'$ is $R_1^2/D$.
	\item[(c)] The distance between the center of $S_3'$ and the center
of $S_2'$ is $R_2^2/D$.
\end{enumerate}
\label{lemma:ford-circles}
\end{lemma}

\begin{proof}
Consider the universal cover.  The isometry gluing $S_1$ to $S_1'$
takes the point on $\CC$ at the center of $S_1$ to infinity.  It
therefore takes the third edge of the triangle, which lifts to a
geodesic in $\HH^3$ running from the center of $S_1$ to the center of
$S_2$, to a geodesic running from infinity down to the center of
$S_3$.

We may assume without loss of generality that the center of $S_1$ is
$0$ and the center of $S_2$ is $D$.  The isometry taking $S_1$ to
$S_1'$ is an inversion in $S_1$, followed by a Euclidean reflection
and rotation \cite{maskit}.  Since Euclidean reflection and rotation
do not affect radii of hemispheres or distance on $\CC$, the lengths
are given by determining the corresponding lengths under the inversion
in $S_1$.

Note under this inversion, $D$ maps to $R_1^2/D$, proving part (b).  A
symmetric argument, reversing the roles of $S_1$ and $S_2$, gives part
(c).

Finally, to show that the size of the radius is as claimed, consider
the point of intersection of $S_1$ and $S_2$ which lies over the real
line.  It has coordinates $(x, 0, z)$, say.  Since this is a point on
$S_1$, it will be taken to itself under the inversion.  However, this
point is on the edge of the Ford domain where the three faces $S_1$,
$S_2$ and $S_3$ meet.  Thus it will also lie on $S_3$ after the
inversion.  So to find the radius of $S_3$, we only need to determine
the Euclidean distance between this point of intersection $(x, 0, z)$
and the center $(R_1^2/D, 0, 0)$ of $S_3$.

The square of this distance is $x^2 -2xR_1^2/D + z^2$.  Using the fact
that $x^2 + z^2 = R_1^2$ (since $(x,0,z)$ lies on $S_1$) we simplify
this formula to $R_1^2/D^2(D^2-2Dx+R_1^2)$.  Now using the fact that
$(x,0,z)$ lies on $S_2$, we know $x^2 -2Dx + D^2 + z^2 = R_2^2$, or
$R_1^2 -2Dx + D^2 = R_2^2$.  Hence the square of the radius is
$R_1^2R_2^2/D^2$.  
\end{proof}

Finally, we prove a general estimate about the lengths of edges that
are not dual to the Ford domain.

\begin{lemma}\label{lemma:length-non-canonical}
Let $e$ be a geodesic from cusp to cusp in a hyperbolic manifold $M$.
Fix a choice of horoball neighborhoods.  If $e$ is not an edge of the
canonical polyhedral decomposition (with respect to this horoball
neighborhood), then the length of $e$ is at least $\ln(2)$.
\end{lemma}

\begin{proof}
Suppose not.  Suppose there exists a geodesic from cusp to cusp which
is not a canonical edge yet has length less than $\ln(2)$.  Lift to
$\HH^3$.  The geodesic lifts to a geodesic $\gamma$.  Conjugate such
that one endpoint of $\gamma$ is infinity, and such that the
horosphere of height $1$ about infinity projects to the cusp.  Then
the other endpoint of $\gamma$ runs through a horosphere $H$ of
diameter greater than $1/2$.

The set of all points equidistant from $H$ and from the horosphere
about infinity is a hemisphere $S$ of radius at least $1/\sqrt{2}$.
This is not a face of the Ford domain, hence there must be some face
of the Ford domain $F_r$ of radius $r$, say, which overlaps the
highest point of $S$.  Thus $1/\sqrt{2} < r \leq 1$, and the distance
$d$ between the center of $F_r$ and the center of $S$ is at most
$\sqrt{r^2-1/2}$.

On the other hand, there must be a horosphere under the hemisphere
$F_r$ of diameter $r^2$.  The distance $d$ between the center of the
horosphere of diameter $r^2$ and that of diameter $1/2$ is at least
$r/\sqrt{2}$, with equality when the two horospheres are tangent.

Then we have
$$ \frac{r}{\sqrt{2}} \leq d \leq \sqrt{r^2 - \frac{1}{2}}.$$
This is possible only when $r = 1$ and $d=1/\sqrt{2}$.  However, in that case the
highest point of $S$ will not be overlapped by $F_r$.
\end{proof}

\subsection{Parameterization by radii of Ford domain faces}\label{subsec:param}

\begin{figure}%[h]
\centerline{\input{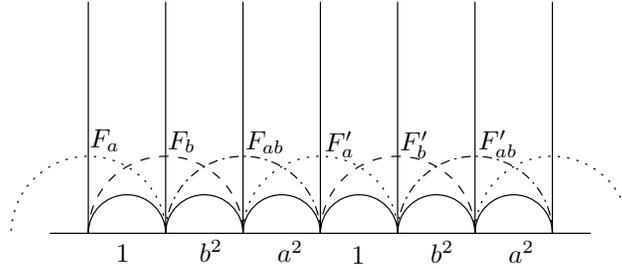}}
\caption{Euclidean distances in the universal cover of a zigzag.  The
translation from the left-most to the right-most edge represents one
meridian in a punctured torus, or two meridians in a $4$--punctured
sphere.}
\label{fig:torus-lengths}
\end{figure}

We now restrict our attention to the case of Farey manifolds. Suppose,
for the moment, that $M$ is a punctured torus bundle.  Consider one
zigzag of $M$; this is a punctured torus $T$.  From the canonical
triangulation on $M$, $T$ inherits a triangulation.  Edges are dual to
faces of the Ford domain of $M$.  Since $T$ is a punctured torus,
there are only three edges in a triangulation of $T$, and two
triangles.  Thus the zigzag of $T$ meets six hemispheres of the Ford domain,
which are identified in pairs.  See Figure \ref{fig:torus-lengths}.

Let $F_a$ be a face of the Ford domain whose radius is largest among
the faces dual to the pleating locus of $T$. (In other words, $F_a$ is
dual to the edge of the pleating that is shortest outside the maximal
cusp.)  Conjugate $\HH^3$ such that the distance between the center of
$F_a$ and the center of the nearest adjacent face of the Ford domain
to the right ($F_b$, say) is $1$.  Let $a$ denote the radius of $F_a$,
$b$ the radius of $F_b$.

By Lemma \ref{lemma:ford-circles}, the other circle of the Ford domain
which is met by $T$ has radius $ab$.  Call this face $F_{ab}$.  By
following the triangulation of a once punctured torus, we see the
Euclidean lengths between centers of horospheres must be as in Figure
\ref{fig:torus-lengths}.  

This parameterization for a pleated punctured torus extends easily to
$4$--punctured spheres.  In an ideal triangulation of a $4$--punctured
sphere $S$, there are six edges and four ideal triangles --- double
the complexity above.  However, recall that we have chosen the cusp
neighborhoods and the canonical triangulation equivariantly.  As a
result, the zigzag of $S$ will look the same when viewed from each
cusp.  When viewed from any puncture of $S$, the zigzag crosses three
faces of the Ford domain, whose radii will be $a$, $b$, and $ab$.

We are interested in the sizes of horospheres at the bottom of each
edge in Figure \ref{fig:torus-lengths}.

\begin{lemma}
Suppose that when we lift to $\HH^3$, the maximal cusp of $M$ lifts to
a horosphere at height $h$, while the zigzag has Euclidean distances
and radii as above.  Then the distances between horospheres along the
edges dual to $F_a$, $F_b$ and $F_{ab}$ are $2\log(h/a)$,
$2\log(h/b)$, and $2\log(h/(ab))$, respectively.

Thus if we conjugate again such that the maximal cusp of $M$ lifts to
a horosphere of height $1$, then we see horospheres of diameter
$a^2/h^2$, $b^2/h^2$, and $a^2b^2/h^2$, respectively.
\label{lemma:horo-sizes}
\end{lemma}

\begin{proof}
Recall that the face $F_a$ is equidistant from the horosphere of
height $h$ about infinity and another horosphere which lies under
$F_a$.  Thus the distance between the face $F_a$ and the horosphere
below it must equal the distance between the face of radius $a$ and
the horosphere of height $h$ above it.  Thus the distance between the
two horospheres is $2\log(h/a)$.

Now, if we conjugate such that the maximal cusp of $M$ lifts to a
horosphere of height $1$, we do not change hyperbolic lengths, so the
distance between horospheres is still $2\log(h/a)$.  But now, if the
diameter of the horosphere centered on $\CC$ is $d$, this implies
$\log(1) - \log(d) = 2\log(h/a)$, or $d= a^2/h^2$.

The argument is the same for horospheres under $F_b$ and $F_{ab}$.
\end{proof}

By Lemma \ref{lemma:horo-sizes}, the largest horosphere has diameter
the maximum of $a^2/h^2$, $b^2/h^2$, and $a^2b^2/h^2$.  But we chose
$F_a$ so that $a$ was the maximum of $a$, $b$, and $ab$.  So the
largest horosphere has diameter $a^2/h^2$.

To improve estimates, we may use the fact that faces of the Ford
domain meet in a certain pattern in the three dimensional manifold $M$
as well as in the surface $S$.  We will need the following lemma about
angles between faces of the Ford domain.  This lemma was first
observed in a slightly different form by Gu{\'e}ritaud \cite[Page
29]{gueritaud:thesis}.

\begin{lemma}\label{lemma:cos-bound}
Let $F_A$, $F_B$, $F_C$ and $F_E$ be faces of the Ford domain
corresponding to a single zigzag, with $F_A$ adjacent to $F_B$, $F_B$
to $F_C$, and $F_C$ to $F_E$.  Suppose also that $F_A$, $F_B$, and
$F_C$ are dual to a canonical tetrahedron.  Denote by $A$ the
Euclidean radius of the hemisphere $F_A$ (which is also the radius of
$F_E$), and denote by $C$ the Euclidean radius of $F_C$.  Denote the
distance between the centers of $F_E$ and $F_C$ by $D$.  Let $\alpha$
denote the angle between the line segments from the center of $F_B$ to
the center of $F_C$, and from the center of $F_B$ to the center of
$F_A$.  Then the angle $\alpha$ satisfies
$$\cos{\alpha} >
 \frac{ A^4 + C^4 + D^4 - 2A^2 D^2 - 2 C^2 D^2}
			{2 A^2 C^2}.$$  
\end{lemma}

\begin{proof}
Note that $\alpha$ is the dihedral angle of a tetrahedron in the
canonical triangulation.  That tetrahedron is dual to the point of
intersection of faces $F_A$, $F_B$, and $F_C$.  The key fact that we
will use is that these three faces must overlap.

Consider the circles given by the points where the spheres of $F_A$
and $F_B$ meet the boundary at infinity.  We will abuse notation and
call these circles $F_A$ and $F_B$.  Consider a third circle
$C(\beta)$ with radius $C$ such that the line between the center of
this circle and the center of $F_B$ makes an angle $\beta$ with the
line between the center of $F_B$ and the center of $F_A$.  When
$\beta=\alpha$, this circle $C(\beta)$ is the circle of $F_C$.  See
Figure \ref{fig:circle-lemma}(a).

However, we want to consider varying $\beta$.  The angle $\beta$ can
lie anywhere in the interval $(0,\pi)$.  For large $\beta$, the
circle $C(\beta)$ may not meet $F_A$.  We can decrease $\beta$ until
these two circles overlap.  Since $F_A$, $F_B$, and $F_C$ are dual to
a tetrahedron of the canonical triangulation, when $\beta= \alpha$,
$C(\alpha)$ and $F_A$ must overlap enough that the interiors of the
regions bounded by these circles and by $F_B$ intersect nontrivially.
Thus $\alpha$ must be strictly less than the value of $\beta$ for
which the three circles meet in a single point.  We will find this
value of $\beta$.  See Figure \ref{fig:circle-lemma}(b).

\begin{figure}
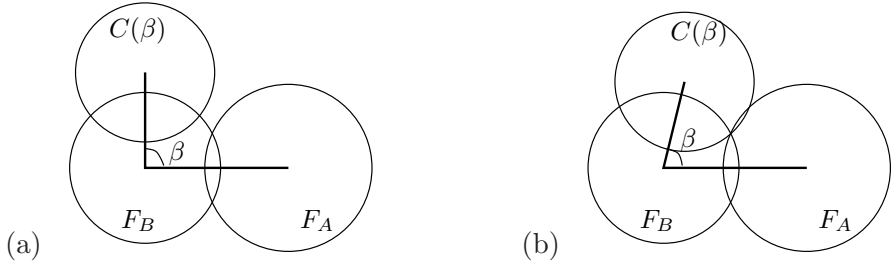

\begin{tabular}{ccccc}
	(a) &
	\input{figures/beta-a.pstex_t}
	& \hspace{.5in} & (b) &
	\input{figures/beta-b.pstex_t}
\end{tabular}
\caption{(a) The circles of $F_A$, $F_B$, $C(\beta)$.  (b) The value of
$\beta$ for which the faces meet in a single point.}
\label{fig:circle-lemma}
\end{figure}

Now, given the distance $D$ and the radii $A$ and $C$, we can compute
all the other distances and radii of the zigzag, using Lemma
\ref{lemma:ford-circles}.  In particular, the radius of $F_B$ is
$AC/D$.  The distance between centers of $F_A$ and $F_B$ is $A^2/D$,
and the distance between the centers of $F_B$ and $F_C$ is $C^2/D$.

Without loss of generality, suppose $F_B$ has center $(0,0)$, and
$F_A$ has center $(A^2/D, 0)$. Here we are writing points in
$\CC$ as points in $\RR^2$.  Then the center of $C(\beta)$ is
$$((C^2/D)\cos{\beta}, (C^2/D)\sin{\beta}).$$

The value of $\beta$ for which the three circles meet in a single
point will be determined as follows.  The circles of $F_A$ and of
$F_B$ intersect in two points which lie on a line $\ell_{AB}$ between
the circles.  Similarly, the circles of $F_B$ and of $C(\beta)$
intersect in two points which lie on a line $\ell_{BC}$.  Notice that
the three circles meet in a single point exactly when the lines
$\ell_{AB}$ and $\ell_{BC}$ intersect in a point which lies on the
circle of $F_B$.  We therefore compute these lines and their
intersection.

The line $\ell_{AB}$ is given by the intersection of the circles
$(x-A^2/D)^2 + y^2 = A^2$, and $x^2 + y^2 = (AC/D)^2$.  This has
equation:
$$ x = \frac{A^2 + C^2 - D^2}{2D}.$$
Similarly, the line $\ell_{BC}$ has equation:
$$ (\cos\beta)x +(\sin\beta)y = \frac{A^2 + C^2 - D^2}{2D}$$
Their intersection is therefore the point
\begin{equation}
\left( \frac{A^2 + C^2 - D^2}{2D}, \frac{A^2 + C^2 - D^2}{2D}
	\left(\frac{1-\cos\beta}{\sin\beta}\right)\right).
\label{eqn:point}
\end{equation}

We want this point to lie on the circle $x^2 + y^2 = (AC/D)^2$.
Plugging the point (\ref{eqn:point}) into the equation of the circle,
we find $\beta$ satisfies
$$ \left( \frac{A^2 + C^2 - D^2}{2D} \right)^2 \left(1 +
\left(\frac{1-\cos\beta}{\sin\beta}\right)^2 \right) =
\frac{A^2 C^2}{D^2},$$
which simplifies to
$$
 \frac{2}{1+ \cos \beta}
= \left(\frac{2 AC}{A^2 + C^2 - D^2}\right)^2.
$$
Thus
$$\cos \beta \: = \:  \frac{(A^2 + C^2 - D^2)^2}{2A^2 C^2} - 1 \: = \: 
\frac{ A^4 + C^4 + D^4 - 2A^2 D^2 - 2 C^2 D^2}{2 A^2 C^2}.$$

%	To simplify notation and the calculation, let $K$ be the right hand
%	side of equation (\ref{eqn:K}).  Then solving for $\cos\beta$, we
%	find:
%	$$K\cos^2\beta - 2\cos\beta + 2-K = 0,$$ or
%	$$\cos\beta = \frac{1}{K} \pm \frac{|1-K|}{K}.$$
%	%
%	Since $\beta \neq 0$, the only possible solution is $\cos\beta = 2/K -
%	1$.

Since $\alpha$ is strictly less than this $\beta$, and $0<\alpha<\pi$,
$\cos\alpha$ must be strictly greater than $\cos\beta$.  
This completes the proof.
%	Thus putting
%	back the value of $K$ and simplifying, we find:
%	$$\cos\alpha >
%	\frac{ A^4 + C^4 + D^4 - 2A^2 D^2 - 2 C^2 D^2}{2 A^2 C^2}$$
%
%	\vspace{-3ex}
\end{proof}

\subsection{Horosphere estimate}\label{subsec:large-horosphere}

We can now show that each zigzag contains a large horosphere.

\begin{prop}\label{prop:horosize}
Let $M$ be a Farey manifold. If $M$ is a $4$--punctured sphere bundle
or two--bridge link complement, denote its meridian length by $\mu$;
if $M$ is a punctured torus bundle, denote its meridian length by
$2\mu$.  (See Definition \ref{def:fibered-meridian}.)  Then every
zigzag in $M$ contains a horosphere of diameter at least $\mu^2/7$.
\end{prop}

\begin{proof}
Let $S$ be a zigzag in $M$.  As at the beginning of
\S\ref{subsec:param}, let $a$ denote the radius of the largest face of
the Ford domain of $S$.  Call this face $F_a$.  Rescale such that the
distance between the center of $F_a$ and the center of the face
directly to its right is $1$.  Call the face to its right $F_b$, and
let $b$ denote the radius of the face $F_b$.  The third face, which we
will call $F_{ab}$, will then have radius $ab$, and have center
distance $b^2$ from $F_b$, and distance $a^2$ from $F_a$, by Lemma
\ref{lemma:ford-circles}. Then the length $\mu$ is equal to $d/h$,
where $d$ is the minimal distance between centers of faces $F_a$ and
$h$ is the height of the maximal cusp in $M$.

By Lemma \ref{lemma:triangle-ineq}, $a$, $b$, and $1$ satisfy the
triangle inequality.  Additionally, because $F_a$ was chosen to have
radius larger than that of $F_b$ and $F_{ab}$, we have the following
inequalities.
\begin{enumerate}
	\item $a\geq b$, and $a\geq ab$, hence $1\geq b$.
	\item $b>-a+1$ and $b>a-1$.
\end{enumerate}

This forces values of $a$ and $b$ to lie within the region shown in
Figure \ref{fig:abregion}.

\begin{figure}[b]
	\includegraphics{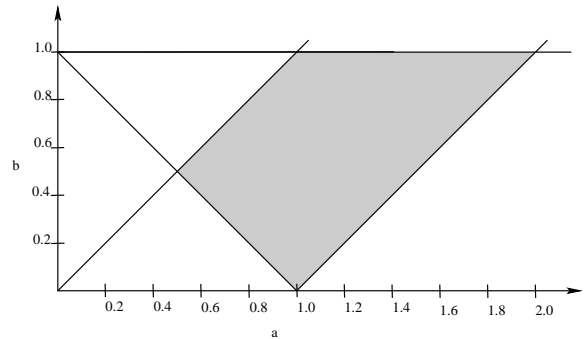}
\caption{The region of allowable values for $a$ and $b$ in Proposition \ref{prop:horosize}.}
\label{fig:abregion}
\end{figure}

Label the angles of the zigzag as follows.  Let $\theta$ denote the
acute angle between the edges of the zigzag of length $a^2$ and $b^2$.
Let $\eta$ denote the acute angle between edges of the zigzag of
length $1$ and $b^2$.  Note this means that the angle between edges of
length $1$ and $a^2$ is $\pi-\theta+\eta$.

By considering orthogonal projections to the edge of length $b^2$, we
find that
\begin{equation}
d^2 = 1+a^4 + b^4 - 2a^2b^2\cos\theta - 2b^2\cos\eta
-2a^2\cos(\pi-\theta+\eta).
\label{eqn:d}
\end{equation}
See Figure \ref{fig:zigzag} for an example.  Note in Figure
\ref{fig:zigzag}, the angles $\theta$ and $\eta$ correspond to angles
of tetrahedra in the canonical decomposition.  Because we chose $F_a$
to be the largest face, this will not necessarily be the case, but two
of the three angles $\theta$, $\eta$, $\pi-\theta+\eta$ will be
canonical (or, if $2\pi-(\pi-\theta+\eta)$ happens to be acute rather
than $\pi-\theta+\eta$, then exactly two of the three angles $\theta$,
$\eta$, and $2\pi-(\pi-\theta+\eta)$ will be canonical).

\begin{figure}[h]
\centerline{\input{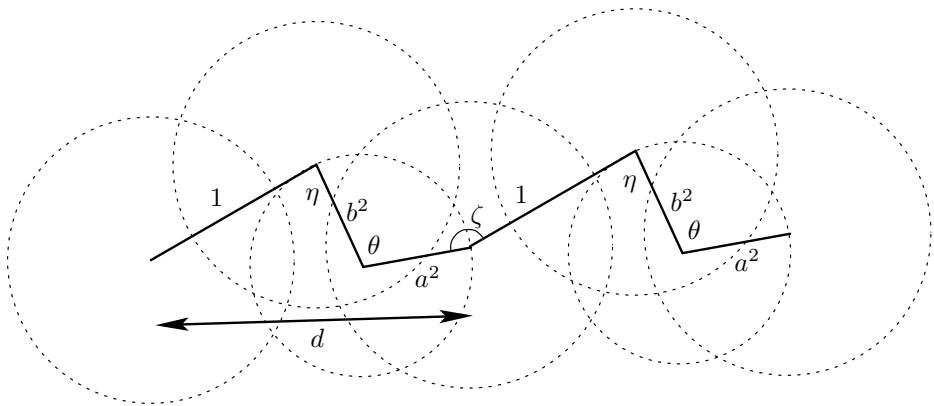}}
\caption{A zigzag.  Here the dotted circles correspond to faces of the
	Ford domain.  Recall the corners of the zigzag are geometric duals
	of these faces.  Reading left to right following the zigzag, the
	faces of the Ford domain have radius $a$, $b$, $ab$, $a$, $b$, $ab$,
	and $a$.  Here $\eta$ and $\theta$ are angles of canonical
	tetrahedra.  The angle $\pi-\theta+\eta$ is not an angle of
	a canonical tetrahedron.}
\label{fig:zigzag}
\end{figure}

By Lemma \ref{lemma:horo-sizes}, and because we chose the face $F_a$
to be largest, the largest horosphere in the zigzag $S$ has
diameter $a^2/h^2$.

Write:
$$\frac{a^2}{h^2} = \frac{a^2\mu^2}{d^2} =
	\mu^2 \, \frac{a^2}{d^2}.$$
We minimize the quantity $a^2/d^2$.

Note that if $\theta$ is an angle of a tetrahedron in the canonical
polyhedral decomposition of $M$, then by Lemma \ref{lemma:cos-bound},
$\cos\theta$ satisfies:
$$\cos\theta > \frac{1+a^4+b^4 - 2a^2-2b^2}{2a^2b^2}.$$
Similarly, Lemma \ref{lemma:cos-bound} implies that if $\eta$ is an
angle of a tetrahedron in the canonical polyhedral decomposition of
$M$, then $\cos\eta$ satisfies:
$$\cos\eta > \frac{1+a^4+b^4 - 2a^2b^2-2a^2}{2b^2},$$
and if $\pi-\theta+\eta$ (or $2\pi-(\pi-\theta+\eta)$) is an angle of
a tetrahedron in the canonical polyhedral decomposition, then
$\cos(\pi-\theta+\eta)$ satisfies:
$$\cos(\pi-\theta+\eta) > \frac{1+a^4+b^4 -2a^2b^2-2b^2}{2a^2}.$$

Two of the three will be canonical.  The third will not, since all
three angles cannot be canonical at the same time.  However, we know
the cosine in that case will be at least $-1$.  Hence combining the
cosine inequalities above with the formula for $d^2$ of (\ref{eqn:d}),
we will have one of the three inequalities:
\begin{enumerate}
	\item If $\theta$ and $\eta$ are canonical:
$$ \frac{a^2}{d^2} >
		\frac{a^2}{6a^2+2a^2b^2+2b^2-1-a^4-b^4} =: f_1(a,b).$$ 
	\item If $\theta$ and $\pi-\theta+\eta$ are canonical:
$$\frac{a^2}{d^2} >
		\frac{a^2}{6b^2 + 2a^2b^2 + 2a^2 -1-a^4-b^4}=:f_2(a,b).$$
	\item If $\eta$ and $\pi-\theta+\eta$ are canonical:
$$\frac{a^2}{d^2} >
		\frac{a^2}{6a^2b^2+2b^2+2a^2 -1-a^4-b^4} =:f_3(a,b).$$
\end{enumerate}

To complete the proof, we minimize all three of these functions in the
region of Figure \ref{fig:abregion}.  This is a calculus problem.

For each $f_j(a,b)$, $j=1, 2, 3$, we find the only critical point of
$f_j$ in the region of Figure \ref{fig:abregion} is the point $a=1,
b=0$.  For all positive $a$, the function $f_j$ is decreasing on the
line $b=a$, increasing on the line $b=1$, increasing or constant on
$b=-a+1$, and decreasing or constant on $b=a-1$.  This implies that
$f_j$ takes its minimum value in the region at the point $a=1, b=1$.

At this value, $f_j(1,1) = 1/7$.  Hence $a^2/h^2 \geq \mu^2/7$.  
\end{proof}

\begin{remark}
The proof of Proposition \ref{prop:horosize} does not require the
zigzag $S$ to be embedded.  In other words, the proposition applies
even to the terminal pleated surfaces $S_1$ and $S_c$ that are folded
in the construction of a 2--bridge link.  When $S$ is a folded surface
$S_1$ or $S_c$, one of the angles $\theta$, $\eta$, or $(\pi - \theta
+ \eta)$ is actually $0$, hence its cosine is even larger than
claimed, which only improves the estimate.
\end{remark}

Proposition \ref{prop:horosize} should be compared to previous work of
J{\o}rgensen \cite[Lemma 4.3]{jorgensen}, which was carefully written
down by Akiyoshi, Sakuma, Wada, and Yamashita \cite[Lemma
8.1.1]{aswy-book}.  After adjusting for slightly different choices of
normalization, J{\o}rgensen's Lemma 4.3 can be summarized as saying
that, whenever a zizag of a quasifuchsian punctured torus group is dual 
to six faces of the Ford domain, one of those faces has radius at least 
$$\mu/(4+2\sqrt{5}).$$
It follows from Minksy's classification of punctured torus groups
\cite{minsky}, that given a punctured torus bundle $M$, the Kleinian
subgroup of $\pi_1(M)$ that corresponds to the fiber can be obtained
as a geometric limit of quasifuchsian groups.  As a result,
J{\o}rgensen's estimate extends to punctured torus bundles.  Because a
Ford domain face of radius $r$ corresponds to a horosphere of diameter
$r^2$, J{\o}rgensen's Lemma 4.3 implies that every zigzag in a
punctured torus bundle contains a horosphere of diameter at least
$$\frac{\mu^2}{(4+2\sqrt{5})^2} \: \approx \: \frac{\mu^2}{71.777}.$$
Proposition \ref{prop:horosize}, which is proved by a direct geometric
argument without reference to quasifuchsian groups, improves this
estimate by a factor of about $10.25$.

This improvement becomes highly significant in Section
\ref{sec:cusp-area}.  In Theorem \ref{thm:ptorus-area}, we estimate
the area of a maximal cusp by packing the horospherical torus with
disjoint disks that are shadows of large horospheres.  As a result, a
10--fold increase in the estimate for the diameter of a horosphere
turns into a 100--fold increase in the estimate for the area of its
shadow.
Since our applications in Sections \ref{sec:volume3braids} and
\ref{sec:jones3braids} rely on these explicit estimates for cusp area,
the 100--fold improvement becomes particularly important for
applications.  

\subsection{The length of a meridian}

To make the estimate of Proposition \ref{prop:horosize} independent of
$\mu$, we prove a bound on the value of $\mu$.  We note that the
following lemma is the \emph{only} result in this section that does
not apply to all Farey manifolds: it fails for 2--bridge links.

\begin{lemma}
In an equivariant cusp of a 4--punctured sphere bundle $M$, $\mu \geq
\sqrt{2}$.
\label{lemma:meridian}
\end{lemma}

\begin{proof}
Lift the hyperbolic structure on $M$ to $\HH^3$.  The cusps lift to
collections of horoballs.  Conjugate such that the horoball about
infinity of height $1$ projects to a cusp neighborhood.  Since we took
a maximal cusp neighborhood of $M$, that is, since we expanded cusps
until they bumped, there must be some full--sized horoball $H$
projecting to a cusp of $M$, tangent to the horoball of height $1$
about infinity.

There is an isometry of $\HH^3$ corresponding to the slope of length
$\mu$ which is a covering transformation of $M$.  It takes $H$ to
another full--sized horoball $H'$.  The Euclidean distance between $H$
and $H'$ is the length $\mu$.

Consider the geodesic $\gamma$ running from the center of $H$ to the
center of $H'$.  This projects to a geodesic in $M$ running from one
puncture of the fiber back to the same puncture.  Note that under the
equivariant cusp expansion, any canonical edge runs between two
distinct punctures of the fiber.  Hence by Lemma
\ref{lemma:length-non-canonical}, the length of the portion of
$\gamma$ outside $H$ and $H'$ is at least $\ln(2)$.

Now, recall the following formula for lengths along ``right angled
hexagons'' (see, for example, \cite[Lemma 3.4]{gmm:smallest-cusped}).  
Let $H_\infty$ denote the horosphere about infinity, and
let $H_p$ and $H_q$ be disjoint horospheres not equal to $H_\infty$,
centered over $p$ and $q$ in $\CC$, respectively.  Denote by $d_p$ the
hyperbolic distance between $H_p$ and $H_\infty$, by $d_q$ the
hyperbolic distance between $H_q$ and $H_\infty$, and by $d_r$ the
hyperbolic distance between $H_p$ and $H_q$.  Then the Euclidean
distance between $p$ and $q$ is given by
\begin{equation}
d(p,q) = \exp( (d_r - (d_p+d_q))/2 ).
\label{eqn:right-hex}
\end{equation}

In our case, $d_p = d_q = 0$, since the corresponding horospheres ($H$
and $H'$) are tangent to $H_\infty$, and $d_r$ is at least $\ln(2)$.
So $\mu = d(p,q)$ is at least $\sqrt{2}$.
\end{proof}

By commensurability, the meridian in a punctured torus bundle has
length $2 \mu \geq 2 \sqrt{2}$.  In the setting of two-bridge links,
on the other hand, Lemma \ref{lemma:meridian} fails because a meridian
of the link is spanned by a single edge of the canonical
triangulation. For two--bridge links, the best available estimate is
Adams's result that $\mu \geq \sqrt[4]{2}$, which works for all links
except the figure--8 and $5_2$ knots \cite{adams:waist2}.

\section{Cusp area estimates}\label{sec:cusp-area}

In this section, we apply the results of Section \ref{sec:horosphere}
to prove quantitative estimates on the cusp area of Farey
manifolds. For most of the section, we shall focus on punctured torus
bundles.  At the end of the section, we will generalize these results
to $4$--punctured sphere bundles and 2--bridge links.

\subsection{Punctured--torus bundles}

We shall prove the following result:

\begin{theorem}\label{thm:ptorus-area}
Let $M$ be a punctured--torus bundle with monodromy
$$\Omega = \pm R^{p_1}L^{q_1} \cdots R^{p_s}L^{q_s}.$$
Let $C$ be a maximal horoball neighborhood about the cusp of $M$.
Then
$$0.1885 \, s \: \approx \: \frac{16 \sqrt{3}}{147}  \, s \: \leq \:  \area(\bdy C) \:<\: 2
\sqrt{3} \, \frac{v_8}{v_3} \, s \: \approx \: 12.505 \, s.$$
Furthermore, if $\gamma$ is any simple closed curve on $\bdy C$ that
is transverse to the fibers, then its length $\ell(\gamma)$ satisfies
$$ \ell(\gamma) \: \geq \: \frac{4\sqrt{6}}{147} \,  s.$$
\end{theorem}

\begin{remark} 
Extensive numerical experiments support the conjecture that
$\area(\bdy C)/s$ is monotonic under the operation of adding more
letters to existing syllables of the monodromy word $\Omega$.  (It is
not hard to show using the method of angled triangulations
\cite{gf:two-bridge} that the volume of $M$ behaves in a similarly
monotonic fashion.)  This conjecture would imply that the quantity
$\area(\bdy C)/s$ is lowest when all syllables have length 1 and $M$
is a cover of the figure--8 knot complement, while $\area(\bdy C)/s$
approaches its upper bound as the syllable lengths approach $\infty$
and the geometry of $M$ converges to a cover of the Borromean rings.
The cusp area of the figure--8 knot complement is $2\sqrt{3}$,
and the cusp area of one component of the Borromean rings is $8$. Thus, if the monotonicity
conjecture is correct, it would follow that
$$2\sqrt{3} \, s \: \leq \: \area(\bdy C) \:<\: 8 \, s.$$
Compared to the values in Theorem \ref{thm:ptorus-area}, this represents a modest improvement of the upper bound but a dramatic improvement of the lower bound.
\end{remark}

The main idea of the proof of Theorem \ref{thm:ptorus-area} is to pack
the horospherical torus $\bdy C$ with disjoint disks that are shadows
of large horospheres.  Recall from Definition
\ref{def:fibered-meridian} that we denote the length of a meridian in
the maximal cusp of a punctured--torus bundle by $2\mu$.  By
Proposition \ref{prop:horosize}, every zigzag on $\bdy C$ will contain
two horospheres of diameter at least $\mu^2/7$, corresponding to the
two endpoints of the same edge of the zigzag.  When we project one of
these horospheres to $\bdy C$, we obtain a disk whose radius is at
least $\mu^2/14$.

To turn this into an effective estimate on the area of $\bdy C$, we
need to employ a somewhat subtle procedure for choosing which
horospheres to count and which ones to discard. We choose the
horospheres in the following manner:

\begin{enumerate}
\item Let $E$ be the set of all edges of $M$ whose length outside the
maximal cusp is at most $\ln(7/\mu^2)$. These are exactly the edges
that lead to horospheres of diameter $ \geq \mu^2/7$. Thus, by
Proposition \ref{prop:horosize}, every pleated surface in $M$ contains
an edge in $E$.

\item Order the letters of the monodromy word $\Omega$: $\alpha_1,
\ldots, \alpha_m$. Recall, from Section \ref{sec:ptorus}, that each
$\alpha_i$ corresponds to a pleated surface $T_{\alpha_i}$.

\item Find the smallest index $i$ such that all three edges in the
pleating of $T_{\alpha_i}$ belong to $E$. (It is possible that such an
$i$ does not exist.) If such a $T_{\alpha_i}$ occurs, remove the
longest of the those three edges from $E$, breaking ties at random.

\item Repeat step (3) inductively. In the end, the set $E$ will
contain at most two edges from each pleated surface.
\end{enumerate}

At the end of step (4), if a pleated surface $T$ contains one edge of
$E$, that edge is the shortest in $T$. If $T$ contains two edges of
$E$, they are the two shortest edges in $T$.

\begin{lemma}\label{lemma:edge-count}
The set $E$, constructed as above, contains at least $2s/3$ distinct
edges.
\end{lemma}

\begin{proof}
By Proposition \ref{prop:horosize}, every pleated surface in $M$
contains an edge whose length is at most $\ln(7/\mu^2)$.  Thus, at the
end of step (1) in the selection procedure above, the set $E$
contained at least one edge from every pleated surface.  Now, observe
that two different pleated surfaces $T_\alpha$ and $T_\beta$,
corresponding to letters $\alpha$ and $\beta$ in $\Omega$, will share
an edge if and only if the corresponding triangles in the Farey graph
share a vertex.  As Figure \ref{fig:farey} illustrates, this can only
happen if the letters $\alpha$ and $\beta$ come from the same
syllable, neighboring syllables, or syllables that share a neighbor.
Therefore, at the end of step (1), the set $E$ contained at least one
edge for every consecutive string of three syllables, hence at least
$2s/3$ distinct edges in total.

\begin{figure}
\begin{center}
  \input{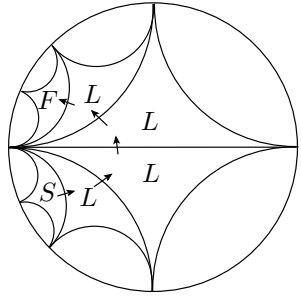}
\end{center}
\caption{The only possible words between two triangles that share a
	vertex are $S\,RR\cdots R\,F$ and $S\,LL\cdots L\,F$, where $S$ and
	$F$ (start and finish) can both be either $L$ or $R$.  }
\label{fig:farey}
\end{figure}

Now, consider what happens when we begin removing edges in step (3).
Suppose that all three edges in a pleated surface $T_\alpha$ belong to
$E$.  Then, just as above, for any pleated surface $T_\beta$ that
shares an edge with $T_\alpha$, the letters $\alpha$ and $\beta$ must
come from the same syllable, neighboring syllables, or syllables that
share a neighbor.  There are at most five such syllables altogether
(the syllable containing $\alpha$, plus two on each side).  Thus,
after we remove the longest edge of $T_\alpha$ from $E$, the set $E$
still contains two edges from a string of five consecutive syllables.

At the end of the selection procedure, every pleated surface in $M$
belongs either to a string of 3 syllables containing at least one edge
of $E$, or to a string of 5 syllables containing at least two edges of
$E$.  In either scenario, there are at least $2s/3$ edges belonging to
$E$.
\end{proof}

\begin{lemma}\label{lemma:disjoint-disks}
Let $M$ be a punctured--torus bundle with monodromy 
$$\Omega = \pm R^{p_1}L^{q_1} \cdots R^{p_s}L^{q_s},$$
Then the maximal cusp boundary $\bdy C$ contains $4s/3$ disjoint
disks, each of radius at least
$$\min \left\{ \frac{1}{4}, \, \frac{\sqrt{2} \, \mu^2}{14}
\right\}.$$
\end{lemma}

\begin{proof}
Consider the edge set $E$, as above.  By Lemma \ref{lemma:edge-count},
$E$ contains at least $2s/3$ edges of length at most $\ln(7/\mu^2)$.
Now, lift everything to the universal cover $\HH^3$, in such a way
that $\bdy C$ lifts to the horizontal plane at height $1$. In a single
fundamental domain for $\bdy C$, each edge $e \in E$ corresponds to
two horospheres: one horosphere for each endpoint of $e$.  Hence $\bdy
C$ contains $4s/3$ shadows of disjoint horospheres, each of which has
radius at least $\mu^2/14$.  If two disjoint horospheres have the same
size, then they also have disjoint projections.  Thus, by shrinking
each horosphere to radius $\mu^2/14$, we conclude that $\bdy C$
contains $4s/3$ disjoint disks of radius $\mu^2/14$.

Next, we claim that the disks on $\bdy C$ can be enlarged considerably
while staying disjoint.  Let $x$ and $y$ be the centers of two of
these disks.  In other words, $x \in \bdy C \cap e_i$ and $y \in \bdy
C \cap e_j$ for some $e_i, e_j \in E$.  The two edges $e_i, e_j$ lead
to horospheres $H_i$ and $H_j$. Let $f$ be the geodesic that connects
$H_i$ directly to $H_j$.  Consider the length of $f$ outside $H_i$ and
$H_j$.  There are two cases:

\smallskip
\underline{Case 1: the length of $f$ is at least $\ln(2)$.}  In this
case, the midpoint of $f$ lies at distance at least $\ln(2)/2$ from
both $H_i$ and $H_j$.  If we apply an isometry $I$ that sends $\bdy
H_i$ to the horosphere at Euclidean height $1$, the midpoint of $I(f)$
will lie at height at most $1/\sqrt{2}$.  In other words, the
horosphere $I(H_i)$ can be expanded by a factor of $\sqrt{2}$ without
hitting the midpoint of $I(f)$, and similarly for $H_j$.  Of course,
this still holds true before applying the isometry $I$: each of $H_i$
and $H_j$ can be expanded by a factor of $\sqrt{2}$ while staying
disjoint from each other.  Since each of $H_i$ and $H_j$ has radius at
least $\mu^2/14$, the disks of radius $\sqrt{2}\mu^2/14$ centered at
$x$ and $y$ in $\bdy C$ are disjoint from each other.

\smallskip
\underline{Case 2: $f$ is shorter than $\ln(2)$.}  Then, by Lemma
\ref{lemma:length-non-canonical}, $f$ must be an edge of the canonical
triangulation.  Since $e_i$ and $e_j$ are also edges of the canonical
triangulation, these three edges bound an ideal triangle contained in
some pleated surface $T_\alpha$.  Now, recall that at the end of our
selection procedure for the set $E$, if two distinct edges of
$T_\alpha$ belong to $E$, then they are the shortest edges in
$T_\alpha$.  Thus both $e_i$ and $e_j$ are shorter than $\ln(2)$.
Since the edges $e_i$ and $e_j$ are already vertical in $\HH^3$ and
meet the cusp at Euclidean height $1$, the horospheres $H_i$ and $H_j$
must have diameter at least $1/2$.  Thus $H_i$ and $H_j$ project to
disjoint disks of radius at least $1/4$ centered at $x$ and $y$ on
$\bdy C$.

\smallskip
In every case, the points $x,y \in \bdy C$ are the centers of disjoint
disks of radius at least
$$\min \left\{ \frac{1}{4}, \, \frac{\sqrt{2} \, \mu^2}{14}
\right\}.$$
There are $4s/3$ such disks, completing the proof.
\end{proof}

We may now estimate the area of $\bdy C$.

\begin{lemma}\label{lemma:cusp-area-mu}
Let $M$ be a punctured--torus bundle with monodromy
$$\Omega = \pm R^{p_1}L^{q_1} \cdots R^{p_s}L^{q_s}.$$
Let $C$ be a maximal horoball neighborhood about the cusp of $M$.
Then
$$\frac{16 \sqrt{3}}{147} \, s \; \leq \; \sqrt{3}\, s \, \min\left\{
\frac{1}{6}, \, \frac{4 \mu^4}{147} \right\} \; \leq \; \area(\bdy C)
\;<\; 2 \sqrt{3} \,\, \frac{v_8}{v_3} \, s.$$
\end{lemma}

\begin{proof}
There are three inequalities in the statement, and we consider them in
turn.

\smallskip
\underline{First inequality.}  This follows immediately from Lemma
\ref{lemma:meridian}, which gives $\mu \geq \sqrt{2}$.  Note that with
our definition of $\mu$ (see Definition \ref{def:fibered-meridian}),
the conclusion of Lemma \ref{lemma:meridian} transfers perfectly from
$4$--punctured sphere bundles to punctured torus bundles.

\smallskip
\underline{Second inequality.}  By Lemma \ref{lemma:disjoint-disks},
$\bdy C$ contains $4s/3$ disjoint disks of equal radius, whose total
area is at least
$$ \frac{4s}{3} \cdot \pi \min\left\{\frac{1}{16}, \, \frac{\mu^4}{98}
\right\}.$$
Now, a classical result (see e.g. \cite[Theorem 1]{boroczky}) states
that a packing of the plane by circles of equal size has density at
most $\pi/(2 \sqrt{3})$.  This gives the desired inequality.

\smallskip
\underline{Third inequality.}  A result of Agol gives that $\vol(M) <
2v_8 \, s$ (see \cite[Theorem B.1]{gf:two-bridge} for a direct proof).
Also, a horosphere packing theorem of B{\"o}r{\"o}czky \cite[Theorem
4]{boroczky} states that a maximal cusp in a hyperbolic $3$--manifold
contains at most $\sqrt{3}/(2v_3)$ of the volume of $M$.  Putting
these results together gives
$$  \vol(C) \:<\: \sqrt{3} \,\, \frac{v_8}{v_3} \, s, \quad
\mbox{hence} \quad 
\area(\bdy C) \:<\: 2 \sqrt{3} \,\, \frac{v_8}{v_3} \, s.$$

\vspace{-3ex}
\end{proof}

\begin{remark}
In the proof of Lemma \ref{lemma:disjoint-disks}, we also showed that
$\bdy C$ contains $4s/3$ disjoint disks of radius $\mu^2/14$.
Plugging this estimate into the proof of Lemma
\ref{lemma:cusp-area-mu} gives
$$ \area(\bdy C) \: \geq  \: \frac{2 \sqrt{3} \, \mu^4}{147} \, s.$$
This statement, although apparently weaker than Lemma
\ref{lemma:cusp-area-mu}, will prove useful for estimating the lengths
of slopes on $\bdy C$.
\end{remark}

\begin{lemma}\label{lemma:ptorus-slope-length}
Let $M$ be a punctured--torus bundle with monodromy
$$\Omega = \pm R^{p_1}L^{q_1} \cdots R^{p_s}L^{q_s}.$$
Let $C$ be a maximal horoball neighborhood about the cusp of $M$.  If
$\gamma$ is any simple closed curve on $\bdy C$ that is transverse to
the fibers, $ \ell(\gamma) \geq 4\sqrt{6} \, s/147$.
\end{lemma}

\begin{proof}
Define the \emph{height} of the cusp to be $h := \area(\bdy C)/2\mu$.
Then $\ell(\gamma) \geq h$.  Note that by Lemma \ref{lemma:meridian},
$\mu \geq \sqrt{2}$.  Also, since a maximal horocycle in a punctured
torus has length at most $6$, it follows that $\mu \leq 3$.  We
consider three possibilities for the values of $\mu$ in the range
$[\sqrt{2}, 3]$.

If $\sqrt{2} \leq \mu \leq \sqrt{7}/ 2^{3/4}$, then $4 \mu^4/147 \leq
1/6$.  Thus, by Lemma \ref{lemma:cusp-area-mu},
$$ \area(\bdy C) \, \geq \, \frac{4 \sqrt{3} \, \mu^4}{147} \, s,
 \quad \mbox{hence} \quad \ell(\gamma) \, \geq \, \frac{2\sqrt{3} \,
 \mu^3}{147} \, s \, \geq \, \frac{4\sqrt{6} \,s}{147} \, \approx \,
 0.066652 \, s.$$

If $\sqrt{7}/ 2^{3/4} \leq \mu \leq 2$, then $1/6 \leq 4
\mu^4/147$. Thus, by Lemma \ref{lemma:cusp-area-mu},
$$ \area(\bdy C) \, \geq \, \frac{\sqrt{3} \, s}{6}, \quad
 \mbox{hence} \quad \ell(\gamma) \, \geq \, \frac{\sqrt{3} \, s}{12
 \mu} \, \geq \, \frac{\sqrt{3} \, s}{24} \, \approx \, 0.072168 \,
 s.$$

If $2 \leq \mu \leq 3$, then by the remark following Lemma
\ref{lemma:cusp-area-mu},
$$ \area(\bdy C) \, \geq \, \frac{ 2\sqrt{3} \, \mu^4}{147} \, s,
 \quad \mbox{hence} \quad \ell(\gamma) \, \geq \, \frac{\sqrt{3} \,
 \mu^3}{147} \, s \, \geq \, \frac{8 \sqrt{3} \,s}{147} \, \approx \,
 0.094261 \, s.$$

Therefore, for all possible values of $\mu$, we have $ \ell(\gamma)
\geq 4\sqrt{6} \, s/147$.
\end{proof}

\noindent Lemmas \ref{lemma:cusp-area-mu} and
\ref{lemma:ptorus-slope-length} complete the proof of Theorem
\ref{thm:ptorus-area}.
% \end{proof}

Combining the results of Theorem \ref{thm:ptorus-area} with our work
in \cite[Theorem 1.1]{fkp-07}, we obtain the following immediate
corollary for volumes of Dehn fillings of punctured torus bundles.

\begin{corollary}\label{cor:ptorus-volume}
Let $M$ be a punctured--torus bundle with monodromy of length $s>94$.
Let $C$ be a maximal horoball neighborhood about the cusp of $M$.  For
any simple closed curve $\gamma$ on $\bdy C$ that is transverse to the
fibers, let $M_{\gamma}$ denote the 3--manifold obtained from $M$ by
Dehn filling $\bdy C$ along $\gamma$.  Then $M_{\gamma}$ is
hyperbolic, and
$$ \left(1-\dfrac{7203 \, \pi^2 }{8 \, s^2} \right)^{3/2} 2 v_3 \, s\:
 \leq \: \vol(M_{\gamma}) \: < \: 2v_8 \, s,
$$
where $v_3 = 1.0149...$ is the volume of a regular ideal tetrahedron
and $v_8 = 3.6638...$ is the volume of a regular ideal octahedron.
\end{corollary}

\begin{proof}
By Theorem \ref{thm:ptorus-area}, the slope length of $\gamma$ will be
at least $2\pi$ when $s\geq 95$.  For such slopes, by \cite[Theorem
1.1]{fkp-07} we know the volume of the manifold obtained by Dehn
filling along the slope of length $\ell(\gamma)$ is at least
$$\vol(M_{\gamma}) \geq \left( 1-\left(\frac{2\pi}{\ell(\gamma)}
\right)^2 \right)^{3/2}\vol(M).$$
Hence, using the volume bound for such manifolds given by
\cite[Theorem B.1]{gf:two-bridge}, and the estimate on slope length of
Theorem \ref{thm:ptorus-area}, we have
$$\vol(M_{\gamma}) \geq \left(1-\dfrac{7203\,\pi^2}{8\, s^2}\right)^{3/2}
2 v_3 \, s.$$

For the upper bound, recall that volume only decreases under Dehn
filling \cite{thurston:notes}, and so the result follows immediately
from \cite[Theorem B.1]{gf:two-bridge}.
\end{proof}

\subsection{4--punctured sphere bundles}

\begin{theorem}\label{thm:4ps-area}
Let $N$ be a $4$--punctured sphere bundle with monodromy
$$\Omega = R^{p_1}L^{q_1} \cdots R^{p_s}L^{q_s},$$
and with the property that the monodromy fixes one preferred boundary
circle of the $4$--holed sphere.  Let $D$ be the maximal horoball
neighborhood of the cusp corresponding to this preferred puncture, and
let $\gamma$ be any simple closed curve on $\bdy D$ that is transverse
to the fibers.  Then
$$\area(\bdy D) \: \geq \: \frac{16 \sqrt{3}}{147} \, s \quad \mbox{and}
\quad \ell(\gamma) \: \geq \: \frac{8\sqrt{3}}{147} \, s.$$
\end{theorem}

\begin{proof}
As described in Section \ref{sec:4ps}, the $4$--punctured sphere
bundle $N$ is commensurable to a punctured--torus bundle $M$ with the
same monodromy $\Omega$.  (The common cover is a $4$--punctured torus
bundle $P$.)  Let $C$ be the maximal cusp neighborhood in $M$. Then,
by lifting $C$ to a cusp neighborhood in $P$ and projecting down to
$N$, we obtain a maximal \emph{equivariant} neighborhood of the cusps
of $N$.

Let $B$ be the cusp neighborhood of the preferred puncture in the
equivariant expansion of the cusps of $N$.  Because the cusp
neighborhood $C \subset M$ lifts to 4 distinct cusps in $P$, and one
of those cusps double--covers $B$, Theorem \ref{thm:ptorus-area}
implies that
$$ 2 \cdot \area(\bdy B) \: = \: \area(\bdy C) \: \geq \: \frac{16
\sqrt{3}}{147} \, s.$$

Observe that in the canonical triangulation of $N$, every edge lies in
a pleated fiber, and connects two distinct punctures of the
$4$--punctured sphere.  Thus no edge of the canonical triangulation
has both endpoints inside $B$.  By Lemma
\ref{lemma:length-non-canonical}, this means that the shortest arc
from $B$ to $B$ has length at least $\ln(2)$, and we may expand $B$ by
a factor of at least $\sqrt{2}$ before it bumps into itself.
Therefore, every linear measurement on $\bdy D$ is at least a factor
of $\sqrt{2}$ greater than on $\bdy B$, and
$$ \area(\bdy D) \: \geq \: 2 \cdot \area(\bdy B) \: = \: \area(\bdy
C) \: \geq \: \frac{16 \sqrt{3}}{147} \, s.$$
By the same argument, every simple closed curve on $\bdy D$ is at
least a factor of $\sqrt{2}$ longer than the corresponding loop on
$\bdy B$.  Thus, if $\gamma$ is transverse to the fibers of $N$, $
\ell(\gamma) \geq 8\sqrt{3} \, s/147$.
\end{proof}

We remark that by Theorem \ref{thm:4ps-area}, an analogue of Corollary
\ref{cor:ptorus-volume} also holds for fillings of 4--punctured sphere
bundles.  One important class of manifolds obtained by Dehn filling
(one cusp of) a $4$--punctured sphere bundle is the class of closed
3--braids in $S^3$.  We shall focus on these manifolds below, in
Section \ref{sec:volume3braids}.

\subsection{2--bridge links}

\begin{theorem}\label{thm:2bridge}
Let $K$ be a 2--bridge link in $S^3$, whose reduced alternating
diagram has twist number $t$.  Let $C$ be a maximal neighborhood about
the cusps of $S^3 \setminus K$, in which the two cusps have equal
volume if $K$ has two components.  Then
$$\frac{8\sqrt{3}}{147}\,(t-1) \: \leq \: \area(\bdy C) \: <\:
2\sqrt{3}\,\, \frac{v_8}{v_3}\,(t-1).$$
Furthermore, if $K$ is a knot, let $\gamma$ be any non-trivial arc
that starts on a meridian and comes back to the same meridian (for
example, a non-meridional simple closed curve). Then its length
satisfies
$$\ell(\gamma) \: \geq \: \frac{4\sqrt{6\sqrt{2}}}{147}\, (t-1).$$
\end{theorem}

\begin{proof}
Let $\mu$ denote the length of a meridian of $K$ on $\bdy C$.  By
Proposition \ref{prop:horosize}, every pleated surface $S_i$ in $S^3
\setminus K$ contains at least one edge of length at most
$\ln(\mu^2/7)$.  Furthermore, opposite edges in $S_i$ have the same
length, because the geometry of each pleated surface is preserved by
the full symmetry group of its triangulation.  Thus, if $S_i$ is
embedded in $S^3 \setminus K$, it contains \emph{two} edges of length
at most $\ln(\mu^2/7)$.  The only pleated surfaces that are not
embedded are the folded surfaces $S_1$ or $S_c$ at the ends of the
product region of $K$; each of these surfaces will contain at least
one short edge. (See \cite[Figure 19]{gf:two-bridge} for a description
of how surfaces are folded in the construction of a 2--bridge link.)

Now, we retrace the proof of Theorem \ref{thm:ptorus-area}.  We
construct the set of short edges $E$ exactly as above, except that we
are now counting \emph{pairs} of edges.  Thus, if all three pairs of
edges in a pleated surface are initially part of $E$, we remove the
longest pair.

By the same argument as in Lemma \ref{lemma:edge-count}, $E$ contains
at least $t/3$ distinct edge pairs.  The two paired edges on a pleated
surface $S$ will be distinct unless $S$ is $S_1$ or $S_c$.  Thus, if
both $S_1$ and $S_c$ contribute edges to $E$, the minimum possible
number of edges (not pairs) is $2(t-1)/3$.

The proof of Lemma \ref{lemma:disjoint-disks} goes through without
modification.  As a result, $\bdy C$ contains $4(t-1)/3$ disjoint
disks, of radius at least
$$r \geq \min \left\{ \frac{1}{4}, \, \frac{\sqrt{2} \, \mu^2}{14}
\right\}, \qquad \mbox{hence total area} \quad \geq \frac{4(t-1)}{3}
\cdot \pi \min\left\{\frac{1}{16}, \, \frac{\mu^4}{98} \right\}.$$
Dividing by the maximal density $\pi/2\sqrt{3}$ of a circle packing in
the plane gives
\begin{equation}\label{eqn:2bridge-area}
 \area(\bdy C) \geq 8(t-1) \, \frac{\sqrt{3}}{3} \,
\min\left\{\frac{1}{16}, \, \frac{\mu^4}{98} \right\}. 
\end{equation}

To complete the proof of the lower bound, we note that for all
2--bridge links except the figure--8 knot and $5_2$ knot, the meridian
$\mu$ is at least $\sqrt[4]{2}$, by work of Adams \cite{adams:waist2}.
Thus the area is at least $8\sqrt{3}\,(t-1)/147$.  Meanwhile, the
figure--8 and $5_2$ knots have twist number $t=2$, hence the estimate
$8\sqrt{3}/147$ is vastly lower than their true cusp area (note a
standard horosphere packing argument implies the area of any cusp is
at least $\sqrt{3}$).

For the upper bound, Futer and Gu\'eritaud \cite[Theorem
B.3]{gf:two-bridge} found that the volume of a hyperbolic 2--bridge
knot satisfies $\vol(S^3\setminus K) < 2v_8(t-1)$.  Again, combining
this with the theorem of B{\"o}r{\"o}czky \cite[Theorem 4]{boroczky},
that a maximal cusp contains at most $\sqrt{3}/(2v_3)$ of the volume
of $M$, we find
$$\vol(C) \: < \: \sqrt{3}\,\, \frac{v_8}{v_3}\,(t-1), \quad
\mbox{hence} \quad \area(\bdy C) \:<\: 2 \sqrt{3}\,\,
\frac{v_8}{v_3}\,(t-1).$$

Finally, for the result on arc length, note that the shortest non-trivial arc
$\gamma$ that starts and ends on the same meridian has length equal to
$\area(\bdy C)/\mu$.  Again using the estimate of Adams
\cite{adams:waist2, adams:2-gen}, the length of a meridian of $K$
satisfies $\sqrt[4]{2} \leq \mu \leq 2$, except if $K$ is the
figure--8 or $5_2$ knot.  Combining Adams's estimates with inequality
(\ref{eqn:2bridge-area}) and arguing as in Lemma
\ref{lemma:ptorus-slope-length} gives the desired lower bound on
$\ell(\gamma)$.  (As above, the figure--8 and $5_2$ knots need to be
checked separately.)
\end{proof}

\section{Volume estimates for closed 3--braids}\label{sec:volume3braids}

In this section, we give a complete characterization of the closed
$3$--braids whose complements are hyperbolic.  Then, we apply Theorem
\ref{thm:4ps-area} from Section \ref{sec:cusp-area} to obtain volume
estimates for closed $3$--braids.

\subsection{A normal form for $3$--braids}
We begin with some notation.  Let $\sigma_1$ and $\sigma_2$ be
generators for $B_3$, the braid group on three strands, as in Figure
\ref{fig:3braid-generators}.  Thus a positive word in $\sigma_1$ and
$\sigma_2^{-1}$ represents an alternating braid diagram.  Let $C =
(\sigma_1 \sigma_2)^3$ denote a full twist of all three strands; $C$
generates the center of $B_3$.  For a braid $w \in B_3$, let $\hat{w}$
denote the link obtained as the closure of $w$.  Note that $\hat{w}$
only depends on the conjugacy class of $w$.  We denote the conjugacy
relation by $\sim$.

In the 1920s, Schreier developed a normal form for this braid group
\cite{schreier}. In particular, he showed the following.

\begin{theorem}[Schreier]\label{thm:schreier-normal}
Let $w \in B_3$ be a braid on 3 strands.  Then $w$ is conjugate to a
braid in exactly one of the following forms:
\begin{enumerate}
\item $C^k \sigma_1^{p_1} \sigma_2^{-q_1} \cdots \sigma_1^{p_s}
\sigma_2^{-q_s}$, where $k \in \ZZ$ and $p_i$, $q_i$, and $s$ are all
positive integers,
\item $C^k \sigma_1^p$, \quad for $k, p \in \ZZ$
\item $C^k \sigma_1 \sigma_2$, \quad for $k \in \ZZ$
\item $C^k \sigma_1 \sigma_2 \sigma_1$, \quad for $k \in \ZZ$, or
\item $C^k \sigma_1 \sigma_2 \sigma_1 \sigma_2$, \quad for $k \in \ZZ$.
\end{enumerate}
This form is unique up to cyclic permutation of the word following
$C^k$.  Braids in form $(1)$ above are called \emph{generic}.
\end{theorem}

Birman and Menasco have shown that nearly every oriented link obtained
as the closure of a 3--braid can be represented by a unique conjugacy
class in $B_3$, with an explicit list of exceptions
\cite{birman-menasco:3braids}.  Thus their theorem, combined with
Schreier's normal form, gives a classification of closed oriented
$3$--braids.  Their paper also contains a modern exposition of
Schreier's algorithm for placing braids in normal form.

Let $K = \hat{w}$ be a closed 3--braid defined by the word $w$, and
let $A$ be the braid axis of $K$. That is, $A$ is an unknot with the
property that the solid torus $S^3 \setminus A$ is swept out by
meridian disks, with each disk intersecting $K$ in 3 points. Then $M_w
:= S^3 \setminus (K \cup A)$ is a 4--punctured sphere bundle over the
circle.  It is well--known, essentially due to work of Thurston
\cite{thurston:hypII}, that the Schreier normal form of $w$ predicts
the geometry of $M_w$.  We include a proof for completeness.

\begin{theorem}\label{thm:hyperbolic-bundle}
$M_w$ is hyperbolic if and only if $w$ is generic.  Moreover, $M_w$
has nonzero Gromov norm if and only if $w$ is generic.
\end{theorem}

\begin{proof}
The braid generators $\sigma_1$ and $\sigma_2^{-1}$ act on the
4--punctured sphere as the standard generators $L$ and $R$ of
$SL_2(\ZZ)$:
$$\sigma_1 \: \mapsto \: L :=
\begin{bmatrix} 1 & 0 \\ 1 & 1 \end{bmatrix}, \qquad
\sigma_2^{-1} \: \mapsto \: R :=
\begin{bmatrix} 1 & 1 \\ 0 & 1 \end{bmatrix}, \qquad
C \: \mapsto \: I =
\begin{bmatrix} 1 & 0 \\ 0 & 1 \end{bmatrix}.$$
Thus generic 3--braids with normal form (1) correspond to positive
words employing both letters $L$ and $R$, hence to pseudo-Anosov
monodromies. Thurston showed that a bundle over $S^1$ with
pseudo-Anosov monodromy is hyperbolic \cite{thurston:hypII}.  More
concretely, Gu\'eritaud showed how to construct the hyperbolic metric
from a positive word in $L$ and $R$ \cite{gf:two-bridge}.  (See
Sections \ref{sec:ptorus} and \ref{sec:4ps} above for a review of the
connection between the monodromy word and the canonical ideal
triangulation of $M_w$.)

The braids with normal form (2) correspond to reducible monodromies of
the form $L^p$.  In this case, $M_w$ is a graph manifold obtained by
gluing two 3--punctured sphere bundles along a torus.  The braids with
normal forms (3--5) correspond to periodic monodromies, hence $M_w$ is
Seifert fibered.  Thus all non-generic normal forms yield
non-hyperbolic manifolds with Gromov norm 0.
\end{proof}
 
\subsection{Hyperbolic $3$--braids}
Our goal in this subsection is to show that the $3$--braids whose closure is a hyperbolic link can be easily identified from their Schreier normal form.  
%	Theorem \ref{thm:hyperbolic-3braid} below characterizes exactly those  3--braids whose closures are hyperbolic. 
%	 The main step of the proof, Lemma \ref{lemma:seifert-braids},
%	describes the Schreier normal forms of 3--braids whose closures are
%	Seifert fibered. 
Because a closed braid presentation of a link $K$ comes with a natural orientation, we need to consider all possible orientations on components of $K$ that are consistent with $K$ being a 3--braid. 

In the lemmas that lead up to Theorem \ref{thm:hyperbolic-3braid}, we rely on two classical invariants that are insensitive to orientation changes on a component of $K$: the (absolute value of) the linking number between components of $K$, and the determinant $\det(K)$. Recall that that $\det(K)$ is the absolute value of the Alexander polynomial of $K$, evaluated at $t=-1$, or  equivalently the absolute value  of the Jones polynomial of $K$, also evaluated at $t=-1$.
It is well--known that reversing the orientation on a component of $K$ leaves the determinant unchanged: from the point of view of the Jones polynomial, this follows because changing the orientation
of some component of $K$ changes the Jones polynomial $J_K(t)$ by a power of $t$.
% This orientation insensitivity of $\det(K)$ is used repeatedly in the proof of Lemma \ref{lemma:seifert-braids}.

\begin{lemma}\label{lemma:determinant}
Let $K = \hat{w}$ be the closure of a generic $3$--braid 
$$w =  C^k
\sigma_1^{p_1} \sigma_2^{-q_1} \cdots \sigma_1^{p_s} \sigma_2^{-q_s},$$
where $p_i$, $q_i$ are all positive. Suppose that $K$ has two or three components and $\det(K) \leq 4$. Then one of two possibilities holds:
\begin{itemize}
\item[(a)] $\det(K) =2$, and $w = C^k a$, where $k$ is even and $a \in \left\{  \sigma_1^2 \sigma_2^{-1} , \,  \sigma_1
\sigma_2^{-2} \right\}$, or
\smallskip
\item[(b)] $\det(K) =4$, and $w = C^k a$, where $k$ is even and $a \in \left\{ \sigma_1^2 \sigma_2^{-2}, \,  \sigma_1^4
\sigma_2^{-1}, \,  \sigma_1 \sigma_2^{-4} \right\}$.
\end{itemize}
\end{lemma}
 
\begin{proof}
The proof uses a result of Murasugi \cite[Proposition 5.1]{murasugi:3braids}. There, Murasugi shows that a generic $3$--braid must have strictly positive determinant. (In his notation, the class of generic $3$--braids is denoted $\Omega_6$.) The same proposition states that, if $\hat{a}$ denotes the closure of $a=\sigma_1^{p_1} \sigma_2^{-q_1} \cdots \sigma_1^{p_s} \sigma_2^{-q_s}$, then
$$\det(K) = \left\{ 
\begin{array}{r l}
\det(\hat{a}), & \mbox{if $k$ is even,} \\
\det(\hat{a})+4, &  \mbox{if $k$ is odd.}
\end{array} \right. $$
% (The above equation can also be derived by substituting $t=-1$ in Lemma \ref{lemma:generic} below.)

If $w$ is generic, then $a$ is also generic, hence $\det(\hat{a})>0$. Thus, if $\det(K) \leq 4$, we must have $k$ even and $\det(K) = \det{\hat{a}}$.
Since $\hat{a}$ is an
alternating link, the minimum crossing number is bounded above by the
determinant (see, for example, Burde and Zieschang
\cite{burde-zieschang:knots}).  Thus the crossing number of $\hat{a}$ is at 
most 4, and we may list the
possibilities for $a$.

Recall that the crossing number of an alternating link is
realized by any alternating diagram without nugatory crossings, and
the only way an alternating 3--braid can have nugatory crossings is if
the braid word is $\sigma_1^r \sigma_2^{-1}$ or $\sigma_1
\sigma_2^{-r}$.  Thus alternating closed 3--braids with crossing
number at most $4$ consist of words of the form $\sigma_1^p$ for
appropriate $p$, $\sigma_1^{p} \sigma_2^{-q}$, for appropriate $p$,
$q$, and $\sigma_1\sigma_2^{-1}\sigma_1 \sigma_2^{-1}$.  All others
will have higher crossing numbers.

Since $\sigma_1^p$ is not generic, we need not consider these.
Since the closed braid corresponding to $C^k \sigma_1 \sigma_2^{-1}
\sigma_1 \sigma_2^{-1}$ has just one component, and we are assuming
$K$ has at least two components, we need not consider these words
either.  Finally, the braids $\sigma_1^p \sigma_2^{-q}$ have
 the appropriate number of crossings for $(p,q) =
(1,1)$, $(1,2)$, $(1,3)$, $(1, 4)$, $(2, 1)$, $(2,2)$, $(3, 1)$, and
$(4, 1)$.  Of these, $(1,1)$, $(1,3)$, and $(3,1)$ have only one link
component. The remaining possibilities are
$$a \,  \in \,  \{  \sigma_1^2 \sigma_2^{-1} , \,  \sigma_1
\sigma_2^{-2},  \, \sigma_1^2 \sigma_2^{-2}, \,  \sigma_1^4
\sigma_2^{-1}, \,  \sigma_1 \sigma_2^{-4} \}.$$
If $a = \sigma_1^2 \sigma_2^{-1}$ or $a=  \sigma_1
\sigma_2^{-2}$, one easily computes that $\det(K) = \det(\hat{a})=2$, and conclusion (a) holds. If $a$ is one of $\sigma_1^2 \sigma_2^{-2}$, $ \sigma_1^4
\sigma_2^{-1}$,  or $\sigma_1 \sigma_2^{-4}$, then $\det(K) = \det(\hat{a})=4$, and conclusion (b) holds.
\end{proof}

We can now restrict the $3$--braids that correspond to Seifert fibered links.

\begin{lemma}\label{lemma:seifert-braids}
Let $K = \hat{w}$ be the closure of a 3--braid $w$, and suppose that
$S^3 \setminus K$ is Seifert fibered. Then $w$ is either non-generic,
or else conjugate to $\sigma_1^p \sigma_2^{\pm 1}$, $\sigma_1^{\pm 1}
\sigma_2^{q}$, or $\sigma_1^{2} \sigma_2^{-2}$.
\end{lemma}

\begin{proof}
A theorem of Burde and Murasugi \cite{burde-murasugi} states that if
$S^3 \setminus K$ is Seifert fibered, then $K$ consists of finitely
many fibers in a (possibly singular) Seifert fibration of $S^3$.  In
case the Seifert fibration of $S^3$ is not singular, the fibration
extends to $S^3$.  The Seifert fibrations of $S^3$ were classified by
Seifert \cite{seifert:fibering} (see also Orlik
\cite{orlik:seifert-manifolds}).  As a consequence, $K$ must be a
\emph{generalized torus link}: a link that is embedded on an unknotted
torus $T$ (this is called an $(m,n)$ torus link) plus possibly one or
both cores of the solid tori in the complement of $T$.

The singular fibration does not extend to $S^3$: it is the product
fibration on a solid torus, in which each fiber is a meridian of the
complementary unknot.  However, note that in this case the result is
again a generalized torus link, with the unknot making up the core of
the $(m, 0)$ torus link.
The question of which closed 3-braids represent generalized torus
links has been studied by Murasugi in \cite{murasugi:3braids}.

In an $(m,n)$ torus link, we may assume without loss of generality
that $m>0$, and that either $n=0$ or $\abs{n} \geq m$.  With this
normalization, a theorem of Schubert \cite{schubert:braids} implies
that the bridge number of the $(m,n)$ torus link is $m$. Since the
braid index is greater than or equal to the bridge number, any choice
of orientation on the components of an $(m,n)$ torus link must yield a
braid index of at least $m$.
Thus, if we add $c$ cores of solid tori and obtain a $3$--braid, $1
\leq m \leq 3-c$.  There are three cases, conditioned on the value of
$c$.

\smallskip
\underline{Case 0: $c=0$.}  Then $K$ is an $(m, n)$ torus link, where
$1 \leq m \leq 3$.
%% Much of the argument in this case is handled by
%% Murasugi in \cite[Section 12]{murasugi:3braids}.
Murasugi classifies the torus links that can be written as closed
3--braids in \cite[Section 12]{murasugi:3braids}.  However, a careful
reading of his proofs indicates that he is assuming multiple
components of torus links are always oriented consistently.
For our purposes, we will also need to consider the torus links in 
which the orientation of some component is reversed.
In fact, the machinery developed by Murasugi in this monograph
is sufficient to handle all choices of orientation.  For ease of reading, we will include the
arguments for all the cases not immediately apparent from
\cite{murasugi:3braids}.

If $m=1$, then $K$ is the unknot, and Theorem 12.1 in Murasugi's
monograph shows $w \sim \sigma_1^{\pm 1} \sigma_2^{\pm 1}$.  If $m=2$,
then $K$ is a $(2,n)$ torus link.  If the (one or two) components of $K$ are 
oriented consistently, then it is an \emph{elementary}
torus link in Murasugi's terminology.  Then, Theorem 12.3 of his
monograph implies $w \sim \sigma_1^p \sigma_2^{\pm 1}$ or $w\sim
\sigma_1^{\pm 1} \sigma_2^{q}$.

Now, suppose that $n=2\ell$, and $K$ is a $(2, 2\ell)$ torus link whose
components have opposite orientations.  Then $K$ is the oriented
boundary of a Seifert surface that is an annulus.  Using this annulus,
we calculate that the Alexander polynomial of $K$ is $\Delta_K(t)=
\ell(1-t)$; thus $\deg{\Delta_K(t)}=1$.  By Proposition 8.1 of
\cite{murasugi:3braids}, $K$ can be the closure of a 3--braid only if
$\abs{\ell} \leq \frac{1}{2}{\deg \Delta_K(t)}+2=\frac{5}{2}$.  Since the
determinant of the $(2, 2\ell)$ torus link is $\abs {2\ell}$, and this fact
remains true under orientation reversal of components, $\det(K) =
\abs{2\ell} \leq 4$.  Thus, if $w$ is a generic $3$--braid that
represents $K$ with this orientation, Lemma \ref{lemma:determinant}
applies.

If $\det(K)= 2$, Lemma \ref{lemma:determinant} says the generic braid
$w$ representing $K$ must be conjugate to $C^k \sigma_1^2
\sigma_2^{-1}$ or $w \sim C^k \sigma_1 \sigma_2^{-2}$, where $k$ is even.  
Without loss
of generality, say $w \sim C^k \sigma_1^2 \sigma_2^{-1}$.  Then the
two components of $\hat{w}$ have linking number equal to $\abs{2k+1}$.
But since $K$ is the $(2, \pm2)$ torus link, the linking number of the
two components is $1$.  Thus $k=0$, and $w \sim \sigma_1^2
\sigma_2^{-1}$, as desired.

If $\det(K)= 4$, Lemma \ref{lemma:determinant} says the only generic
braids representing $K$ must be conjugate to $C^k \sigma_1^2
\sigma_2^{-2}$, $C^k \sigma_1^4 \sigma_2^{-1}$, or $C^k \sigma_1
\sigma_2^{-4}$.  The closure of $C^k \sigma_1^2 \sigma_2^{-2}$ has
three components, so we don't need to consider this case.  If $w \sim
C^k \sigma_1^{4} \sigma_2^{-1}$, then the two components of $\hat{w}$
have linking number $\abs{2+2k}$, which is equal to $2$ by hypothesis.
Thus $k$ is $-2$ or $0$.  But if $k = -2$, the closure of $C^{-2}
\sigma_1^4 \sigma_2^{-1}$ is a hyperbolic link, a contradiction.  Thus
$k=0$, and $w \sim \sigma_1^4 \sigma_2^{-1}$, as desired.  The
argument when $w \sim C^k \sigma_1 \sigma_2^{-4}$ is identical.

Finally, suppose that $m=3$, and $K$ is a $(3,n)$ torus link. Then $K$ 
has either one or three components.
If $K$ is a knot, Proposition 12.3 of \cite{murasugi:3braids} shows that $K$ cannot
be represented by a generic $3$--braid. In fact, Murasugi's
argument also works for links, with all orientations, but for completenss
we include the argument here.

If $K$ is a 3--component torus link $(3, 3\ell)$,  the linking
number of any two components of $K$ has absolute value $\abs{\ell}$.
Proposition 3.3 of \cite{murasugi:3braids} implies
that $\det(K)$ is either $0$ or $4$. By Lemma \ref{lemma:determinant}, 
any generic braid
$w\sim C^k a$ representing $K$, must be conjugate to $C^k
\sigma_1^2 \sigma_2^{-2}$, as all the other possibilities in the lemma 
have two components. Now, consider the pairwise linking numbers between
components in the closure of $w\sim C^k \sigma_1^2 \sigma_2^{-2}$.
These pairwise linking numbers are equal in absolute value to
$\abs{k}$, $\abs{k+1}$, and $\abs{k-1}$.  But all these numbers must be
equal to $\abs{\ell}$, a contradiction.

\smallskip
\underline{Case 1: $c=1$.}  Then $K = L_t \cup L_a$, where $L_t$ is a
$(m,n)$ torus link with $m=1$ or $m=2$, and $L_a$ is the core of one
of the two solid tori.  If $m=1$, then (depending on the choice of
solid torus) $K$ is either a $(2, 2n)$ torus link and we reduce to
case 0, or $K$ is the Hopf link, which is the $(2,2)$ torus link, and
we again reduce to case 0.  Thus we may suppose that $L_t$ is a
$(2,n)$ torus link.  Then, as an unoriented link, $K$ admits a diagram
in one of two possible forms, shown in Figure
\ref{fig:augmented-torus}, depending on which solid torus $L_a$ came
from.

\begin{figure}%[h]
\psfrag{n}{$n$}
\begin{center}
\includegraphics{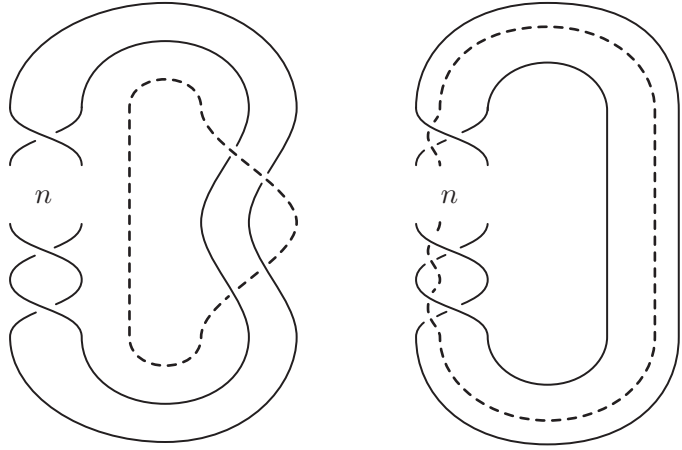}
\end{center}
\caption{Case 1 of Lemma \ref{lemma:seifert-braids}: two ways to add
the core of a solid torus to a $(2,n)$ torus link.}
\label{fig:augmented-torus}
\end{figure}

\smallskip
\underline{Subcase 1a:}  Suppose that the link $K$ is depicted in the left panel of Figure 10.
Then we can characterize the linking number $lk(L_t, L_a) $ as follows.  If $n$ is odd and $L_t$ is a knot, then its
linking number with $L_a$ is $\pm 2$; if $n$ is even and $L_t$ has two
components, then each component of $L_t$ has linking number $\pm 1$
with $L_a$. Note that the absolute value of the linking number is insensitive to changes of orientation.

If all the strands in Figure \ref{fig:augmented-torus} are oriented counterclockwise, the link $K$ in the 
left panel of the figure can
be represented by the braid word $v = \sigma_1^n \sigma_2 \sigma_1^2
\sigma_2$.  Of course, \emph{a priori}  there may be other braid
representatives, possibly corresponding to other choices of orientation on one or more components of $K$.  
Nonetheless, knowing that $K$ can be obtained as  the closure
$\hat{v}$ of the braid represented by $v$ allows us to compute link
invariants.  Because the normal form of $v$ is $C \,
\sigma_1^{n-2}$, Murasugi's Proposition 3.6 gives that $\det(K)=4$. 
Suppose that, with some choice of orientation, $K$ is represented by the generic braid $w$. Since $\det(K)=4$, Lemma \ref{lemma:determinant} implies the normal form of $w$ must be one of  $C^k \sigma_1^2 \sigma_2^{-2}$, $C^k \sigma_1^4
\sigma_2^{-1}$, or $C^k \sigma_1 \sigma_2^{-4}$.

If $w \sim C^k \sigma_1^2 \sigma_2^{-2}$, then $K = \hat{w}$ is a
3--component link.  The $2$--component links contained in $K$ have
pairwise linking numbers equal in absolute value to $\abs{k}$,
$\abs{k+1}$, and $\abs{k-1}$.   By hypothesis,  two of these linking numbers must be equal to 1.
It follows that the only possibility is $k = 0$, hence
$w \sim \sigma_1^2 \sigma_2^{-2}$, as desired.

If $w \sim C^k \sigma_1^{4} \sigma_2^{-1}$, then $K = \hat{w} = L_t
\cup L_a$ is a 2--component link.  In this case, we compute that
$\abs{lk(L_t, L_a)} = \abs{2+2k}$, which is equal to $2$ by
hypothesis.  Thus $k$ is $-2$ or $0$.  But if $k = -2$, the closure of
$C^{-2} \sigma_1^4 \sigma_2^{-1}$ is a hyperbolic link, a
contradiction.
Thus $k=0$, and $w \sim \sigma_1^4 \sigma_2^{-1}$, as desired. The
case when $w \sim C^k \sigma_1 \sigma_2^{-4}$ is identical.

\smallskip
\underline{Subcase 1b:} Suppose that the unoriented link $K$ is depicted in the right panel of Figure 10.
It follows that  $\, lk(L_t, L_a) = \pm n,$ depending on choices
of orientations of the components.
If all the strands of $K$ are oriented counterclockwise, $K$ can be represented by the braid word 
$v = (\sigma_1 \sigma_2 \sigma_1)^n$. (Just as in Subcase 1a, there may be other braid representatives, 
but knowing one braid representative $v$ allows us to compute invariants that are insensitive to orientation.) 
Using the
braid relation $\sigma_1 \sigma_2 \sigma_1 = \sigma_2 \sigma_1
\sigma_2$, we can rewrite $v$ as
$$v = \left\{ \begin{array}{l l} C^{n/2}, & \mbox{if $n$ is even,} \\
C^{(n-1)/2} \sigma_1 \sigma_2 \sigma_1, & \mbox{if $n$ is odd.}
\end{array} 
\right. $$
If $n$ is even, then $\hat v$ is a torus link. Thus the unoriented link $K$ is a torus link, and we reduce to Case 0. 
If $n$ is odd, then Murasugi's Proposition 3.5 gives $\det(K)=2$. 
Thus, by Lemma \ref{lemma:determinant}, any generic braid $w$ that also represents $K$ must have normal form
$C^k \sigma_1^{2} \sigma_2^{-1}$ or $C^k \sigma_1 \sigma_2^{-2}$.

If $w \sim C^k \sigma_1^{2} \sigma_2^{-1}$, then $K = \hat{w} = L_t
\cup L_a$ is a 2--component link.  We may immediately compute that one
component $L_t$ is the $(2, 2k-1)$ torus knot, the other component
$L_a$ is the unknot, and $\abs{lk(L_t, L_a)} = \abs{2k+1}$.  Because
by assumption, $L_t$ is the $(2, n)$ torus link, this implies
$\abs{n}=\abs{2k-1}$.  Additionally, since by assumption $lk(L_t,
L_a)=\pm n$, we may conclude that $\abs{n}=\abs{2k+1}$.  This is possible
only if $k=0$.  So $w \sim \sigma_1^2 \sigma_2^{-1}$, as desired.  The
case when $w \sim C^k \sigma_1 \sigma_2^{-2}$ is identical. 

\smallskip
\underline{Case 2: $c=2$.} Then $K = L_t \cup L_a \cup L_b$, where
$L_t$ is an unknot on $T$, and $L_a$ and $L_b$ are cores of the two
solid tori.  Since $L_t$ is a $(1,n)$ curve on the torus $T$, one of
the cores $L_a$ or $L_b$ (say, $L_b$) can be isotoped to lie on $T$,
disjointly from $L_t$. Thus, as an unoriented link, $L_t \cup L_b$ is a  torus link on
$T$, and this case reduces to Case 1.
\end{proof}

The next theorem characterizes the 3-braids whose closures  represent hyperbolic links.

\begin{theorem}\label{thm:hyperbolic-3braid}
Let $w \in B_3$ be a word in the braid group, and let $K \subset S^3$
be the link obtained as the closure of $w$. Then $S^3 \setminus K$ is
hyperbolic if and only if $w$ is generic and not conjugate to
$\sigma_1^p \sigma_2^q$ for arbitrary integers $p, q$.
\end{theorem}

\begin{proof}
First, we check the ``only if'' direction.  If $w$ is non-generic,
then by Theorem \ref{thm:hyperbolic-bundle}, $M_w = S^3 \setminus (K
\cup A)$ is a graph manifold with Gromov norm 0.  Since the Gromov
norm of a manifold cannot increase under Dehn filling
\cite[Proposition 6.5.2]{thurston:notes}, $S^3 \setminus K$ also has
Gromov norm 0, and is not hyperbolic.  If $w$ is generic and conjugate
to $\sigma_1^p \sigma_2^{\pm 1}$, then $K$ is a $(2,p)$ torus link
(similarly for $\sigma_1^{\pm 1} \sigma_2^{q}$).  Finally, if $w$ is
conjugate to $\sigma_1^p \sigma_2^q$, where $|p|, |q| \geq 2$, then
$K$ is the connected sum of $(2,p)$ and $(2,q)$ torus links, hence
cannot be hyperbolic.

For the ``if'' direction, suppose that $S^3 \setminus K$ is not
hyperbolic.  Then, by Thurston's hyperbolization theorem
\cite{thurston:survey}, it is reducible, toroidal, or Seifert fibered.
If $S^3\setminus K$ is Seifert fibered, then Lemma
\ref{lemma:seifert-braids} implies $w$ is non-generic or conjugate to
$\sigma_1^p \sigma_2^{q}$.  Meanwhile, if $S^3 \setminus K$ is reducible,
then $K$ is a split link.  By a theorem of Murasugi \cite[Theorem
5.1]{murasugi:3braids}, this can only happen if $w \sim \sigma_1^p$:
hence, $w$ is not generic.

Finally, suppose that $S^3 \setminus K$ contains an essential torus
$T$.  If $w$ is not generic, then we are done.  If $w$ is generic, a
theorem of Lozano and Przytycki \cite[Corollary 3.3]{lozano-przytycki}
says that $T$ always has meridional compression disks, i.e. there is
some disk $D \subset S^3$ such that $D\cap T = \partial D$ and $D\cap
K$ is a point.\protect\footnote{Lozano and Przytycki's result is
stated for ``hyperbolic'' 3--braids. However, their definition of
\emph{hyperbolic} is the same as our definition of \emph{generic}.}
After meridionally compressing $T$, i.e. after replacing a
neighborhood of $\partial D$ on $T$ with two parallel copies of the
annulus $D\setminus K$, we obtain an essential, meridional annulus
that splits $K$ into connected summands.  But by a theorem of Morton
\cite{morton:composite-braids}, a braid $w \in B_3$ represents a
composite link if and only if $w \sim \sigma_1^p \sigma_2^q$, where
$|p|, |q| \geq 2$. See Birman and Menasco \cite[Corollary 1]{birman-menasco:toroidal}
for another way to identify the toroidal $3$--braids.
\end{proof}

\subsection{Volume estimates}
For sufficiently long generic $3$--braids, the methods of the previous
sections estimate hyperbolic volume.

\begin{theorem}\label{thm:3braid-volume}
Let $K= \hat{w}$ be the closure of a generic 3--braid $w \sim C^k
\sigma_1^{p_1} \sigma_2^{-q_1} \cdots \sigma_1^{p_s} \sigma_2^{-q_s}$,
where $p_i$, $q_i$ are all positive and $w \nsim
\sigma_1^p\sigma_2^{-q}$.  Then $K$ is hyperbolic, and
\begin{equation}\label{eq:3braid-volume}
4v_3 \, s - 276.6 \: < \: \vol(S^3 \setminus K) \: < \: 4v_8 \, s.
\end{equation}
Furthermore, the multiplicative constants in both the upper and lower
bounds are sharp.
\end{theorem}

\begin{proof}
Let $A$ be the braid axis of $K$. Then $M_w = S^3 \setminus (K \cup
A)$ is a 4--punctured sphere bundle with monodromy
$$\Omega = L^{p_1} R^{q_1} \cdots L^{p_s} R^{q_s}.$$
Futer and Gu\'eritaud showed \cite[Corollary B.2]{gf:two-bridge} that
the length of $\Omega$ coarsely determines the volume of $M$:
\begin{equation}\label{eq:bundle-volume}
4 v_3 \, s \: \leq \: \vol(M_w) \: < \: 4v_8 \, s,
\end{equation}
where both the upper and lower bounds are sharp.  That is: there exist
$4$--punctured sphere bundles that realize the lower bound, and other
bundles that are $\varepsilon$--close to the upper bound.  Since $S^3
\setminus K$ is obtained by Dehn filling on $M_w$, the same upper
bound applies to the volume of $S^3 \setminus K$.  Furthermore, by
choosing an extremely long filling slope (which will happen when
$\abs{k} \to \infty$), one can arrange for $\vol(S^3 \setminus K)$ to
be arbitrarily close to $4 v_8 \, s$.

For the lower bound on volume, we rely on Theorem \ref{thm:4ps-area}.
That theorem states that the meridian of $A$ (which will be transverse
to the fibers) has length at least $8\sqrt{3} \, s/147$.  In
particular, when $s \geq 67$, the meridian will be longer than $2\pi$.
Thus we may apply Theorem 1.1 of \cite{fkp-07}, which estimates the
change in volume under Dehn filling along slopes longer than $2\pi$.
For all $s \geq 67$, we obtain
$$
\begin{array}{l l l l l}
\vol(S^3 \setminus K) 
&\geq&  \left(1-\left(\dfrac{2\pi}{8\sqrt{3} \, s/147}\right)^2\right)^{3/2} \vol(M_w), 
& \mbox{by \cite[Thm 1.1]{fkp-07} and Thm \ref{thm:4ps-area} }\\
&\geq&  \left(1-\dfrac{7203  \, \pi^2 }{16 \, s^2} \right)^{3/2} 4 v_3 \, s,
& \mbox{by inequality (\ref{eq:bundle-volume})}
\end{array}
$$
Note that by calculus,
$$\left(1-\dfrac{7203 \, \pi^2}{16\,s^2}\right)^{3/2} 4 v_3 \, s -
4v_3\,s$$
has a minimum of $-276.52\cdots$ for $s\geq 67$.  Thus $\vol(S^3
\setminus K) > 4v_3\,s - 276.6$.

On the other hand, if $s \leq 67$, then $4v_3 \, s - 276.6 < 0$, hence
the volume estimate is trivially true.  Thus the lower bound on volume
holds for all hyperbolic $3$--braids.

Finally, to show sharpness of the multiplicative constant in the lower
bound, consider $3$--braids of the form $w=(\sigma_1
\sigma_2^{-1})^s$.  In the proof of Theorem B.1 of
\cite{gf:two-bridge}, it was shown that for the closures of these
braids, $\vol(M_w) = 4v_3s$.  Since $\vol(S^3\setminus K) < \vol(M_w)
= 4v_3s$ for this sequence of knots, the multiplicative constant
$4v_3$ must be sharp.
\end{proof}

We close this section with an interesting side comment.  Theorem
\ref{thm:3braid-volume} compares in an intriguing way to prior
results that estimate the volume of a link complement in terms of the
twist number of a diagram. (See Section \ref{sec:diagram-intro} and
the introduction of \cite{fkp:coils} for definitions and background.)
In the braid word $w$, each term $\sigma_1^{p_i}$ or $\sigma_2^{-q_i}$
corresponds to a twist region involving a pair of strands of $K$.
Meanwhile, when $k \neq 0$, the term $C^k$ defines a single \emph{generalized twist
region}, in which we perform $k$ full twists on all three strands of
the braid.
Altogether, the braid word $w$ defines a diagram with either $2s$ or $2s+1$
generalized twist regions -- including $2s$ ordinary twist regions
twisting on two strands of $K$. As a result, Theorem
\ref{thm:3braid-volume} can be reformulated in the following way.

\begin{corollary}\label{cor:3braid-twists}
Let $K$ be a hyperbolic closed 3--braid, and let $D(K)$ be the braid 
diagram corresponding to the Schreier normal form for $K$. If
$\tgen$ denotes the number of generalized twist regions in the diagram
$D$, then
$$2v_3 \, \tgen - 279 \: < \: \vol(S^3 \setminus K) \: < \: 2v_8 \,
\tgen.$$
\end{corollary}

\begin{proof}
By Theorem \ref{thm:hyperbolic-3braid}, $K$ must be represented by a
generic word $w = C^k \sigma_1^{p_1} \sigma_2^{-q_1} \cdots
\sigma_1^{p_s} \sigma_2^{-q_s}$. Substituting $2s \leq \tgen \leq 2s+1$ into
Theorem \ref{thm:3braid-volume} gives the desired volume estimate.
\end{proof}

\begin{remark}
In \cite[Corollary 3.2]{fkp:coils}, we show that the twist number alone, as opposed
to the generalized twist number, is not a good measure of the volume
of 3--braids.  Thus the single generalized twist region from the term
$C^k$ is important in the corollary above.
\end{remark}

\section{The Jones polynomial  and volume of closed
	3--braids}\label{sec:jones3braids}

% \subsection{Volume and Jones Polynomial relations}

In this section, we will apply the previous results to the Jones
polynomial of a closed 3--braid.  We begin by relating the Jones
polynomial of a closed $3$--braid to the Schreier normal form of the
braid.  By applying Theorem \ref{thm:3braid-volume}, we will show in
Theorem \ref{thm:relations} that certain coefficients of the Jones
polynomial are bounded in terms of the volume.  At the end of the
section, we will prove Theorem \ref{thm:counterexample}, which shows
that no function of $\beta_K$ and $\beta'_K$ can coarsely
predict the volume of all hyperbolic knots.

\subsection{The Jones polynomials of generic 3--braids}\label{subsec:generic-jones}

In the case that $K$ is the closure of a 3--braid, we need to relate
the Jones polynomial to the Schreier normal form of the braid. (See
Theorem \ref{thm:schreier-normal}.)  Here, we will concern ourselves
with 3--braids whose Schreier normal forms are \emph{generic}.  That
is, we will consider braids $b\in B_3$ written in the form
$$b= C^k \sigma_1^{p_1} \sigma_2^{-q_1} \cdots \sigma_1^{p_s}
\sigma_2^{-q_s},$$
where $p_i, q_i, k \in \ZZ$, with $p_i, q_i>0$, and $C:=(\sigma_1
\sigma_2 \sigma_1 )^2$.  We set
$${\bf p}:=\sum_{i=1}^s p_i, \ \ {\rm and} \ \ {\bf q}:=\sum_{i=1}^s
q_i.$$

The \emph{exponent} $e_b$ of a braid $b$ is the signed sum of its
powers.  Thus for an alternating braid $a$, $e_a = {\bf p} - {\bf q}$,
and if $b=C^k a$, then $e_b = 6k + {\bf p} - {\bf q} = 6k + e_a$.  The
exponent $e_b$ is closely related to the \emph{writhe} of a diagram,
namely the algebraic sum of oriented crossings. Because both of the
generators $\sigma_1$ and $\sigma_2$ depicted in Figure
\ref{fig:3braid-generators} are negative crossings, the writhe of the
standard diagram of a closed 3--braid is $w(D_b) = - e_b$.

For a braid $b\in B_3$, let $\hat {b}$ denote the closure of $b$. Let
$K$ denote the link type represented by ${\hat b}$ and let $J_{K}(t)$
denote the Jones polynomial of $K$. We write
\begin{equation}
J_{K}(t)= \alpha_K t^{M(K)}+ \beta_K t^{M(K)-1}+ \ldots + \beta'_K
t^{m(K)+1}+ \alpha'_K t^{m(K)},
\label{eqn:jp}
\end{equation}
so that $M(K)$ is the highest power of $t$ in $J_K(t)$ and $m(K)$ is
the lowest power of $t$ in $J_K(t)$.  Now the second and next-to-last
coefficients of $J_{K}(t)$ are $\beta_K$ and $\beta'_K$, respectively.

\begin{figure}[ht] 
\centerline{\input{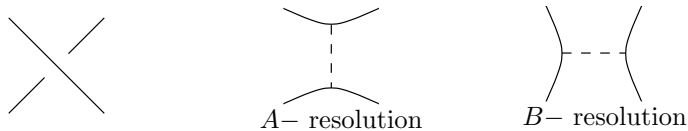}}
\caption{Resolutions of a crossing}
\label{fig:splicing}
\end{figure}

We will also need the following definitions.  Associated to a link
diagram $D$ and a crossing of $D$ are two link diagrams, each with one
fewer crossing than $D$, called the \emph{$A$--resolution} and
\emph{$B$--resolution} of the crossing.  See Figure
\ref{fig:splicing}.
	
Starting with any $D$, let $s_A(D)$ (resp. $s_B(D)$) denote the
crossing--free diagram obtained by applying the $A$--resolution (resp.
$B$--resolution) to all the crossings of $D$.  We obtain graphs
$\GA(D)$, $\GB(D)$ as follows: The vertices of $\GA (D)$ are in
one-to-one correspondence with the components of $s_A(D)$.  For every
crossing of $D$, we add an edge between the vertices of $\GA(D)$ that
correspond to the components of $s_A(D)$ at that crossing.  In a
similar manner, construct the $B$--graph $\GB$ by considering
components of $s_B(D)$.  A link diagram $D$ is called \emph{adequate}
if the graphs $\GA(D)$, $\GB(D)$ contain no 1--edge loops, i.e. if
there are no edges with both ends at the same vertex.

Let $v_A(D)$, $e_A(D)$ (resp. $v_B(D)$, $e_B(D)$) denote the number of
vertices and edges of $\GA(D)$ (resp. $\GB(D)$).  The reduced graph
$\GRA(D)$ is obtained from $\GA(D)$ by removing multiple edges
connected to the same pair of vertices; similarly one has the reduced
graph $\GRB(D)$.  Let $e'_A(D)$ (resp. $e'_B(D)$) denote the number
edges of $\GRA(D)$ (resp. $\GRB(D)$).

The following results about Jones polynomials of adequate knots are
well known.

\begin{lemma}\label{lemma:adequate-jones}
Let $D$ be an adequate diagram of a link $K$, whose Jones polynomial
is written as in equation (\ref{eqn:jp}). Then the top and bottom
coefficients of $J_K(t)$ satisfy
$$\alpha_K = (-1)^{v_B(D) + w(D) -1}, \qquad
  \beta_{K} = (-1)^{v_B(D) + w(D)} (e'_B(D) - v_B(D) + 1),$$
$$\alpha'_K = (-1)^{v_A(D) + w(D) -1} , \qquad
	\beta'_{K} = (-1)^{v_A(D) + w(D)} (e'_A(D) - v_A(D) + 1).$$
\end{lemma}

\begin{proof}
Let $\ZZ [A,\ A^{-1}]$ denote the ring of Laurent polynomials in a
variable $A$, with integer coefficients.  Recall that the Kauffman
bracket of the diagram $D$, denoted by $\left < D\right >$, is an
element in $\ZZ [A, \ A^{-1}]$ such that
\begin{equation}
  J_{K_a}(t)= (-A)^{-3 w(D)} \left< D\right>
  \left|_{\displaystyle{A=t^{-1/4}}}\right., 
  \label{eqn:Jones-Kauffman}
\end{equation}
where $w(D)$ is the writhe of $D$, or the algebraic sum of crossings.
Now, Kauffman showed that the first and last coefficients of $\left< D
\right>$ are $(-1)^{v_A(D)-1}$ and $(-1)^{v_B(D)-1}$, respectively.
(See \cite[Theorem 6.1]{dasbach-futer2...} for a proof.) Meanwhile,
Stoimenow showed \cite{stoimenow:amphichiral} that the second coefficient is given by
$(-1)^{v_A(D)}(e'_A(D) - v_A(D) + 1)$, and similarly the
next--to--last coefficient is given by $(-1)^{v_B(D)}(e'_B(D) - v_B(D)
+ 1)$.  See \cite[Corollary 6.3]{dasbach-futer2...} for an alternate proof.

Next, we multiply $\left< D \right>$ by $(-A)^{-3w(D)}$. As a result,
all the coefficients are multiplied by $(-1)^{-3w(D)} =
(-1)^{w(D)}$. Finally, to recover the Jones polynomial, we substitute
$A = t^{-1/4}$. As a result, the highest powers of $A$ will correspond
to the lowest powers of $t$, and vice versa. Thus the top and bottom
coefficients of $J_K(t)$ are as claimed.
\end{proof}

Given a condition $R$, let $\delta_R$ be the characteristic function of $R$: its value is $1$ when 
$R$ is true, and $0$ when $R$ is false. The characteristic functions $\ep: = \delta_{\p \leq 2}$ and 
$\eq: = \delta_{\q \leq 2}$ will be particularly useful for expressing the Jones polynomials of $3$--braids.

\begin{lemma}\label{lemma:alternating}
Suppose that a link $K_a$ is the closure of an alternating 3-braid
$$a= \sigma_1^{p_1} \sigma_2^{-q_1} \cdots \sigma_1^{p_s}
\sigma_2^{-q_s},$$
with $p_i, q_i>0$.  Suppose as well that $\p > 1$ and $\q > 1$.
Then the following hold:
\begin{enumerate}
  \item[(a)] The highest and lowest powers of $t$ in $J_{K_a}(t)$ are
$$M(K_a)=\frac{3\q-\p}{2} \quad \mbox{and} \quad
		m(K_a)= \frac{\q-3\p}{2}.$$
  \item[(b)] The first two and last two coefficients in $J_{K_a}(t)$
  are
$$\alpha_{K_a} = (-1)^{\p}, \quad \beta_{K_a} = (-1)^{\p+1} (s - \eq), \quad 
\beta'_{K_a} = (-1)^{\q+1}  (s - \ep), \quad \alpha'_{K_a} =(-1)^{\q},$$
where $\ep$ equals $1$ if $\p = 2$ and $0$ if $\p > 2$, and similarly for $\eq$.
\smallskip

  \item[(c)] The third and third-to-last coefficients in $J_{K_a}(t)$ satisfy
$$1 \leq (-1)^{\p} \gamma_{K_a} \leq \frac{s(s+3)}{2} \quad \mbox{and} \quad
1 \leq (-1)^{\q} \gamma'_{K_a} \leq \frac{s(s+3)}{2}.$$
% except if $K_a$ is the figure--$8$ knot (and $\p + \q =4$) or if $K_a$ is the Whitehead link (and $\p + \q =5$).
These upper bounds are attained when $s>1$ and $p_i, q_i > 1$ for all $i$. See Equation (\ref{eq:gamma-exact})
for a precise formula for $\gamma_{K_a}$.
\end{enumerate}
\end{lemma}

\begin{proof}
Consider the link diagram $D:={\hat a}$ obtained as the closure of $a$.
Note that the diagram $D$ is 
alternating and reduced, i.e.\ it
contains no nugatory crossings (here, we are using the hypothesis that if
$s=1$, then $p_1 = \p >1$ and $q_1 = \q >1$).  This implies that $D$ is an adequate
diagram (compare \cite[Proposition 5.3]{lickbook}).  Thus we may use Kauffman's work 
to find the highest and lowest powers of $A$, and 
Lemma \ref{lemma:adequate-jones} to find the coefficients.

Kauffman showed that the highest and
lowest powers of $A$ in the bracket polynomial $\left< D\right>$ are $c(D)+2v_A(D) -2$ and
$-c(D)-2v_B(D)+2$, respectively.  (See \cite{lickbook} for an
exposition.  A proof from the graph theoretic viewpoint can be found
in \cite[Proposition 7.1]{dasbach-futer...}.) In our setting, the crossing number is $c(D)={\bf p}+{\bf q}$,  $v_A(D)={\bf p}+1$, and $v_B(D)={\bf q}+1$.
Thus the highest and lowest powers of $A$ in $\left< D \right>$ are
$3\p+\q$ and $-3\q-\p$, respectively.

By equation (\ref{eqn:Jones-Kauffman}) we multiply $\left< D \right>$
by $(-A)^{-3w(D)} = (-A)^{-3\q+3\p}$. (Recall that for a closed 3--braid, $w(D)
= -e_a={\bf q}-{\bf p}$.) Thus the highest power of $A$
becomes $6\p-2\q$, and the lowest becomes $2\p-6\q$. Then to obtain
$J_{K_a}(t)$, replace $A$ by $t^{-1/4}$.  Thus the highest power of $t$ in
$J_{K_a}(t)$ corresponds to $-1/4$ times the lowest power of $A$, and vice versa:
$$M(K_a) = \frac{-(2\p-6\q)}{4} = \frac{3\q-\p}{2}
\quad \mbox{and} \quad
m(K_a) = \frac{\q-3\p}{2}.$$

Next, we turn our attention to the top and bottom coefficients of $J_{K_a}$. In the following calculation, we 
focus on the first three coefficients $\alpha_{K_a}$ $\beta_{K_a}$, and $\gamma_{K_a}$. 
By Lemma  \ref{lemma:adequate-jones} and \cite{dasbach-lin:head-tail, 
stoimenow:amphichiral}, these top coefficients only depend on the $B$--resolution of $D$. 
To find the last three coefficients, one merely needs to interchange $\p$ with $\q$ in all the formulas.

From the $B$--resolution of the diagram $D$, one can readily compute that
$v_B(D)={\bf q}+1$. Thus, by Lemma \ref{lemma:adequate-jones},
$\alpha_{K_a} = (-1)^{2\q-\p}=(-1)^{\p}$. 
The total number of edges of $\GB$ is $e_B(D) = \p + \q$, the number of crossings in $D$. Out of this total number, the $\p$ edges corresponding to the powers of $\sigma_1$ will become identified to $s$ classes of edges in $e'_B$ (with one edge for each $\sigma_1$--twist region. The $\q$ edges corresponding to the powers of $\sigma_2$ all survive in $e'_B$, except in the special case when $\q=2$, when a loop of two edges running all the way around the braid counts for only a single reduced edge in $e'_B$.
Putting it together, the number of reduced edges will be $e'_B(D)= \q - \eq+ s$, hence
\begin{equation}\label{eq:beta}
\beta_{K_a} =(-1)^{2\q - \p+1} (e'_B - v_B + 1) = (-1)^{\p+1} (s - \eq).
\end{equation}

To find $\gamma_{K_a}$, we must calculate
the third-to-last
coefficient of the Kauffman bracket $\langle D \rangle$, and then multiply by 
$(-1)^{w(D)} = (-1)^{\q-\p}$. Thus $\gamma_{K_a}$ is the third-to-last coefficient of $\langle D \rangle$, which 
was computed in closed form by Dasbach and Lin \cite[Theorem 4.1]{dasbach-lin:head-tail} and {ozawa:adequate}
 \cite[Proposition 3.3]{stoimenow:amphichiral}.
According to their formula,
\begin{equation}\label{eq:gamma}
\gamma_{K_a}\:=\:(-1)^{\p} \left( \frac{ \abs{ \beta_{K_a}} (\abs{\beta_{K_a}}+1)}{2}-\theta+\mu-\tau \right) \!,
\end{equation}
where $\theta$, $\mu$, and $\tau$ are defined as follows.
The quantity $\theta$ is always zero for reduced alternating
 diagrams; this is because the circles $s_B$ do not nest on
the projection plane.  The quantity $\mu$ is 
the number of edges in the reduced graph $\GRB$ whose multiplicity in
$\GB$ is greater than one. By the argument preceding equation (\ref{eq:beta}),
\begin{equation}\label{eq:mu}
\mu \: = \: \#\{i : p_i > 1 \} + \eq \: = \: s -   \#\{i : p_i = 1 \} + \eq.
\end{equation}

Finally, the quantity $\tau$ is
defined to be number of  loops  in $\GRB$ that consist of exactly $3$ edges. 
In our context, one of these loops (surrounding a region of the diagram $D$) 
typically arises for any $i$ where $q_i= 1$. However, in the special case when 
$s=\q=2$, two of these regions involve the same triple of vertices in $\GB$, and 
account for the same loop in $\GRB$. Furthermore, a loop of length $3$ goes all 
the way around the braid when $\q = 3$. All together,
\begin{equation}\label{eq:tau}
\tau =  \#\{i : q_i = 1 \} - \delta_{s=\q=2} + \delta_{\q = 3}.
\end{equation}
%	Combining the information on $\theta$, $\mu$, and $\tau$ gives
%	$$ -\theta+\mu-\tau  \quad = \quad 
%	s \, - \,   \#\{i : p_i = 1 \}  \,-\,   \#\{i : q_i = 1 \} \, + \,  \eq  \, + \, \delta_{s=\q=2} \, - \, \delta_{\q = 3}. $$

Thus, plugging equations (\ref{eq:beta}), (\ref{eq:mu}), (\ref{eq:tau}) into equation 
(\ref{eq:gamma}) gives
\begin{eqnarray}
(-1)^{\p} \gamma_{K_a} &=&
\frac{ (s-\eq) (s-\eq+1)}{2}  +\mu - \tau \notag \\
 &=&
\frac{s^2 + s - 2s \eq}{2}  + s  -     \#\{i : p_i = 1 \}   -    \#\{i : q_i = 1 \}   +    \eq    +   \delta_{s=\q=2}   -   \delta_{\q = 3} \notag \\
 &=&
\frac{s^2 + 3s}{2}   - (s-1)\eq +   \delta_{s=\q=2} -  \#\{i : p_i = 1 \}   -    \#\{i : q_i = 1 \}         -   \delta_{\q = 3} \notag \\
 &=&
\frac{s^2 + 3s}{2}     -  \#\{i : p_i = 1 \}   -    \#\{i : q_i = 1 \}       -   \delta_{\q = 3}.
\label{eq:gamma-exact}
\end{eqnarray}
The last equality holds because both $(s-1)\eq$ and  $\delta_{s=\q=2}$ must be $0$ for $s \neq 2$, while on the other hand $\delta_{s=\q=2} = \eq$ for $s = 2$. Thus, since all terms of (\ref{eq:gamma-exact}) after $(s^2+3s)/2$ are non-positive, $(-1)^{\p} \gamma_{K_a}$ is always bounded above by $(s^2+3s)/2$. This upper bound will be 
attained whenever $s>1$ and $p_i, q_i >1$ for all $i$, for then all correction terms must be $0$.

For the lower bound on $(-1)^{\p} \gamma_{K_a}$, observe that each of $ \#\{i : p_i = 1 \}$ and $\#\{i : q_i = 1 \} $ is at most $s$. In fact, each of these quantities is $0$ when $s=1$ (by hypothesis). Furthermore, $\#\{i : q_i = 1 \}   +   \delta_{\q = 3} \leq s$ when $s=2$. Thus
$$
{\setlength{\jot}{1.3ex}
\begin{array}{r c c c c c l}
(-1)^{\p} \gamma_{K_a}  &\geq &  \dfrac{s^2 + 3s}{2}    -1 &\geq &  1 & \qquad \mbox{when} & s = 1, \\
(-1)^{\p} \gamma_{K_a}  &\geq &  \dfrac{s^2 + 3s}{2}   - 2s &\geq &  1 & \qquad \mbox{when} & s = 2,\\
(-1)^{\p} \gamma_{K_a}  &\geq &  \dfrac{s^2 + 3s}{2}   - 2s -1 &\geq &  2 & \qquad \mbox{when} & s \geq 3.
\end{array}}
$$

\vspace{-3.5ex}
\end{proof}

Let $\ZZ [t, \ t^{-1}]$ denote the ring of Laurent polynomials with
integer coefficients and $G(2, t)$ the group of $2\times 2$ matrices
with entries in $\ZZ [t, t^{-1}]$.  The Burau representation $\psi\co
B_3 \to G(2,t)$ is defined by
\begin{equation}
  \psi(\sigma_1^{-1}) = \left[\begin{array}{cc}
			-t&1\\0&1\end{array}\right], \quad 
\psi(\sigma_2^{-1}) = \left[\begin{array}{cc}
		1&0\\t&-t\end{array}\right]. 
\label{eqn:burau} %%\eqno(5)
\end{equation}
See \cite{birman:book, jones} for more details.\protect\footnote{Our 
definition of the braid generator $\sigma_i$, depicted 
in Figure \ref{fig:3braid-generators}, corresponds to Jones' definition of 
$\sigma_i^{-1}$ \cite{jones}. The literature contains many examples of both conventions: 
compare \cite{birman:3braid-jones, burde-zieschang:knots} to \cite{kanenobu, 
murasugi}. Replacing $\sigma_i$ with $\sigma_i^{-1}$ produces the mirror 
image of a link, and affects the Jones polynomial by replacing $t$ with $t^{-1}$.}

For a braid $b\in B_3$, let $\hat {b}$ denote the closure of $b$ and
let $e_b$ denote the exponent of $b$. As calculated in \cite{jones}
(the formula is also given and used in the papers
\cite{birman:3braid-jones} and \cite{kanenobu} where properties of the
Jones polynomial of 3--braids are discussed), the Jones polynomial of
${\hat b}$ is given by
\begin{equation}
  J_{{\hat b}}(t)= (-\sqrt t)^{-e_b}\cdot (t+t^{-1}+ {\rm trace}(\psi
  (b))).
\label{eqn:jones-burau} %%\eqno(6)
\end{equation}

\begin{lemma}
Suppose that a link $K_b$ is the closure of a generic 3--braid
$$b= C^k \sigma_1^{p_1} \sigma_2^{-q_1} \cdots \sigma_1^{p_s}
\sigma_2^{-q_s},$$
with $p_i, q_i>0$.  Let $K_a$ denote the alternating link represented
by the closure of the alternating braid $a:=\sigma_1^{p_1}
\sigma_2^{-q_1} \cdots \sigma_1^{p_s} \sigma_2^{-q_s}$.  If
$J_{K_a}(t)$ and $J_{K_b}(t)$ denote the Jones polynomials of $K_a$
and $K_b$ respectively, then
$$ J_{{K_b}}(t)= t^{-6k}\, J_{{K_a}}(t)+ (-\sqrt t)^{-e_a}\,
\left(t+t^{-1}\right) \, \left(t^{-3k}-t^{-6k}\right),$$
\noindent where $e_a$ is the braid exponent of $a$.
\label{lemma:generic}
\end{lemma}

\begin{proof} 
An easy calculation, using equation (\ref{eqn:burau}), will show that
$$\psi(C)=\psi((\sigma_1 \sigma_2 \sigma_1)^{2})=
\left[\begin{array}{cc} t^{-3}&0\\0&t^{-3}\end{array}\right],$$
and thus
\begin{equation}
 {\rm trace}(\psi (b))=t^{-3k}\,  {\rm trace}(\psi (a)).
\label{eqn:trace} %%\eqno(7)
\end{equation}
The braid exponents $e_b$ and $e_a$ satisfy
$e_b=6k + e_a$.
Thus, by equations (\ref{eqn:jones-burau}) and (\ref{eqn:trace}),
\begin{equation}
  J_{{K_b}}(t)= t^{-3k}(-\sqrt t)^{-e_a}\, \left(t+t^{-1}+ t^{-3k}\,
  {\rm trace}(\psi (a))\right),
\label{eqn:jones-Kb} %%\eqno(8)
\end{equation}
and
\begin{equation}
  J_{{K_a}}(t)= (-\sqrt t)^{-e_a}\, \left(t+t^{-1}+ {\rm trace}(\psi
  (a))\right). 
\label{eqn:jones-Ka} %%\eqno(9)
\end{equation}

By eliminating ${\rm trace}(\psi (a))$ from equations
(\ref{eqn:jones-Kb}) and (\ref{eqn:jones-Ka}) we obtain
\begin{equation}
 J_{{K_b}}(t)= t^{-6k}\, J_{{K_a}}(t)+ (-\sqrt t)^{-e_a}\,
\left(t+t^{-1}\right) \, \left(t^{-3k}-t^{-6k}\right),
\label{eqn:lemma-generic} % \eqno(10)
\end{equation}
as desired.
\end{proof}
\smallskip

We are now ready to estimate certain outer coefficients of the Jones polynomial for
any generic closed 3--braid.

\begin{define}\label{def:zeta}
Let $K$ be a link in $S^3$. From the Jones polynomial $J_K(t)$, we define the following quantities. Let
$$\zeta_K = \left\{ \begin{array}{r l}
\beta_K, & \mbox{if } \abs{\alpha_K} = 1 \\
0, & \mbox{otherwise}
\end{array} \right.
\qquad \mbox{and} \qquad
\zeta'_K = \left\{ \begin{array}{r l}
\beta'_K, & \mbox{if } \abs{\alpha'_K} = 1 \\
0, & \mbox{otherwise.}
\end{array} \right.$$
% Let $\zeta = \max \left\{ \abs{\zeta_K}, \abs{\zeta'_K}  \right\} $.
 Note that by Lemma \ref{lemma:adequate-jones}, adequate links will satisfy $\zeta_K = \beta_K$ and 
$\zeta'_K = \beta'_K$. This definition gives a way to generalize approximately the same quantity.
\end{define}

\begin{prop}\label{prop:3-braids}
Let $K_b$ be the closure of a generic 3--braid
$b= C^k \sigma_1^{p_1} \sigma_2^{-q_1} \cdots \sigma_1^{p_s}
\sigma_2^{-q_s},$
with $p_i, q_i>0$.  Define $\zeta_{K_b}$ and $\zeta'_{K_b}$ as in Definition \ref{def:zeta}. Then
$$s-1 \leq  \max \left\{ \abs{\zeta_{K_b}}, \abs{\zeta'_{K_b}} \right\} \leq s+1.$$
\end{prop}

\begin{proof}
For the majority of this argument, we will work under the hypothesis that $\p >1$ and $\q >1$. 
At the end of the proof, we will consider the (simpler) case when $\p = 1$ or $\q = 1$.

So: assume that $\p >1$ and $\q >1$, and let $a = \sigma_1^{p_1} \sigma_2^{-q_1} \cdots \sigma_1^{p_s}$ be the alternating part of $b$. Consider $J_{K_b}(t)$, as expressed as a sum of two terms in equation
(\ref{eqn:lemma-generic}).  By Lemma \ref{lemma:alternating}, the
first term, $t^{-6k}\,J_{K_a}(t)$, is
\begin{eqnarray}
  t^{-6k} \, J_{{K_a}}(t) & = 
   (-1)^{{\bf p}} \ t^{(3\q-\p)/2-6k} +
   (-1)^{\p +1}(s- \eq) \, t^{(3\q-\p)/2-6k-1} + \cdots
		\label{eqn:jones-1st} %%\eqno(11)
 \\
  &  + (-1)^{\q +1}(s - \ep) \, t^{(\q-3\p)/2-6k+1} +
   (-1)^{{\bf q}} \ t^{(\q-3\p)/2-6k}. \notag
\end{eqnarray}
Meanwhile, the second term on the right hand
side of equation (\ref{eqn:lemma-generic}) expands out to
\begin{equation}
(-1)^{e_a} \left( t^{-e_a/2-3k +1} +
t^{-e_a/2-3k-1}- t^{-e_a/2-6k+1}-
t^{-e_a/2-6k-1}\right).
\label{eqn:jones-2nd} %%\eqno(12).
\end{equation}

If $k=0$, then the expression in (\ref{eqn:jones-2nd}) vanishes, and
the link $K_b$ is alternating. Thus  $\zeta_{K_b} = \beta_{K_b}$ and 
$\zeta'_{K_b} = \beta'_{K_b}$, and the desired result is true by
Lemma \ref{lemma:alternating}.

Next, suppose that $k \neq 0$. We claim that no generality is lost by assuming that $k > 0$. Otherwise, if $k < 0$, the mirror image $K_d$ of the link $K_b$ can be represented by the braid word $d = C^{-k} \sigma_1^{q_s} \sigma_2^{-p_1} \sigma_1^{q_1} \cdots \sigma_2^{-p_s} ,$ so the power of $C$ will now be positive. The Jones polynomial $J_{K_d}$ can be obtained from $J_{K_b}$ by interchanging $t$ and $t^{-1}$, so $\zeta_{K_b} = \zeta'_{K_d}$ and $\zeta'_{K_b} = \zeta_{K_d}$, with the maximum of the two values unaffected. Thus we may assume $k > 0$.

If $k > 0$, the monomials of (\ref{eqn:jones-2nd}) are listed in order of decreasing powers of $t$, and each monomial has a  coefficient of $\pm 1$. Furthermore, we claim that the degree of any term in (\ref{eqn:jones-2nd}) is strictly higher than the degree of the lowest term in (\ref{eqn:jones-1st}). This is because the last monomial of (\ref{eqn:jones-2nd}) has degree $-e_a/2-6k-1 = (\q-\p)/2 - 6k -1$, and 
\begin{equation}\label{eq:powers}
\p \geq 2 \qquad \mbox{is equivalent to} \qquad   \frac{\q-\p}{2} - 6k -1 \: \geq \: \frac{\q-3\p}{2} - 6k +1.
\end{equation}
Thus the lowest--degree term of $J_{K_b}(t)$ is $ (-1)^{\q} \ t^{(\q-3\p)/2-6k}$, and $\alpha'_{K_b} = (-1)^{\q}$.

From equation (\ref{eq:powers}), we can also conclude that the last monomial of (\ref{eqn:jones-2nd}) only affects the next-to-last monomial of (\ref{eqn:jones-1st})  if $\p=2$. When $\p$ is even, the signs of these two monomials are $(-1)^{\q +1}$ and $(-1)^{e_a + 1} = (-1)^{\p - \q + 1} = (-1)^{\q + 1}$: they have the same sign. Thus when $\p = 2$, the last monomial of (\ref{eqn:jones-2nd}) will contribute $1$ to $|\beta'_{K_b}|$; when $\p > 2$, no monomial of (\ref{eqn:jones-1st}) affects $|\beta'_{K_b}|$ at all. We conclude that when $k >0$,
$$\zeta'_{K_b} \: = \: \beta'_{K_a} + (-1)^{\q + 1}\ep \: = \: (-1)^{\q + 1} (s - \ep + \ep) \: = \: (-1)^{\q + 1} s.$$

Next, consider how the top two terms of (\ref{eqn:jones-1st}) might interact with the monomials of (\ref{eqn:jones-2nd}). If the top degree of (\ref{eqn:jones-2nd}) is lower than the top degree of (\ref{eqn:jones-1st}), we will have $\alpha_{K_b} = (-1)^{\p}$ and $\zeta_{K_b} = \beta_{K_b}$ will be off by at most $1$ from $\beta_{K_a}$. In particular, $s-2 \leq \abs{\zeta_{K_b}} \leq s+1$. If the top degree of (\ref{eqn:jones-2nd}) is higher than the top degree of (\ref{eqn:jones-1st}) by $2$ or more, we will have $\alpha_{K_b} = (-1)^{e_a}$ and $\zeta_{K_b} = \beta_{K_b} = 0$. If the top degree of (\ref{eqn:jones-2nd}) is higher by exactly $1$, then we will have $\alpha_{K_b} = (-1)^{e_a}$ and $\zeta_{K_b} = \beta_{K_b} = (-1)^\p$. If the top degree of (\ref{eqn:jones-2nd}) is exactly equal to the top degree of (\ref{eqn:jones-1st}), then the two monomials either add or cancel. If the top monomials add, $\alpha_{K_b} = (-1)^\p \cdot 2$ and $\zeta_{K_b} = 0$ by Definiton \ref{def:zeta}. Thus, in all cases when the top monomials of  (\ref{eqn:jones-1st}) and (\ref{eqn:jones-2nd}) do not cancel, we have $0 \leq \abs{\zeta_{K_b}} \leq s+1$.

If the top monomials of  (\ref{eqn:jones-1st}) and (\ref{eqn:jones-2nd}) cancel (which can occur: see Proposition \ref{prop:cancelation}), we have
$$\frac{3\q-\p}{2}-6k = \frac{-e_a}{2}-3k +1, \quad \mbox{which simplifies to} \quad
\q = 3k + 1,$$
since $e_a = \p - \q$. In particular, it follows that $\q \geq 4$. Then
 $\alpha_{K_b} = \beta_{K_a} = (-1)^{\p + 1}s$. So, if $s \geq 2$, we have $\abs{\alpha_{K_b}} \geq 2$, and again $\zeta_{K_b} = 0$ by Definiton \ref{def:zeta}. Thus when $s \geq 2$, we get $0 \leq \abs{\zeta_{K_b}} \leq s+1$ in all cases. Since $|\zeta'_{K_b}|$ is always equal to $s$ if $k>0$, we conclude that
$$s-1 \leq  \max \left\{ \abs{\zeta_{K_b}}, \abs{\zeta'_{K_b}} \right\} \leq s+1
\qquad \mbox{whenever} \qquad s \geq 2.$$
Note that the lower bound of $s-1$ is only achieved when $k = 0$ and $s = \p = \q = 2$, i.e., when $K_b$ is the figure--$8$ knot.

If $s = 1$ and the top monomials of  (\ref{eqn:jones-1st}) and (\ref{eqn:jones-2nd}) cancel, then $\alpha_{K_b} = \beta_{K_a} = (-1)^{\p + 1}$ and $\zeta_{K_b} = \beta_{K_b} = \gamma_{K_a} + (-1)^{\p+1}$,
where $\gamma_{K_a}$ is the third coefficient of $J_{K_a}(t)$ and $(-1)^{\p+1}$ comes from the second monomial of (\ref{eqn:jones-2nd}). By Lemma \ref{lemma:alternating}(c),
$$\abs{\zeta_{K_b}} \: = \: \abs{\gamma_{K_a}} - 1 \: \leq \: \frac{s(s+3)}{2} - 1\: = \: 1.$$
Since we get $0 \leq \abs{\zeta_{K_b}} \leq s+1$ whenever the top monomials do not cancel, and since $|\zeta'_{K_b}|$ is always equal to $s$ if $k>0$, we conclude that
$$s-1 \leq  \max \left\{ \abs{\zeta_{K_b}}, \abs{\zeta'_{K_b}} \right\} \leq s+1
\qquad \mbox{whenever} \qquad s = 1 \quad \mbox{and} \quad \p, \q \geq 2.$$

To complete the proof, we need to consider the case when $\p = 1$ or $\q = 1$, hence $s=1$. Without loss of generality, we may assume that $\p \geq \q = 1$ (otherwise, take the mirror image of $K_b$, as above). Then $K_a$ is the $(2,\p)$ torus link. By direct computation (for example, using \cite[Theorem 14.13]{lickbook} or \cite[Proposition 2.1 and Example 3]{ck:cyclotomic}), it follows that
$$J_{K_a}(t) = (-1)^{\p+1} \left( t^{-(\p-1)/2} + t^{-(\p+3)/2} - t^{-(\p+5)/2}  +  t^{-(\p+7)/2} - \ldots + (-1)^\p t^{-(3\p-1)/2} \right).$$
In particular, the second coefficient of $J_{K_a}(t)$ is $0$ and every other coefficient is $\pm 1$. If $k \neq 0$, we may still use Lemma \ref{lemma:generic} to compute the Jones polynomial $J_{K_b}(t)$. Thus, after multiplying $J_{K_a}(t)$ by $t^{-6k}$ and adding in the four monomials of (\ref{eqn:jones-2nd}), it will follow that every coefficient of $J_{K_b}(t)$ is $0$, $\pm 1$, or $\pm 2$. Thus
$$0 \leq  \max \left\{ \abs{\zeta_{K_b}}, \abs{\zeta'_{K_b}} \right\} \leq 2
\qquad \mbox{whenever} \qquad  \p = 1 \quad \mbox{or} \quad \q = 1,$$
which is exactly what the proposition requires for $s=1$.
\end{proof}

\subsection{Connections to volume}\label{subsec:volume-jones}

Proposition \ref{prop:3-braids} and Theorem \ref{thm:3braid-volume}
immediately imply the following.

\begin{theorem}
Let $K$ be a hyperbolic closed $3$--braid. 
From the Jones polynomial $J_K(t)$, define ${\zeta_K}, {\zeta'_K}$ as in Definition \ref{def:zeta}.
Let $\zeta = \max \left\{ \abs{\zeta_K}, \abs{\zeta'_K} \right\}$. Then
$$4v_3 \cdot \zeta - 281 \: < \: \vol(S^3 \setminus K) \: < \:
4v_8 \, ( \zeta +1 ).$$
Furthermore, the multiplicative constants in both the upper and lower bounds are sharp.
\label{thm:relations}
\qed
\end{theorem}

% 	In fact, by considering the Jones polynomials of non-generic 3--braids, one can show 
%	that Theorem \ref{thm:relations} 
%	also holds for a non-hyperbolic link $K$, with $\vol(S^3 \setminus K) = 0$.

In contrast with Proposition \ref{prop:3-braids}, there exist
$3$--braids for which the second coefficient of the Jones polynomial
is quite different from $s$.

\begin{prop}\label{prop:cancelation}
For every $s>1$ there is a knot $K = K_s$, represented by the 3--braid word
$$b= C^k \sigma_1^{p_1} \sigma_2^{-q_1} \cdots \sigma_1^{p_s}
\sigma_2^{-q_s},\ \ \mbox{with} \ \ p_i, q_i>0,$$
such that the second and next-to-last coefficients of the Jones
polynomial $J_K(t)$ satisfy
$$\beta_K \: = \: \frac{s(s+3)}{-2} + 1 ,  \qquad \beta'_K = s.$$
\end{prop}

The head and tail of the Jones polynomial for several values of $s$ is computed in Table \ref{table:jones}.

\begin{table}
$$
{\setlength{\jot}{1.5ex}
\begin{array}{|c | r c l | }
\hline
s & \multicolumn{3}{c|}{\mbox{Jones polynomial } J_{K_s} } \\
\hline
2 &  2 \,t ^{-8} -4 \,t ^{-9} +9 \,t ^{-10} -14 \,t ^{-11} + & \!\! \ldots \!\! & +\, 9 \,t ^{-22} -5 \,t ^{-23} +2 \,t ^{-24} -t ^{-25} \\
3 & 3 \,t ^{-12} -8 \,t ^{-13} +20 \,t ^{-14} -39 \,t ^{-15} + & \!\! \ldots \!\!  & +\,20 \,t ^{-34} -9 \,t ^{-35} +3 \,t ^{-36} -t ^{-37}  \\
4 & 4 \,t ^{-16} -13 \,t ^{-17} +37 \,t ^{-18} -85 \,t ^{-19} + & \!\! \ldots \!\!  & +\,37 \,t ^{-46} -14 \,t ^{-47} +4 \,t ^{-48} -t ^{-49}  \\
5 & 5 \,t ^{-20} -19 \,t ^{-21} +61 \,t ^{-22} -160 \,t ^{-23} + & \!\! \ldots \!\!  & +\,61 \,t ^{-58} -20 \,t ^{-59} +5 \,t ^{-60} -t ^{-61} \\
6 & 6 \,t ^{-24} -26 \,t ^{-25} +93 \,t ^{-26} -273 \,t ^{-27} + & \!\! \ldots \!\!  & +\,93 \,t ^{-70} -27 \,t ^{-71} +6 \,t ^{-72} -t ^{-73} \\
7 & 7 \,t ^{-28} -34 \,t ^{-29} +134 \,t ^{-30} -434 \,t ^{-31} + & \!\! \ldots \!\!  & +\,134 \,t ^{-82} -35 \,t ^{-83} +7 \,t ^{-84} -t ^{-85} \\
8 & 8 \,t ^{-32} -43 \,t ^{-33} +185 \,t ^{-34} -654 \,t ^{-35} + & \!\! \ldots \!\!  & +\,185 \,t ^{-94} -44 \,t ^{-95} +8 \,t ^{-96} -t ^{-97} \\
9 & 9 \,t ^{-36} -53 \,t ^{-37} +247 \,t ^{-38} -945 \,t ^{-39} + & \!\! \ldots \!\!  & +\,247 \,t ^{-106} -54 \,t ^{-107} +9 \,t ^{-108} -t ^{-109} \\
10 & 10 \,t ^{-40} -64 \,t ^{-41} +321 \,t ^{-42} -1320 \,t ^{-43} + & \!\! \ldots \!\!  & +321 \,t ^{-118} -65 \,t ^{-119} +10 \,t ^{-120} -t ^{-121} \\\hline
\end{array}
}
$$
\caption{The head and tail of the Jones polynomial for the first few values of $s$ in Proposition \ref{prop:cancelation}. These values were computed independently by two pieces of software: {\tt jpclosed.c} by Nathan Broaddus 
and the Mathematica package {\tt KnotTheory`} by Dror Bar-Natan.
\label{table:jones}}
\end{table}

\begin{proof} 
Fix $s>1$, and let
$p_1 = \ldots = p_{s-1} = 2, \quad p_s = 3, \quad
q_1 = \ldots = q_{s-1} = 6, \quad q_s = 7.$

Notice that $\p = \sum p_i=2s+1$ and $\q = \sum
q_i=6s+1$; thus, in particular, they are both odd.  Let $k=2s$
and let $K$ be the closure of the generic $3$--braid
$$b= C^k \sigma_1^{p_1} \sigma_2^{-q_1} \cdots \sigma_1^{p_s}
\sigma_2^{-q_s}.$$
Since $C^k$ is a pure braid, and each of $\sigma_1^{p_1}, \ldots,
\sigma_2^{-q_{s-1}}$ is also a pure braid (because their powers are
even), $K$ will have the same number of components as the closure of
$\sigma_1^{3} \sigma_2^{-7}$.  Hence $K$ is a knot.
 
The Jones polynomial $J_K(t)$ will be the sum of equations
(\ref{eqn:jones-1st}) and (\ref{eqn:jones-2nd}). Note that the leading
terms of these equations are
$$(-1)^{\p} \ t^{(3\q-\p)/2-6k} \: = \: -t^{(16s+2)/2-12s} =\: - t^{-4s+1}, $$
and
$$(-1)^{\p-\q} \  t^{(\q-\p)/2-3k+1} \: = \: t^{(4s)/2-6s+1} \: = \: t^{-4s+1},$$
which will cancel.  Thus, continuing along equations
(\ref{eqn:jones-1st}) and (\ref{eqn:jones-2nd}), one can see that the
leading monomial of $J_{K_b}(t)$ will be $(-1)^{\p+1} s\, t^{-4s}$,
and the second term will be
$$(\gamma_{K_a}) t^{-4s-1} + (-1)^{(\p-\q)}\ t^{-4s-1} \: =
\: \left( \gamma_{K_a} + 1 \right) t^{-4s-1}.$$
Since $s>1$ and $p_i, q_i > 1$ for all $i$, Lemma \ref{lemma:alternating}(c) implies that
$$ \gamma_{K_a} \:=\: (-1)^{\p}\, \left( \frac{s(s+3)}{2} \right) 
\:=\: \frac{s(s+3)}{-2},
\qquad \mbox{hence} \qquad
\beta_K \: = \:  \frac{s(s+3)}{-2}+1.$$

Now we consider the term $\beta'_K$.
The two lowest--degree terms of (\ref{eqn:jones-1st})  are
$$(-1)^{\q} \ t^{(\q-3\p)/2 -6k} \: = \: - t^{(-2)/2 - 12s }  \: = \: -t^{-12s-1} \qquad
\mbox{and} \qquad
(-1)^{\q +1}s \, t^{-12s} \: = \: s \, t^{-12s}.
$$
The remaining terms of  (\ref{eqn:jones-2nd}) are
$$(-1)^{\p-\q+1} t^{(\q-\p)/2-6k\pm 1} \: =\: 
t^{(4s)/2 - 12s \pm 1} \: = \:  -t^{-10s \pm 1} .$$
Since  $-12s < -10s-1$ for all positive $s$, no term of (\ref{eqn:jones-2nd}) can affect either
of the two lowest--degree terms of $J_{K_b}(t)$. Therefore,
$$\alpha'_K = (-1)^{\q} = -1 \qquad \mbox{and} \qquad
\beta'_K = (-1)^{\q+1} s = s.$$

\vspace{-3ex}
\end{proof}

\begin{theorem} \label{thm:counterexample}
There does not exist a function $f(\cdot, \cdot)$ of two variables,
together with positive constants $C_1, \ldots, C_4$ such that
$$C_1 f(\beta_K, \beta'_K) - C_2 \: < \: \vol(S^3\setminus K) \:<\: C_3
f(\beta_K, \beta'_K) + C_4$$
for every hyperbolic knot $K$.
In other words, the second and
next-to-last coefficients of the Jones polynomial do not coarsely
predict the volume of a knot.
\end{theorem}

\begin{proof}
Suppose, for a contradiction, that such a function $f(\cdot, \cdot)$
does exist. Then it will follow that for every pair of knots $K$ and
$L$, such that $\beta_K = \beta_L$ and $\beta'_K = \beta'_L$,
\begin{equation}\label{eqn:volume-ratio}
\vol(S^3 \setminus L) \: < \: \frac{C_3}{C_1} \, \vol(S^3 \setminus K)
+ \frac{C_2 C_3}{C_1} + C_4 .
\end{equation}
On the other hand, we shall construct an infinite sequence of such
pairs $K_s$ and $L_s$, such that the volumes
$\vol(S^3 \setminus
L_s)$,  $\vol(S^3 \setminus K_s)$, as well as the ratio 
$\vol(S^3 \setminus
L_s) / \vol(S^3 \setminus K_s)$ go to infinity as $s$ goes to infinity. This will contradict
the existence of $f$.

The sequence $K_s$ is the one given by Proposition
\ref{prop:cancelation}, where we require that $s \equiv 3 \mod 4$ and
$s >10$.  By Proposition \ref{prop:cancelation}, these knots
satisfy
$$\beta_{K_s} = \frac{s(s+3)}{-2} + 1 , \qquad \mbox{and} \qquad
\beta'_{K_s} = s.$$

The sequence $L_s$ will consist of $(3, \ldots, 3, -3, \ldots, -3)$
pretzel knots, where there are
$$n = \left( \frac{s(s+3)}{2} - 2 \right) \quad \mbox{positive 3's} \quad \mbox{and}
\quad m = \left( s-1 \right) \quad \mbox{negative 3's}. $$
When $s \equiv 3 \mod 4$, it follows that $n$ is odd, $m$ is even, and
thus $L_s$ is indeed a knot.

It is easy to check that the standard pretzel diagram $D_s$ of one
these knots is adequate. Computing the graphs $\GA(D_s)$, $\GB(D_s)$
gives
$$v_A = n+2m, \qquad e'_A = n+3m, \qquad v_B = 2n+m, \qquad e'_B =
3n+m.$$
Also, we may compute the writhe of $D_s$ to be $w(D) = 3m-3n$.  Now,
Lemma \ref{lemma:adequate-jones} gives
$$\beta_{L_s}  \:=\: (-1)^{v_B + w(D)} (e'_B - v_B + 1) 
\:=\: (-1)^{-n + 4m} (n+1)
\:=\: (-1)^n(n+1) \:=\: \frac{s(s+3)}{-2} + 1,$$
since $n$ is odd. Similarly,
$$\qquad \beta'_{L_s} \:=\: (-1)^{v_A + w(D)} (e'_A - v_A + 1) \:=\:
(-1)^m(m+1) \:=\: s.$$
Therefore, $\beta_{K_s} = \beta_{L_s}$ and $\beta'_{K_s} =
\beta'_{L_s}$ for all $s \equiv 3 \mod 4$.

Now, consider the volumes of these links. By Theorem
\ref{thm:3braid-volume},
$$4v_3\,s -277 < \vol(S^3 \setminus K_s) < 4 v_8 \, s.$$
On the other hand, the pretzel diagram $D_s$ contains $(n+m) > s^2/2$
twist regions and is built by joining together $(n+m) > s^2/2$ rational
tangles. Thus, by \cite[Theorem 1.5]{fkp:conway},
$$\vol(S^3 \setminus L_s) \: \geq \: \frac{v_8}{2} \left(1 - \left(
\frac{8\pi}{11.524 + s^2 \sqrt[4]{2} /2 } \right)^2 \right)^{\!3/2}
\!\! \left(\frac{s^2}{2} -3\right) \: > \: \frac{v_8 \, s^2}{8} \qquad
\mbox{for } s >10.
$$
Hence $\vol(S^3 \setminus L_s) / \vol(S^3 \setminus K_s) \geq s/32$.
Since both volumes are eventually large enough that the additive
constants in equation (\ref{eqn:volume-ratio}) do not matter, this
contradicts equation (\ref{eqn:volume-ratio}).  Thus the function $f$
cannot exist.
\end{proof}

\bibliographystyle{hamsplain}
\bibliography{biblio.bib}

\providecommand{\bysame}{\leavevmode\hbox to3em{\hrulefill}\thinspace}
\providecommand{\href}[2]{#2}
\begin{thebibliography}{10}

\bibitem{adams:waist2}
Colin Adams, \emph{Waist size for cusps in hyperbolic 3-manifolds {II}},
  Preprint.

\bibitem{adams:2-gen}
\bysame, \emph{Hyperbolic {$3$}-manifolds with two generators}, Comm. Anal.
  Geom. \textbf{4} (1996), no.~1-2, 181--206.

\bibitem{adams-students:cusp-area}
Colin Adams, A.~Colestock, J.~Fowler, W.~Gillam, and E.~Katerman, \emph{Cusp
  size bounds from singular surfaces in hyperbolic 3-manifolds}, Trans. Amer.
  Math. Soc. \textbf{358} (2006), no.~2, 727--741 (electronic).

\bibitem{agol:bounds}
Ian Agol, \emph{Bounds on exceptional {D}ehn filling}, Geom. Topol. \textbf{4}
  (2000), 431--449 (electronic).

\bibitem{alr:alternating}
Iain~R. Aitchison, E.~Lumsden, and J.~Hyam Rubinstein, \emph{Cusp structures of
  alternating links}, Invent. Math. \textbf{109} (1992), no.~3, 473--494.

\bibitem{akiyoshi}
Hirotaka Akiyoshi, \emph{On the {F}ord domains of once-punctured torus groups},
  S\=urikaisekikenky\=usho K\=oky\=uroku (1999), no.~1104, 109--121, Hyperbolic
  spaces and related topics (Japanese) (Kyoto, 1998).

\bibitem{aswy}
Hirotaka Akiyoshi, Makoto Sakuma, Masaaki Wada, and Yasushi Yamashita,
  \emph{Ford domains of punctured torus groups and two-bridge knot groups},
  S\=urikaisekikenky\=usho K\=oky\=uroku (2000), no.~1163, 67--77, Hyperbolic
  spaces and related topics, II (Japanese) (Kyoto, 1999).

\bibitem{aswy-book}
\bysame, \emph{Punctured torus groups and 2-bridge knot groups. {I}}, Lecture
  Notes in Mathematics, vol. 1909, Springer, Berlin, 2007.

\bibitem{birman:book}
Joan~S. Birman, \emph{Braids, links, and mapping class groups}, Princeton
  University Press, Princeton, N.J., 1974, Annals of Mathematics Studies, No.
  82.

\bibitem{birman:3braid-jones}
\bysame, \emph{On the {J}ones polynomial of closed {$3$}-braids}, Invent. Math.
  \textbf{81} (1985), no.~2, 287--294.

\bibitem{birman-menasco:3braids}
Joan~S. Birman and William~W. Menasco, \emph{Studying links via closed braids.
  {III}. {C}lassifying links which are closed {$3$}-braids}, Pacific J. Math.
  \textbf{161} (1993), no.~1, 25--113.

\bibitem{birman-menasco:toroidal}
\bysame, \emph{Special positions for essential tori in link complements},
  Topology \textbf{33} (1994), no.~3, 525--556.

\bibitem{bleiler-hodgson}
Steven~A. Bleiler and Craig~D. Hodgson, \emph{Spherical space forms and {D}ehn
  surgery}, Knots 90 (Osaka, 1990), de Gruyter, Berlin, 1992, pp.~425--433.

\bibitem{boroczky}
K\'{a}roly B{\"o}r{\"o}czky, \emph{Packing of spheres in spaces of constant
  curvature}, Acta Math. Acad. Sci. Hungar. \textbf{32} (1978), no.~3-4,
  243--261.

\bibitem{burde-murasugi}
Gerhard Burde and Kunio Murasugi, \emph{Links and {S}eifert fiber spaces}, Duke
  Math. J. \textbf{37} (1970), 89--93.

\bibitem{burde-zieschang:knots}
Gerhard Burde and Heiner Zieschang, \emph{Knots}, second ed., de Gruyter
  Studies in Mathematics, vol.~5, Walter de Gruyter \& Co., Berlin, 2003.

\bibitem{ck:cyclotomic}
Abhijit Champanerkar and Ilya Kofman, \emph{On links with cyclotomic {J}ones
  polynomials}, Algebr. Geom. Topol. \textbf{6} (2006), 1655--1668
  (electronic).

\bibitem{ckp:simplest-knots}
Abhijit Champanerkar, Ilya Kofman, and Eric Patterson, \emph{The next simplest
  hyperbolic knots}, J. Knot Theory Ramifications \textbf{13} (2004), no.~7,
  965--987.

\bibitem{dasbach-futer...}
Oliver~T. Dasbach, David Futer, Efstratia Kalfagianni, Xiao-Song Lin, and
  Neal~W. Stoltzfus, \emph{The {J}ones polynomial and graphs on surfaces},
  Journal of Combinatorial Theory Ser. B \textbf{98} (2008), no.~2, 384--399.

\bibitem{dasbach-futer2...}
\bysame, \emph{Alternating sum formulae for the determinant and other link
  invariants}, J. Knot Theory Ramifications (to appear),
  \mbox{arXiv:math/0611025}.

\bibitem{dasbach-lin:head-tail}
Oliver~T. Dasbach and Xiao-Song Lin, \emph{On the head and the tail of the
  colored {J}ones polynomial}, Compos. Math. \textbf{142} (2006), no.~5,
  1332--1342.

\bibitem{dasbach-lin:volumeish}
\bysame, \emph{{A volume-ish theorem for the Jones polynomial of alternating
  knots}}, Pacific J. Math. \textbf{231} (2007), no.~2, 279--291.

\bibitem{epstein-penner}
David B.~A. Epstein and Robert~C. Penner, \emph{Euclidean decompositions of
  noncompact hyperbolic manifolds}, J. Differential Geom. \textbf{27} (1988),
  no.~1, 67--80.

\bibitem{floyd-hatcher}
William Floyd and Allen Hatcher, \emph{Incompressible surfaces in
  punctured-torus bundles}, Topology Appl. \textbf{13} (1982), no.~3, 263--282.

\bibitem{fkp-07}
David Futer, Efstratia Kalfagianni, and Jessica~S. Purcell, \emph{{Dehn
  filling, volume, and the Jones polynomial}}, J. Differential Geom.
  \textbf{78} (2008), no.~3, 429--464.

\bibitem{fkp:conway}
\bysame, \emph{Symmetric links and {C}onway sums: volume and {J}ones
  polynomial}, Math. Res. Lett. \textbf{16} (2009), no.~2, 233--253.

\bibitem{fkp:coils}
\bysame, \emph{On diagrammatic bounds of knot volumes and spectral invariants},
  Geom. Dedicata (to appear), \mbox{arXiv:math/0901.0119}.

\bibitem{futer-purcell}
David Futer and Jessica~S. Purcell, \emph{Links with no exceptional surgeries},
  Comment. Math. Helv. \textbf{82} (2007), no.~3, 629--664.

\bibitem{gmm:smallest-cusped}
David Gabai, Robert Meyerhoff, and Peter Milley, \emph{Minimum volume cusped
  hyperbolic three-manifolds}, J. Amer. Math. Soc. \textbf{22} (2009), no.~4,
  1157--1215.

\bibitem{gueritaud:thesis}
Fran\c{c}ois Gu{\'e}ritaud, \emph{G{\'e}om{\'e}trie hyperbolique effective et
  triangulations id{\'e}ales cononiques en dimension 3.}, Ph.D. thesis,
  L'Universit{\'e} de Paris, 2006.

\bibitem{gf:two-bridge}
Fran\c{c}ois Gu\'eritaud and David Futer~(appendix), \emph{On canonical
  triangulations of once-punctured torus bundles and two-bridge link
  complements}, Geom. Topol. \textbf{10} (2006), 1239--1284.

\bibitem{hk:univ}
Craig~D. Hodgson and Steven~P. Kerckhoff, \emph{Universal bounds for hyperbolic
  {D}ehn surgery}, Ann. of Math. (2) \textbf{162} (2005), no.~1, 367--421.

\bibitem{jones}
Vaughan F.~R. Jones, \emph{Hecke algebra representations of braid groups and
  link polynomials}, Ann. of Math. (2) \textbf{126} (1987), no.~2, 335--388.

\bibitem{jorgensen}
Troels J{\o}rgensen, \emph{On pairs of once-punctured tori}, Kleinian groups
  and hyperbolic 3-manifolds (Warwick, 2001), London Math. Soc. Lecture Note
  Ser., vol. 299, Cambridge Univ. Press, Cambridge, 2003, pp.~183--207.

\bibitem{kanenobu}
Taizo Kanenobu, \emph{Relations between the {J}ones and {$Q$} polynomials for
  {$2$}-bridge and {$3$}-braid links}, Math. Ann. \textbf{285} (1989), no.~1,
  115--124.

\bibitem{kashaev:volume-conj}
Rinat Kashaev, \emph{The hyperbolic volume of knots from the quantum
  dilogarithm}, Lett. Math. Phys. \textbf{39} (1997), no.~3, 269--275.

\bibitem{lackenby:word}
Marc Lackenby, \emph{Word hyperbolic {D}ehn surgery}, Invent. Math.
  \textbf{140} (2000), no.~2, 243--282.

\bibitem{lackenby:punct}
\bysame, \emph{The canonical decomposition of once-punctured torus bundles},
  Comment. Math. Helv. \textbf{78} (2003), no.~2, 363--384.

\bibitem{lickbook}
W.~B.~Raymond Lickorish, \emph{An introduction to knot theory}, Graduate Texts
  in Mathematics, vol. 175, Springer-Verlag, New York, 1997.

\bibitem{lozano-przytycki}
Mar{\'{\i}}a~Teresa Lozano and J{\'o}zef~H. Przytycki, \emph{Incompressible
  surfaces in the exterior of a closed {$3$}-braid. {I}. {S}urfaces with
  horizontal boundary components}, Math. Proc. Cambridge Philos. Soc.
  \textbf{98} (1985), no.~2, 275--299.

\bibitem{maskit}
Bernard Maskit, \emph{Kleinian groups}, Grundlehren der Mathematischen
  Wissenschaften [Fundamental Principles of Mathematical Sciences], vol. 287,
  Springer-Verlag, Berlin, 1988.

\bibitem{minsky}
Yair~N. Minsky, \emph{The classification of punctured-torus groups}, Ann. of
  Math. (2) \textbf{149} (1999), no.~2, 559--626.

\bibitem{morton:composite-braids}
Hugh~R. Morton, \emph{Closed braids which are not prime knots}, Math. Proc.
  Cambridge Philos. Soc. \textbf{86} (1979), no.~3, 421--426.

\bibitem{murakami-squared}
Hitoshi Murakami and Jun Murakami, \emph{The colored {J}ones polynomials and
  the simplicial volume of a knot}, Acta Math. \textbf{186} (2001), no.~1,
  85--104.

\bibitem{murasugi:3braids}
Kunio Murasugi, \emph{On closed {$3$}-braids}, American Mathematical Society,
  Providence, R.I., 1974, Memoirs of the American Mathmatical Society, No. 151.

\bibitem{murasugi}
\bysame, \emph{Knot theory \& its applications}, Modern Birkh\"auser Classics,
  Birkh\"auser Boston Inc., Boston, MA, 2008, Translated from the 1993 Japanese
  original by Bohdan Kurpita, Reprint of the 1996 translated edition
  [MR1391727].

\bibitem{nimershiem:flat-tori}
Barbara~E. Nimershiem, \emph{Isometry classes of flat {$2$}-tori appearing as
  cusps of hyperbolic {$3$}-manifolds are dense in the moduli space of the
  torus}, Low-dimensional topology (Knoxville, TN, 1992), Conf. Proc. Lecture
  Notes Geom. Topology, III, Int. Press, Cambridge, MA, 1994, pp.~133--142.

\bibitem{orlik:seifert-manifolds}
Peter Orlik, \emph{Seifert manifolds}, Springer-Verlag, Berlin, 1972, Lecture
  Notes in Mathematics, Vol. 291.

\bibitem{purcell:cusps}
Jessica~S. Purcell, \emph{Cusp shapes under cone deformation}, J. Differential
  Geom. \textbf{80} (2008), no.~3, 453--500.

\bibitem{sakuma-weeks}
Makoto Sakuma and Jeffrey Weeks, \emph{Examples of canonical decompositions of
  hyperbolic link complements}, Japan. J. Math. (N.S.) \textbf{21} (1995),
  no.~2, 393--439.

\bibitem{schreier}
Otto Schreier, \emph{{\"U}ber die {G}ruppen ${A}^a {B}^b=1$}, Abh.\ Math.\
  Sem.\ Univ.\ Hamburg \textbf{3} (1924), 167--169.

\bibitem{schubert:braids}
Horst Schubert, \emph{Knoten und {V}ollringe}, Acta Math. \textbf{90} (1953),
  131--286.

\bibitem{seifert:fibering}
Herbert Seifert, \emph{Topologie {D}reidimensionaler {G}efaserter {R}\"aume},
  Acta Math. \textbf{60} (1933), no.~1, 147--238.

\bibitem{stoimenow:amphichiral}
Alexander Stoimenow, \emph{Coefficients and non-triviality of the {J}ones
  polynomial}, J. Reine Angew. Math. (to appear), Theorem numbers correspond to
  the version at arXiv:math/0606255.

\bibitem{thurston:hypII}
William~P. Thurston, \emph{Hyperbolic structures on 3-manifolds, {II}: Surface
  groups and 3-manifolds which fiber over the circle},
  \mbox{arXiv:math/9801045}.

\bibitem{thurston:notes}
\bysame, \emph{The geometry and topology of three-manifolds}, Princeton Univ.
  Math. Dept. Notes, 1979.

\bibitem{thurston:survey}
\bysame, \emph{Three-dimensional manifolds, {K}leinian groups and hyperbolic
  geometry}, Bull. Amer. Math. Soc. (N.S.) \textbf{6} (1982), no.~3, 357--381.

\bibitem{weeks:snappea}
Jeffrey~R. Weeks, \emph{{Snappea}}, Available at
  http://www.geometrygames.org/SnapPea/.

\end{thebibliography}

\end{document}